\documentstyle[twoside,amssymb,texdraw]{article}

\oddsidemargin=\evensidemargin
\catcode`\@=11
\long\def\@makefntext#1{
\protect\noindent \hbox to 3.2pt {\hskip-.9pt
$^{{\eightrm\@thefnmark}}$\hfil}#1\hfill}               

\def\ps@myheadings{\let\@mkboth\@gobbletwo              
\def\@oddhead{\hbox{}
\rightmark\hfil\eightrm\thepage}
\def\@oddfoot{}\def\@evenhead{\eightrm\thepage\hfil
\leftmark\hbox{}}\def\@evenfoot{}
\def\sectionmark##1{}\def\subsectionmark##1{}}

\catcode`@=11                                
\def\ps@plain{\let\@mkboth\@gobbletwo
     \def\@oddhead{}\def\@oddfoot{\eightrm\hfil\thepage
     \hfil}\def\@evenhead{}\let\@evenfoot\@oddfoot}

\renewcommand{\thefootnote}{\fnsymbol{footnote}}

\newcounter{sectionc}\newcounter{subsectionc}\newcounter{subsubsectionc}
\renewcommand{\section}[1] {\vspace{12pt}\addtocounter{sectionc}{1}
\setcounter{subsectionc}{0}\setcounter{subsubsectionc}{0}\noindent
        {\tenbf\thesectionc. #1}\par\vspace{5pt}}
\renewcommand{\subsection}[1] {\vspace{12pt}\addtocounter{subsectionc}{1}
        \setcounter{subsubsectionc}{0}\noindent
        {\bf\thesectionc.\thesubsectionc.
        {\kern1pt \bfit #1}}\par\vspace{5pt}}
\renewcommand{\subsubsection}[1] {\vspace{12pt}
        \addtocounter{subsubsectionc}{1}
        \noindent
        {\tenrm\thesectionc.\thesubsectionc.\thesubsubsectionc. {\kern1pt
        \it #1}}\par\vspace{5pt}}
\newcommand{\nonumsection}[1] {\vspace{12pt}\noindent{\tenbf #1}
        \par\vspace{5pt}}

\newcounter{appendixc}
\newcounter{subappendixc}[appendixc]
\newcounter{subsubappendixc}[subappendixc]

\renewcommand{\appendix}[1] {\vspace{12pt}      
        \refstepcounter{appendixc}              
        \setcounter{figure}{0}
        \setcounter{table}{0}
        \setcounter{lemma}{0}
        \setcounter{theorem}{0}
        \setcounter{corollary}{0}
        \setcounter{definition}{0}
        \setcounter{equation}{0}
        \renewcommand{\thefigure}{\Alph{appendixc}.\arabic{figure}}
        \renewcommand{\thetable}{\Alph{appendixc}.\arabic{table}}
        \renewcommand{\theappendixc}{\Alph{appendixc}}
        \renewcommand{\thelemma}{\Alph{appendixc}.\arabic{lemma}}
        \renewcommand{\thetheorem}{\Alph{appendixc}.\arabic{theorem}}
        \renewcommand{\thedefinition}{\Alph{appendixc}.\arabic{definition}}
        \renewcommand{\thecorollary}{\Alph{appendixc}.\arabic{corollary}}
        \renewcommand{\theequation}{\Alph{appendixc}.\arabic{equation}}
        \noindent{\tenbf Appendix \theappendixc #1}\par\vspace{5pt}}

\newcommand{\textlineskip}{\baselineskip=13pt}
\newcommand{\smalllineskip}{\baselineskip=10pt}

\newcommand{\copyrightheading}[1]
        {\vspace*{-2.5cm}\smalllineskip{\flushleft
        {\footnotesize Journal of Knot Theory and Its Ramifications #1}\\
        {\footnotesize \copyright\kern2pt World Scientific
         Publishing Company}\\
         }}

\newcommand{\publisher}[2]{{\begin{center}\footnotesize\smalllineskip
        Received #1\\
        Revised #2
        \end{center}
        }}

\def\abstracts#1#2#3#4{{
        \centering{\begin{minipage}{4.5in}\footnotesize\baselineskip=10pt
        \centerline{ABSTRACT}
        \parindent=15pt #1\par
        \parindent=15pt #2\par
        \parindent=15pt #3\par
        \parindent=15pt #4\par
        \end{minipage}}\par}}

\def\keywords#1{{
        \centering{\begin{minipage}{4.5in}\footnotesize\baselineskip=10pt
        {\footnotesize\it Keywords}\/: #1
        \end{minipage}}\par}}

\renewenvironment{thebibliography}[1]
        {\frenchspacing
         \ninerm\baselineskip=11pt
         \begin{list}{[\arabic{enumi}]}
        {\usecounter{enumi}\setlength{\parsep}{0pt}
         \setlength{\leftmargin 13.7pt}{\rightmargin 0pt} 
         \setlength{\itemsep}{0pt} \settowidth
        {\labelwidth}{[#1]}\sloppy}}{\end{list}}

\newcounter{itemlistc}
\newcounter{romanlistc}
\newcounter{alphlistc}
\newcounter{arabiclistc}

\newcommand{\fcaption}[1]{
        \refstepcounter{figure}
        \setbox\@tempboxa = \hbox{\footnotesize Fig.~\thefigure. #1}
        \ifdim \wd\@tempboxa > 5in
           {\begin{center}
        \parbox{5in}{\footnotesize\smalllineskip Fig.~\thefigure. #1}
            \end{center}}
        \else
             {\begin{center}
             {\footnotesize Fig.~\thefigure. #1}
              \end{center}}
        \fi}

\newcommand{\tcaption}[1]{
        \refstepcounter{table}
        \setbox\@tempboxa = \hbox{\footnotesize Table~\thetable. #1}
        \ifdim \wd\@tempboxa > 5in
           {\begin{center}
        \parbox{5in}{\footnotesize\smalllineskip Table~\thetable. #1}
            \end{center}}
        \else
             {\begin{center}
             {\footnotesize Table~\thetable. #1}
              \end{center}}
        \fi}

\def\pmb#1{\setbox0=\hbox{#1}
        \kern-.025em\copy0\kern-\wd0
        \kern.05em\copy0\kern-\wd0
        \kern-.025em\raise.0433em\box0}

\def\fnt#1#2{\footnotetext{\kern-.3em
        {$^{\mbox{\scriptsize #1}}$}{#2}}}

\def\fpage#1{\begingroup
\voffset=.3in
\thispagestyle{empty}\begin{table}[b]\centerline{\footnotesize #1}
        \end{table}\endgroup}

\def\runninghead#1#2{\pagestyle{myheadings}
\markboth{{\protect\footnotesize\it{\quad #1}}\hfill}
{\hfill{\protect\footnotesize\it{#2\quad}}}}

\font\tenrm=cmr10

\font\tenbf=cmbx10
\font\bfit=cmbxti10 at 10pt
\font\ninerm=cmr9
\font\nineit=cmti9
\font\ninebf=cmbx9
\font\eightrm=cmr8

\newtheorem{theorem}{Theorem}   

\newtheorem{lemma}{Lemma}

\newtheorem{definition}{Definition}
\newtheorem{conjecture}{Conjecture}

\newtheorem{corollary}{Corollary}
\newtheorem{proposition}{Proposition}
\newtheorem{remark}{Remark}
\def\@begintheorem#1#2{\trivlist        
        \item[\hskip\labelsep{\bf #1\ #2.}]}
\def\@opargbegintheorem#1#2#3{\trivlist
        \item[\hskip\labelsep{\bf #1\ #2\ (#3).}]}

\newcommand{\refthm}[1]{Theorem~\ref{#1}}
    
    \newcommand{\refconj}[1]{Conjecture~\ref{#1}}
    \newcommand{\refcor}[1]{Corollary~\ref{#1}}
    \newcommand{\refdefi}[1]{Definition~\ref{#1}}
    \newcommand{\reflem}[1]{Lemma~\ref{#1}}
    \newcommand{\refprop}[1]{Proposition~\ref{#1}}
    \newcommand{\refrem}[1]{Remark~\ref{#1}}

\newenvironment{proof}{\begin{trivlist}
        \item[\noindent]{\it Proof.}}{\hfill \qed\kern0.8pt\end{trivlist}}

\newenvironment{proofa}{\begin{trivlist}
        \item[\noindent]{\it Proof }}{\hfill \qed\kern0.8pt\end{trivlist}}

        {\setcounter{itemlistc}{0}              
         \begin{list}{$\bullet$}                
        {\usecounter{itemlistc}                 
         \leftmargin10pt               
         \setlength{\parsep}{0pt}
         \setlength{\itemsep}{0pt}     
        }}{\end{list}}

        {\setcounter{romanlistc}{0}             
         \begin{list}{$($\roman{romanlistc}$)$} 
        {\usecounter{romanlistc}                
         \leftmargin18pt 
         \setlength{\parsep}{0pt}
         \setlength{\itemsep}{0pt}      
         \settowidth{\labelwidth}{#1}
        }}{\end{list}}

        {\setcounter{enumii}{0}                 
         \begin{list}{$($\alph{enumii}$)$}      
        {\usecounter{enumii}                    
         \leftmargin18pt                
         \setlength{\parsep}{0pt}
         \setlength{\itemsep}{0pt}      
         \settowidth{\labelwidth}{#1}
        }}{\end{list}}

\def\qed{\hbox{${\vcenter{\vbox{                        
   \hrule height 0.4pt\hbox{\vrule width 0.4pt height 6pt
   \kern5pt\vrule width 0.4pt}\hrule height 0.4pt}}}$}}

\renewcommand{\thefootnote}{\fnsymbol{footnote}}  

\def\theequation{\thesectionc.\arabic{equation}}  

\def\undersmile#1{\lower5.7pt\hbox{$\smallsmile$}\kern-0.6em #1}

\newcommand{\C}{{\mathbb{C}}}
\newcommand{\J}{{\mathbb{J}}}

\newcommand{\Pro}{{\mathbb{P}}}
\newcommand{\Q}{{\mathbb{Q}}}
\newcommand{\R}{{\mathbb{R}}}

\newcommand{\Z}{{\mathbb{Z}}}

\newcommand{\frsl}{{\mathfrak{sl}}}

\newcommand{\mA}{{\mathcal A}}
\newcommand{\mC}{{\mathcal C}}
\newcommand{\mD}{{\mathcal D}}
\newcommand{\mM}{{\mathcal M}}
\newcommand{\mN}{{\mathcal N}}
\newcommand{\mS}{{\mathcal S}}

\newcommand{\lieta}{\lim_{\eta \rightarrow 0_{+}}}
\newcommand{\livep}{\lim_{\vep \rightarrow 0_{+}}}
\newcommand{\la}{\langle}
\newcommand{\ra}{\rangle}
\newcommand{\qsum}{\sum_{n \in \Z/|q|\Z}}

\newcommand{\Arg}{\mbox{\rm Arg}}
\newcommand{\Cl}{\mbox{\rm Cl}}
\newcommand{\CS}{\mbox{\rm CS}}
\newcommand{\diag}{\mbox{\rm diag}}
\newcommand{\dS}{\mbox{\rm S}}
\newcommand{\dte}{\mbox{\rm d}}
\newcommand{\Hom}{\mbox{\rm Hom}}

\newcommand{\ima}{\mbox{\rm Im}}
\newcommand{\Log}{\mbox{\rm Log}}
\newcommand{\Li}{\mbox{\rm Li}}
\newcommand{\main}{\mbox{\rm main}}
\newcommand{\nab}{\mbox{\rm nab}}
\newcommand{\nbd}{\mbox{\rm nbd}}

\newcommand{\Os}{\mbox{\rm O}}
\newcommand{\rea}{\mbox{\rm Re}}
\newcommand{\Res}{\mbox{\rm Res}}
\newcommand{\sign}{\mbox{\rm sign}}
\newcommand{\SL}{\mbox{\rm SL}}
\newcommand{\tr}{\mbox{\rm tr}}
\newcommand{\SU}{\mbox{\rm SU}}
\newcommand{\Vol}{\mbox{\rm Vol}}

\newcommand{\I}{\sqrt{-1}}

\newcommand{\ria}{\rightarrow}
\newcommand{\sm}{\setminus}
\newcommand{\vep}{\varepsilon}

\newcommand{\HS}{\noindent \hfill \qed\kern0.8pt}

\begin{document}

\setlength{\textheight}{7.7truein}  
\runninghead{Asymptotics of the quantum invariants for the
figure 8 knot}
{Asymptotics of the quantum invariants for the figure 8 knot}

\normalsize\textlineskip
\thispagestyle{empty}
\setcounter{page}{1}

\copyrightheading{}   

\vspace*{0.88truein}

\fpage{1}
\centerline{\bf ASYMPTOTICS OF THE QUANTUM INVARIANTS FOR
SURGERIES}
\baselineskip=13pt
\centerline{\bf ON THE FIGURE 8 KNOT}
\vspace*{0.37truein}
\centerline{\footnotesize J\O RGEN ELLEGAARD ANDERSEN
\footnote{This research was conducted in part by the first
author for the Clay Mathematics Institute at University of
California, Berkeley.}\footnote{This work was supported
        by  MaPhySto -- A Network in Mathematical Physics and
        Stochastics, funded by The Danish National Research
        Foundation}}
\baselineskip=12pt
\centerline{\footnotesize\it Department of Mathematical
Sciences}
\baselineskip=10pt
\centerline{\footnotesize\it University of Aarhus, Building
530, Ny Munkegade, 8000 Aarhus, Denmark}
\baselineskip=10pt
\centerline{\footnotesize\it andersen@imf.au.dk}

\vspace*{10pt}
\centerline{\footnotesize S\O REN KOLD HANSEN}
\baselineskip=12pt
\centerline{\footnotesize\it Department of Mathematics,
Kansas State University}
\baselineskip=10pt
\centerline{\footnotesize\it 138 Cardwell Hall, Manhattan,
KS 66506, USA}
\baselineskip=10pt
\centerline{\footnotesize\it hansen@math.ksu.edu}

\vspace*{0.225truein}
\publisher{(Leave 1 inch blank space for publisher.)}

\vspace*{0.21truein}
\abstracts{We investigate the Reshetikhin--Turaev invariants
associated to $\SU(2)$ for the $3$--manifolds $M$
obtained by doing any rational surgery along the figure $8$
knot. In particular, we express these invariants in terms of
certain complex double contour integrals. These integral
formulae allow us to propose a formula for the leading
asymptotics of the invariants in the limit of large quantum
level. We analyze this expression using the saddle point
method. We construct a certain surjection from the set of
stationary points for the relevant phase functions onto the
space of conjugacy classes of nonabelian
$\SL(2,\C)$--representations of the fundamental group of $M$
and prove that the values of these phase functions at the
relevant stationary points equals the classical Chern--Simons
invariants of the corresponding flat $\SU(2)$--connections. Our
findings are in agreement with the
asymptotic expansion conjecture. Moreover, we calculate the
leading asymptotics of the colored Jones polynomial of the
figure $8$ knot following Kashaev \cite{[14]}. This leads to a
slightly finer asymptotic description of the invariant than
predicted by the volume conjecture \cite{[24]}.}{}{}{}

\vspace*{10pt}

\vspace*{1pt}\textlineskip    

\vspace*{-0.5pt}

\setcounter{footnote}{0}
\renewcommand{\thefootnote}{\alph{footnote}}


\section{Introduction}

\noindent In this paper we investigate the large level
asymptotics of the Reshetikhin--Turaev invariants of the
$3$--manifolds obtained by doing surgery on the figure $8$
knot with an arbitrary rational surgery coefficient. Let $X$
be a closed oriented $3$--manifold and let $\tau_r(X)$ be the
RT--invariant of $X$ associated to $\frsl_{2}(\C)$ at level
$r$, some integer $\geq 2$. The investigations of this paper
are motivated by the following conjecture.

\begin{conjecture}[Asymptotic expansion conjecture
(AEC)]\label{conj:AEC}
There exist constants (depending on X) $d_j \in \Q$ and
$b_j \in \C$ for $j=0,1, \ldots, n$ and $a_j^l \in  \C$ for
$j=0,1, \ldots, n$,\, $l=1,2,\ldots$ such that the asymptotic
expansion of $\tau_{r}(X)$ in the limit $r \ria \infty$ is
given by
\[\tau_r(X) \sim \sum_{j=0}^n e^{2\pi i r q_j} r^{d_j}
b_j \left( 1 + \sum_{l=1}^\infty a_j^l r^{-l}\right), \]
where $q_0 = 0, q_1, \ldots, q_n$ are the finitely many
different values of the Chern--Simons functional on the space
of flat $\SU(2)$--connections on $X$.
\end{conjecture}

Here {\bf $\sim$} means {\bf asymptotic expansion} in the
Poincar\'{e} sense, which means the following: Let
\[d = \max\{d_0,\ldots,d_n\}.\]
Then, for any non-negative integer $L$, there is a $c_L \in \R$
such that
\[\left| \tau_r(X) - \sum_{j=0}^n e^{2\pi i r q_j} r^{d_j}
b_j \left( 1 + \sum_{l=0}^L a_j^l
r^{-l}\right) \right| \leq c_L r^{d-L-1}\]
for all levels $r$. Of course such a condition only puts
limits on the large $r$ behaviour of $\tau_{r}(X)$.

Given an arbitrary complex function $f(r)$ defined on the
positive integers a little simple argument shows that there
at most exists one list of numbers $n \in \{0,1,\ldots\}$,
$q_0,q_1,\ldots,q_n \in \R \cap [0,1[$, $d_j \in \Q$
and $b_j \in \C$ for $j=0,1,\ldots, n$
and $a_j^l \in  \C$ for $j=0,1, \ldots, n$,\,
$l=1,2,\ldots$ such that $q_0 < q_1 < \ldots < q_n$,
$b_j \neq 0$ for $j>0$ and such
that the large $r$ asymptotic expansion of $f(r)$ is given by
\[f(r) \sim \sum_{j=0}^n e^{2\pi i r q_j} r^{d_j}
b_j \left( 1 + \sum_{l=1}^\infty a_j^l r^{-l}\right). \]
This implies that if the function $r \mapsto \tau_r(X)$ has an
asymptotic expansion of the form stated in \refconj{conj:AEC}
then the $q_j$'s, $d_j$'s, $b_j$'s and the $a_j^l$'s are all
uniquely determined by the set of invariants
$\{ \tau_r(X) \}_{r \geq 2}$, hence they are also topological
invariants of $X$. As stated above the AEC already includes the
claim that the $q_j$'s are the Chern--Simons invariants. There
are also conjectured topological formulae for the $d_j$'s, and
$b_j$'s (see e.g.\ \cite{[2]} and the references given there).

In general, there should be expressions for each of the
$a_j^l$'s in terms of sums over Feynman diagrams of certain
contributions determined by the Feynman rules of the
Chern--Simons theory. This has not yet been worked out in
general, except in the case of an acyclic flat connection and
the case of a smooth non-degenerate component of the moduli
space of flat connections by Axelrod and Singer,
cf.\ \cite{[4]}, \cite{[5]}, \cite{[3]}.

The AEC, \refconj{conj:AEC}, however offers in a sense a
converse point of view, where one seeks to derive the final
output of perturbation theory after all cancellations have
been made (i.e.\ collect all terms with the same Chern--Simons
value). This seems actually rather reasonable in this case,
since the exact invariant is known explicitly.

The AEC was proved by Andersen in \cite{[1]} in the case of
mapping tori of finite order
diffeomorphisms of orientable surfaces of genus at least two
using the gauge-theoretic approach to the quantum invariants.
Later on the AEC was proved by Hansen in \cite{[9]} for all
Seifert manifolds with orientable base by supplementing the
work of Rozansky \cite{[31]} with the required analytic
estimates. In \cite{[10]} the AEC is futher proved for the
Seifert manifolds with nonorientable base of even genus.

Using the approach of Reshetikhin and Turaev to the quantum
invariants, the AEC has not yet been proved for any hyperbolic
$3$--manifold. It is therefore particular interesting to
consider surgeries on the figure $8$ knot. Let $M_{p/q}$ be
the manifold obtained by (rational) Dehn surgery on the
figure $8$ knot with surgery coefficient $p/q$. Then $M_{p/q}$
has a hyperbolic structure if and only if $|p|>4$ or $|q|>1$,
see e.g.\ \cite[Theorem 10.5.10]{[25]} or \cite{[33]}.
(It is well-known that for $|q|=1$ and $|p| \in \{1,2,3\}$
$M_{p/q}$ is a Seifert manifold, and $M_{\pm 4}$ are obtained
by gluing together the complement of the trefoil knot and the
non-trivial $I$--bundle over the Klein bottle, see e.g.\
\cite[p.~95]{[15]}. The manifold $M_0$ is a mapping torus of
a torus, see Appendix C.) We use here the convention of Rolfsen
for surgery coefficients, cf.\ \cite[Chap.~9]{[30]}. In
particular Dehn surgery on a knot $K$ in $S^{3}$ with surgery
coefficient $f \in \Z$ is equal to the boundary of the compact
$4$--manifold obtained by attaching a $2$--handle to the
$4$--ball using the knot $K$ with framing $f$, see
\cite[p.~261]{[30]}. As usual $M_{p/q}$ is given the
orientation induced by the standard right-handed orientation
of $S^{3}$.

The advantage of working with surgeries on the figure $8$ knot
$K$ is that the (normalized) colored Jones polynomial
$J_K'(\lambda)$ is known explicitly. In fact
\[J'_K(\lambda) =  \sum_{m=0}^{\lambda-1} \xi^{-m\lambda}
\prod_{l=1}^{m}(1-\xi^{\lambda-l}) (1-\xi^{\lambda+l}), \]
where $\xi=\exp(2\pi i/r)$ (and the product is $1$ for $m=0$).
The colors $\lambda$ are here dimensions of irreducible
representations of the quantum group associated to
$\frsl_{2}(\C)$ and the root of unity $\xi$, so
$\lambda =1,2,\ldots,r$. By the above expression for
$J_{K}'(\lambda)$ we have an explicit formula for the quantum
invariant $\tau_r(M_{p/q})$ (see formula (\ref{eq:inv2})).
Although this formula is completely explicit, it is not clear
from it what the leading order large $r$ asymptotics of
$\tau_r(M_{p/q})$ is. In order to study this asymptotics, we
observe (generalizing from Kashaev's work) that the product
in the expression for the colored Jones polynomial can be
expressed in terms of a quotient of two evaluations of
Faddeev's quantum dilogarithm $S_\gamma$ ($\gamma = \pi/r$):
\begin{equation}\label{ncJones}
J'_K(\lambda) =  \sum_{m=0}^{\lambda-1}
\frac{\xi^{-m\lambda}}{(1-\xi^{\lambda})}
\frac{S_{\gamma}(-\pi+2\gamma(\lambda - m) -\gamma)}
{S_{\gamma}(-\pi+2\gamma(\lambda + m) +\gamma)}.
\end{equation}
This follows directly from the functional equation
$$
(1+e^{\I \zeta})S_{\gamma}(\zeta+\gamma) =
S_{\gamma}(\zeta-\gamma)
$$
which Faddeev's $S_\gamma$ satisfies. Recall that for
$\rea(\zeta) < \pi +\gamma$ we have the expression
\[S_{\gamma}(\zeta) = \exp \left( \frac{1}{4} \int_{C_{R}}
\frac{e^{\zeta z}}{\sinh(\pi z)\sinh(\gamma z)z} \dte z\right),
\]
which together with the functional equation determines
$S_\gamma$ as a meromorphic function on $\C$. For the so-called
top color, i.e.\ $\lambda=r$, we obtain the sligthly simpler
expression
$$
J'_K(r) =  r\sum_{m=0}^{r-1}
\frac{S_{\gamma}(\pi - (2m+1)\gamma)}
{S_{\gamma}(-\pi + (2m+1)\gamma)}.
$$
Then we simply use the residue formula to convert this
expression into a contour integral
$$
J_K'(r) =  \int_{C_r} \tan(\pi r x) \tilde{g}_{r}(x) \dte x,
$$
where $C_r$ is contained in the strip
$\{ z \in \C \mid 0 < \rea(z) < 1 \}$ and encloses
$(m + \frac12)/r$ for $m=0, \ldots, r-1$, and $\tilde{g}_{r}$
is a holomorphic function in this strip expressed in terms of
the above quotient of $S_\gamma$ functions (see formula
(\ref{eq:contourfigure8})).

Similarly we get for the quantum invariant with the use of
(\ref{ncJones}) and the residue theorem, now a double contour
integral, since the quantum invariant also involves a sum over
colors:
\begin{equation}\label{qidci}
\tau_r(M_{p/q}) = \int_{C_r\times C_r} \cot(\pi r x)
\tan(\pi r y) \tilde{f}_{p,q,r}(x,y) \dte x \dte y,
\end{equation}
where $\tilde{f}_{p,q,r}$ is holomorphic on the double strip
$\{ z \in \C \mid 0 < \rea(z) < 1 \}^{2}$ and given by some
expression involving quotients of evaluations of
$S_\gamma$ functions, and where we furthermore require of
$C_r$ that it also encloses $k/r$ for $k=1, \ldots, r-1$.

From this it is clear that we need to understand the small
$\gamma$ asymptotics of $S_\gamma$. We have that
\[S_\gamma(\zeta) = \exp\left( \frac{1}{2\I \gamma}
\Li _{2} (-e^{\I \zeta} ) + I_\gamma(\zeta)\right),\]
where $\Li_{2}$ is Euler's dilogarithm function, and where we
have certain analytic estimates on $I_\gamma(\zeta)$ (see
\reflem{lem:Sapproximation}).

Let us first explain how we use this to give a proof of the
volume conjecture of Murakami and Murakami \cite{[24]} for
the figure $8$ knot, namely that
$$
\lim_{r \ria \infty} \frac{2\pi\Log(J_{K}'(r))}{r} =
\Vol(K),
$$
where the right-hand side is the hyperbolic volume of the
figure $8$ knot, i.e.\ the hyperbolic volume of the complement
$S^{3} \sm K$. The basic idea in analyzing the above contour
integral expression for $J_K'(r)$ is the following. In the
upper half plane we approximate $\tan$ by $\I$ and by $-\I$ in
the lower half plane. Further we approximate $S_\gamma$ by the
above expression involving only the dilogarithm. In Appendix B
we prove the required estimates which allows us to do these
approximations and we end up with the following formula for the
leading order large $r$ asymptotics
\begin{equation}\label{eq:Asymf8'}
J'_K(r) \sim r^{2}  \int_\epsilon^{1-\epsilon}
e^{r\Phi(x)} \dte x,
\end{equation}
where $\epsilon <1/(4r)$ is a small positive parameter and
\[\Phi(x)= \frac{1}{2\pi\I}
\left( \Li_{2}(e^{-2\pi\I x}) - \Li_{2}(e^{2\pi\I x})\right).\]
Now we simply analyze the integral on the right-hand side of
(\ref{eq:Asymf8'}) by the saddle point method. This consists of
finding the stationary points of $\Phi$ and also the
directions of steepest descent, see e.g.\ \cite{[6]}.
(In this paper we use critical point and stationary point
interchangeably to mean a point in which the derivative is
zero.) This analysis leads to two interesting results. Firstly,
the search for stationary points leads to the hyperbolicity
equation for the complement of the figure $8$ knot in the
$3$--sphere. Recall that this complement can be decomposed into
two ideal hyperbolic tetrahedra each parametrized by
a certain complex number. This decomposition then defines a
hyperbolic structure on the complement exactly when the two
parameters are equal and satisfy the hyperbolicity equation.

Secondly, we find that the value of the phase function in the
relevant stationary point (there is only one such point in
this case) is equal to the hyperbolic volume of the knot
complement (divided by $2\pi$), hence the leading asymptotics
of $J_{K}'(r)$ is determined by this volume.

These phenomena were first observed by Kashaev \cite{[14]} and
have been conjectured by Thurston \cite{[32]} and Yokota
\cite{[35]} to be generally true for hyperbolic knots
(see \refrem{rem:hyperbolicity}).

Ultimately our asymptotic analysis leads to the following

\begin{theorem}
The leading order large $r$ asymptotics of the colored Jones
polynomial evaluated at the top color is given by
$$
J_{K}'(r) \sim 3^{-1/4} r^{3/2}
\exp\left(\frac{r}{2\pi}\Vol(K)\right).
$$
\end{theorem}

As a corollary we obtain the volume conjecture for the figure
$8$ knot. We note that none of the proofs so far given in the
literature for the volume conjecture for the figure $8$ knot
have been able to see the finer details of the asymptotic
behaviour, namely the polynomial part $3^{-1/4}r^{3/2}$.

After completion of this paper we have been informed that
D.\ Zagier has computed the full asymptotic expansion of
$J_{K}'(r)$ by using the Euler--Maclaurin summation formula,
however his techniques seem not applicable to the calculation
of the large $r$ asymptotics of $\tau_r(M_{p/q})$.

Let us now return to the study of the large $r$ asymptotics of
the quantum invariant $\tau_{r}(M_{p/q})$. We expect that an
analysis of the expression (\ref{qidci}) paralleling our
analysis of $J'_K(r)$ should be applicable. I.e.\ $\tan$ and
$\cot$ should be approximated by $\pm \I$ depending on the
signs of $\ima(y)$ and $\ima(x)$ and $\tilde{f}_{p,q,r}(x,y)$
by an appropriate expression involving the dilogarithm for some
deformation of the relevant part of $C_r\times C_r$. We have
partial analytic results supporting this.

To be more specific we propose the following analog of
(\ref{eq:Asymf8'}) for the quantum invariant. Let $d$ be the
inverse of $p \bmod{q}$ and let $(a,b)\in \{0,1\}^2$ and
$n \in \Z$. Define
$$
\Phi_{n}(x,y) = -\frac{dn^{2}}{q} - \frac{p}{4q} x^{2}
  + \frac{n}{q}x-xy + \frac{1}{4\pi^{2}}
  \left( \Li_{2}(e^{2\pi i(x+y)}) - \Li_{2}(e^{2\pi i (x-y)})
\right),
$$
and
\begin{equation}\label{Phi}
\Phi_{n}^{a,b}(x,y) = a(x+y) + b(x-y) + \Phi_{n}(x,y).
\end{equation}
Finally, let
\begin{equation}\label{eq:SUcriticalpoint}
\mS = \{ (x,y) \in \R \times \C \; | \; e^{2\pi iy} \in
]-\infty,0[ \;\}.
\end{equation}
Thus $\mS$ is the union of infinitely many planes,
namely the ones characterized by $x \in \R$ and
$\rea(y) \in \frac{1}{2} + \Z$.

\begin{conjecture}\label{conj:MC}
There exist surfaces
$\tilde{\Sigma}_{a,b}^{\mu,\nu,n} \subset \C^2$
for $(a,b) \in \{0,1\}^2$, $n \in \Z / |q|\Z$ and
$(\mu,\nu) \in \{\pm 1\}^2$ such that the leading order large
$r$ asymptotics of the quantum invariant is given by
\begin{eqnarray}\label{eq:asympqiI}
\bar{\tau}_{r}(M_{p/q}) &\sim& C r \qsum \sum_{(a,b) \in
\{0,1\}^{2}} \sum_{(\mu,\nu) \in \{ \pm 1 \}^{2}} \mu\nu \\ &&
\times \int_{\tilde{\Sigma}_{a,b}^{\mu,\nu,n}} \tilde{g}_n(x)
e^{2\pi ir\Phi_{n}^{a,b}(x,y)} \dte x \dte y,
\end{eqnarray}
where $C$ is a constant only depending on $p$ and $q$ and
$\tilde{g}_n$ is some simple $r$--independent function of
$x \in \C$. If $p/q \neq 0$ the surfaces
$\tilde{\Sigma}_{a,b}^{\mu,\nu,n} \subset \C^2$ can be chosen
such that any critical point of $\Phi_{n}^{a,b}$ belonging to
$\tilde{\Sigma}_{a,b}^{\mu,\nu,n}$ also belongs to the set
$\mS$.
\end{conjecture}

Please see \refconj{conj:mainconj} for the more detailed
version of this conjecture, including the precise formula for
$\tilde{g}_n$. (We have here for sign-reasons switched to the
complex conjugate invariant
$\bar{\tau}_{r}(M_{p/q})=\overline{\tau_{r}(M_{p/q})}
=\tau_{r}(M_{-p/q})$.) In case $p/q=0$ our results show that
we have to include critical points not belonging to $\mS$ in
our surfaces $\tilde{\Sigma}_{a,b}^{\mu,\nu,n}$, see 
Appendix C.

We now proceed by making an asymptotic analysis of the
right-hand side of (\ref{eq:asympqiI}) using the saddle point
method like in our analysis of (\ref{eq:Asymf8'}). Thus we
need to analyze integrals of the form
\begin{eqnarray}\label{eq:Int}
I_{a,b}^{\mu,\nu,n} = \int_{\tilde{\Sigma}_{a,b}^{\mu,\nu,n}}
\tilde{g}_n(x) e^{2\pi ir\Phi_{n}^{a,b}(x,y)} \dte x \dte y.
\end{eqnarray}
Again we have to determine the stationary points of
$\Phi_{n}^{a,b}$ and the values of $\Phi_{n}^{a,b}$ in the
relevant stationary points. The main idea behind the saddle
point method is to deform $\tilde{\Sigma}_{a,b}^{\mu,\nu,n}$ 
so that it contains certain stationary points of 
$\Phi_{n}^{a,b}$ satisfying that the leading order large $r$ 
asymptotics of $I_{a,b}^{\mu,\nu,n}$ is given by integrating 
the integrant of the integral (\ref{eq:Int}) over small 
neighborhoods of these stationary points.

If we let $v = e^{\pi i x}$ and $w = e^{2\pi i y}$, then by
exponentiating the two equations for $(x,y)$ being a
stationary point of $\Phi_{n}^{a,b}$ (see
\refthm{thm:main'} below) we obtain the equations
\begin{eqnarray}\label{eq:stationary'}
v^{-p} &=& \left(\frac{w-v^{2}}{1-v^{2}w}\right)^{q},
\nonumber \\*
v^2w &=& (1-v^{2}w)(w-v^{2}),
\end{eqnarray}
which are independent of the integer parameters $a,b,n$. To
link the asymptotics to the flat connections (as proposed by
the AEC) we then have to relate the relevant stationary points
of the phase functions $\Phi_{n}^{a,b}$ to the classical
$\SU(2)$ Chern--Simons theory on the manifolds $M_{p/q}$.
Fortunately, this Chern--Simons theory has been given a
detailed description by Kirk and Klassen \cite{[18]} using the
work of Riley \cite{[28]}, \cite{[29]} on the $\SL(2,\C)$
representation variety of the knot group of the figure $8$
knot. According to Riley the conjugacy classes of the
nonabelian elements of this variety can be represented
by a certain set of representations $\rho$ parametrized by
$\rho=\rho_{(s,u)}$, where $(s,u) \in \C^{*} \times \C$
satisfies a certain polynomial equation. (To be precise
$\rho_{(s_1,u_1)}$ and $\rho_{(s_2,u_2)}$ are conjugate if and
only if $u_2=u_1 \neq 0$ and $s_2 \in \{ s_1,s_1^{-1} \}$ or
$u_2=u_1=0$ and $s_2=s_1 \in A$, where $A$ is a certain subset
of $\C^{*}$ consisting of $4$ points.)

We show that $\rho_{(s,u)}$ defines a
$\SL(2,\C)$--representation of $\pi_{1}(M_{p/q})$ if and only
if $(v,w)=(s,u+1)$ is a solution to (\ref{eq:stationary'}) and
$v^{2} \neq 1$. By a result of Riley it is known that this
representation is conjugate to a $\SU(2)$--representation if
and only if $(s,u) \in S^{1} \times \R$. We relate the results
of Kirk and Klassen \cite{[18]} on the Chern--Simons invariants
of flat $\SU(2)$--connections on the manifolds $M_{p/q}$ to the
asymptotic analysis of the quantum invariants
$\bar{\tau}_{r}(M_{p/q})$ via a detailed analysis of the
relevant critical values of the phase functions $\Phi_n^{a,b}$.
Ultimately we arrive at

\begin{theorem}\label{thm:main'}
The map
$$
(x,y) \mapsto [\rho_{(e^{\pi ix},e^{2\pi i y} -1)}]
$$
gives a surjection from the set of critical points $(x,y)$ of
the functions $\Phi_{n}^{a,b}$ with $x \notin \Z$ onto the
set of conjugacy classes of nonabelian
$\SL(2,\C)$--representations of $\pi_{1}(M_{p/q})$. Moreover,
$(x,y) \in \C^{2}$ is a critical point of $\Phi_{n}^{a,b}$ if
and only if
\begin{eqnarray}\label{eq:stationary''}
2a + \frac{n}{q} &=& y + \left( \frac{p}{2q} + 1 \right)x
     + \frac{i}{\pi}\Log\left( 1 - e^{2\pi i (x+y)} \right),
\nonumber \\*
2b + \frac{n}{q} &=& y + \left( \frac{p}{2q} - 1 \right)x
     - \frac{i}{\pi}\Log\left( 1 - e^{2\pi i (x-y)} \right).
\end{eqnarray}
Futhermore, if $(x,y)$ is a critical point of $\Phi_{n}^{a,b}$
then $\rho_{(e^{\pi ix},e^{2\pi i y} -1)}$ is equivalent to
a $\SU(2)$--representation $\bar{\rho}$ of $\pi_{1}(M_{p/q})$
if and only if $(x,y) \in \mS$ (see
(\ref{eq:SUcriticalpoint})), and in that case
$$
\CS(\bar{\rho}) = \Phi_{n}^{a,b}(x,y) \pmod{\Z},
$$
where $\CS$ is the Chern--Simons functional on the space of
flat $\SU(2)$--connections on $M_{p/q}$.
\end{theorem}

Continuing our asymptotic analysis of the integral
(\ref{eq:Int}) we point out that not all critical points have
to give contributions to the asymptotics. The simplest case
arise when all the relevant critical points are non-degenerate.
Following \refconj{conj:MC} we claim that the only
critical points giving a contribution to the leading order
large $r$ asymptotics of $\tau_r(M_{p/q})$ are (some of) the
ones belonging to the set $\mS$ in (\ref{eq:SUcriticalpoint}).
In \refprop{prop:Nondeg} in Sec.~4 we prove that all these
critical points are non-degenerate in case $|p/q| < \sqrt{20}$.
In case $|p/q| > \sqrt{20}$ we have the following partial
result:
The surjection described in the first part of
\refthm{thm:main'} restricts to a surjection $\phi$ from the
set of critical points belonging to $\mS$ onto the moduli space
$\mM_{p/q}'$ of flat irreducible $\SU(2)$--connections on
$M_{p/q}$. If $p/q > \sqrt{20}$ we prove that all the
critical points belonging to $\mS \sm \mC_0$ are
non-degenerate, where $\mC_0$ is a certain set of critical
points. To be more precise the set $\mC_0$ is either empty or
it is equal to the preimage of $\phi$ of one point if $p$ is
odd and of two points if $p$ is even. We expect that the
argument equation following from the first of the equations
(\ref{eq:stationary'}) together with the last statement in
\refprop{prop:Nondeg} should establish that $\mC_0$ is empty.

From \refconj{conj:MC} we see that the only
relevant phase functions $\Phi_n^{a,b}$ are the ones with
$a,b \in \{0,1\}$. This, however, does not imply that only a
proper subset of the conjugacy classes of
$\SL(2,\C)$--representations of $\pi_1(M_{p/q})$ are in play.
In fact, by using (\ref{eq:stationary''}), we see that
if $(x,y)$ is a critical point of $\Phi_n^{a,b}$ then
$(x+2k,y+2l)$ is a critical point of
$\Phi_{n+pk}^{a+l+k,b+l-k}$ for any $k,l \in \Z$. Moreover,
the critical points $(x,y)$ and $(x+2k,y+2l)$ correspond to
the same conjugacy class of representations by the surjection
in \refthm{thm:main'}. In particular all points of the
moduli space of flat irreducible $\SU(2)$--connections on
$M_{p/q}$ could potentially give a contribution to the
asymptotics of $\bar{\tau}_r(M_{p/q})$.

One is only likely to succeed in using the saddle point method
to calculate the large $r$ asymptotics of an integral of the
form (\ref{eq:Int}) in case one can deform the surface
$\tilde{\Sigma}_{a,b}^{\mu,\nu,n}$ so that it only contains 
critical point of a certain nice kind. Here we will only 
consider the case of non-degenerate critical points. Therefore,
assume that $(x,y)$ is a non-degenerate critical point of 
$\Phi_n^{a,b}$ belonging to $\tilde{\Sigma}_{a,b}^{\mu,\nu,n}$.
Then to be able to calculate that critical point's contribution
to the large $r$ asymptotics of the integral in (\ref{eq:Int}) 
we assume that $(x,y)$ is positive definite, meaning that the 
imaginary part of a certain ``twisted'' Hessian of 
$\Phi_n^{a,b}$ in $(x,y)$ is positive definite, see 
\refdefi{defi:posdefinite} for the precise definition. 
Ultimatively we arive at the following result.

\begin{theorem}\label{thm:MT}
Assume that $p/q \neq 0$. IF \refconj{conj:MC} is true and if
all the critical points of the phase functions $\Phi_n^{a,b}$
belonging to $\mS$ are non-degenerate and positive definite,
then the leading order large $r$ asymptotics of the quantum
invariant $\bar{\tau}_{r}(M_{p/q})$ is given by
\begin{eqnarray}\label{eq:asympqi}
\bar{\tau}_{r}(M_{p/q}) &\sim&
\frac{\sign(q)}{4\sqrt{|q|}} e^{\frac{3\pi i}{4}\sign(pq)}
\sum_{\bar \rho\in {\mathcal M}'_{p/q}}
e^{2\pi i r CS(\bar \rho)} b_{\bar \rho},
\end{eqnarray}
where $\mM_{p/q}'$ is the moduli space of flat irreducible
$\SU(2)$--connections on $M_{p/q}$. For each
$\bar \rho \in {\mathcal M}'_{p/q}$ there are
$(a,b) \in \{0,1\}^2$, $n \in \Z$ and a stationary point
$(x,y) \in \mS$ for the function $\Phi_n^{a,b}$ given by
(\ref{Phi}), i.e.\ $(x,y)$ solves the equations
(\ref{eq:stationary''}), such that $\bar \rho$ is equivalent to
$\rho_{(e^{\pi ix}, e^{2\pi iy}-1)}$. Moreover
\[CS(\bar\rho) = \Phi_n^{a,b}(x,y) \pmod{\Z},\]
and
\[b_{\bar \rho} = m_{\bar \rho}
e^{\frac{\pi i\sigma_{\bar \rho}}{2}}
\sin\left( \frac{\pi}q(x-2nd)\right)
\left| 1-4\cos(2\pi x)) +
\frac{p}{2q}\sinh(2\pi \ima(y))\right|^{-\frac12}, \]
where $m_{\bar \rho}$ and $\sigma_{\bar \rho}$ are some
integers.
\end{theorem}

The positive integers $m_{\bar \rho}$ arrise from the following
two situations. First of all it can happen that two or more
critical points of the relevant stationary phase functions
correspond to the same irreducible flat $\SU(2)$--connection on
$M_{p/q}$. Secondly, the same stationary point for the same 
phase function $\Phi_n^{a,b}$ can belong to two or more of the 
surfaces $\tilde{\Sigma}_{a,b}^{\mu,\nu,n}$. In fact, we expect
that there is a one to one correspondence between the moduli 
space $\mM_{p/q}'$ of flat irreducible $\SU(2)$--connections on
$M_{p/q}$ and the set of stationary points contributing to the
leading large $r$ asymptotics of $\bar{\tau}_r(M_{p/q})$.
Moreover, we expect that a contributing stationary point for a
relevant phase function $\Phi_n^{a,b}$ belongs to (a 
deformation) of $\Sigma_{a,b}^{\mu,\nu,n}$ for all 
$(\mu,\nu) \in \{\pm 1\}^2$ thus causing $m_{\bar \rho}=4$ for 
all $\bar \rho \in \mM_{p/q}'$, see \refrem{rem:Remark4}. The 
above result, \refthm{thm:MT}, coincides with the one found by 
the second author in \cite{[10]} for the cases 
$|p/q| \in \{1,2,3\}$ with $m_{\bar \rho}=4$ for all
$\bar \rho$.

The invariant $\tau_{r}(M_{0})$ and its full asymptotic
expansion have been calculated by Jeffrey \cite{[13]}. We show
(see Appendix C) that Jeffrey's result is in agreement with the
AEC. In this case the reducible flat $\SU(2)$--connections
contribute to the leading asymptotics. We find that the part of
the leading order large $r$ asymptotics of $\tau_{r}(M_{0})$
associated to the irreducible flat $\SU(2)$--connections on
$M_0$ is given by the right hand side of (\ref{eq:asympqi})
with all $m_{\bar \rho}=4$.

\nonumsection{Acknowledgements}

\noindent The first author thanks the Department of Mathematics
and the Mathematical Sciences Research Institute, University
of California, Berkeley for their hospitality during several
visits, where part of this work was undertaken. The second
author thanks the Universit\'{e} Louis Pasteur, Strasbourg,
the University of Edinburgh, and the Max--Planck--Institut
f\"{u}r Mathematik, Bonn for their hospitality during this work.
He was supported by the European Commission, the Danish
Natural Science Research Council and the Max--Planck--Institut
f\"{u}r Mathematik, Bonn.

\section{The RT-invariant for surgeries on the figure $8$ knot}

\noindent This section is primarily intended to introduce
notation. Moreover, we present some preliminary formulas for
the colored Jones polynomial of the figure $8$ knot and for
the RT--invariants of the $3$--manifolds $M_{p/q}$.

Let $t=\exp(2\pi \I /(4r))$, $r$ an integer $\geq 2$, and let
$U_{t}$ be the modular Hopf algebra considered in
\cite[Sec.~8]{[27]}, i.e.\ $U_{t}$ is a finite-dimensional
factor of the quantum group $U_{\xi}(\frsl_{2}(\C))$,
$\xi=t^{4}$. (In \cite{[27]}, and in most literature on the
subject, $\xi$ is denoted $q$, but we use in this paper $q$
to mean something different.) For an integer $k$ we let
$$
[k] = \frac{t^{2k}-t^{-2k}}{t^{2}-t^{-2}} =
            \frac{\sin(\pi k /r)}{\sin(\pi/r)},
$$
sometimes called a quantum integer. For a knot $K$ in $S^{3}$
we denote by $K^{0}$ the knot $K$ considered as a framed knot
with framing zero. The colored Jones polynomial associated to
$U_{t}$ of a framed oriented knot $K$ with color
$\lambda \in \{1,2,\ldots,r\}$ is denoted $J_{K}(\lambda)$,
and for an oriented knot $K$ in $S^{3}$ we let $
J_{K}'(\lambda)=J_{K^{0}}(\lambda)/[\lambda]$. Here the
colors are the dimensions (as complex vector spaces) of
irreducible $U_{t}$--modules.

Let $N_{f}$ be the $3$--manifold obtained by surgery on
$S^{3}$ along $K$ with surgery coefficient $f \in \Z$. By
\cite{[16]} or \cite{[27]} the RT--invariant (at level $r-2$)
of $N_{f}$ is
\begin{equation}\label{eq:inv1}
\tau_{r}(N_{f}) =
  \alpha\sum_{k=1}^{r-1} \xi^{(k^{2}-1)f/4} [k]^{2}J_{K}'(k),
\end{equation}
where $K$ is given an arbitrary orientation. Here
$\alpha=C^{\sign(f)}\mD^{-2}$, where
\begin{eqnarray*}
\mD &=& \sqrt{\frac{r}{2}}\frac{1}{\sin(\pi/r)}, \\
C &=& \exp\left( \frac{\I\pi}{4} \frac{3(2-r)}{r} \right).
\end{eqnarray*}
We use here the normalization of \cite{[34]}. This is
$\mD^{-1}$ times the normalization of \cite{[16]} and
$C^{-b_{1}(N_{f})}\mD^{-1}$ times the normalization of
\cite{[27]}, where $b_{1}(N_{f})$ is the first Betti
number of $N_{f}$, see \cite[Appendix A]{[8]}. (In the
notation of \cite{[34]}, $C=\Delta \mD^{-1}$.)

Let us next generalize to arbitrary rational surgery. Let
$p,q$ be a pair of coprime integers with $q \neq 0$, and let
$N_{p/q}$ be the $3$--manifold obtained by surgery along $K$
with surgery coefficient $p/q$. Choose $c,d \in \Z$ such that
$B= \left( \begin{array}{ll}
                p & c \\
                q & d
                \end{array}
\right) \in \SL(2,\Z)$. Then (see e.g.\ \cite[Theorem 5.1 and
the proofs of Corollary 8.3 and Theorem 8.4]{[8]}),
$$
\tau_{r}(N_{p/q}) = \left( e^{\frac{i \pi}{4}}
 \exp \left( -\frac{i\pi}{2r} \right) \right)^{\Phi(B)
 -3\sign(pq) } \sqrt{\frac{2}{r}}\sin\left(\frac{\pi}{r}\right)
 \sum_{\lambda=1}^{r-1} [\lambda] J_{K}'(\lambda)
 \tilde{B}_{\lambda,1},
$$
where $\Phi$ is the Rademacher Phi function, see \cite{[26]},
and $\tilde{\cdot}$ is the unitary representation of
$P\SL(2,\Z)$ given by
\begin{eqnarray*}
&&\tilde{B}_{j,k} = \I \frac{\sign (q)}{\sqrt{2r|q|}}
      e^{-\frac{\I \pi}{4} \Phi(B)} \\*
 && \hspace{.5in} \times \sum_{\mu =
    \pm 1} \qsum \mu \exp \left( \frac{\I \pi}{2rq}
    [ p j^{2}- 2 \mu j (k +2rn \mu) + d (k + 2rn \mu)^{2} ]
    \right).
\end{eqnarray*}
By evaluating the sum over $\mu$ we get
\begin{eqnarray}\label{eq:rational}
&&\tau_{r}(N_{p/q}) = a(r)\qsum
       \exp\left( 2\pi ir \frac{dn^{2}}{q} \right) \\*
&& \hspace{.7in} \times \sum_{k=1}^{r-1}
    \sin\left(\frac{\pi}{q}\left[2nd-\frac{k}{r}\right]\right)
    \exp\left( \frac{\pi ir}{2q}
        \left[ p\left(\frac{k}{r}\right)^{2}
               - 4n\frac{k}{r}\right]
      \right) [k]J_{K}'(k), \nonumber
\end{eqnarray}
where
$$
a(r) = -\frac{2\sign(q)}{r\sqrt{|q|}}
        \sin\left(\frac{\pi}{r}\right)
    e^{-\frac{3\pi i}{4}\sign(pq)}\exp\left(\frac{\pi i}{2r}
               \left[ 3\sign(pq) - \frac{p}{q}
          + \dS \left(\frac{p}{q}\right)\right]\right).
$$
Here $\dS$ is the Dedekind symbol, see e.g.\ \cite{[17]}.
We note that the quantum invariant $\tau_{r}$ is independent
of the colored Jones polynomial $J_{K}'(k)$ for the top-color
$k=r$.

In the remaining part of this paper $K$ will denote the
figure $8$ knot unless explicitly stated otherwise. Recall
that $M_{p/q}$ denotes the $3$--manifold obtained by surgery
on $S^{3}$ along $K$ with surgery coefficient $p/q \in \Q$.
By an $R$--matrix calculation (see e.g.\ \cite{[11]}) we find
that
\begin{equation}\label{eq:figure8Le}
J_{K}'(\lambda) = \sum_{m=0}^{\lambda-1} \xi^{-m\lambda}
 \prod_{l=1}^{m}(1-\xi^{\lambda-l}) (1-\xi^{\lambda+l})
\end{equation}
for $\lambda=1,2,\ldots,r$, where
$\prod_{l=1}^{m} (1-\xi^{k-l}) (1-\xi^{k+l})=1$ for $m=0$. Le
and Habiro have obtained a similar formula, cf.\ \cite{[20]}.

\begin{remark}
Unitarity of the TQFT associated to $U_t$ implies that
\begin{equation}\label{eq:unitarity}
{\tau}_{r}(-M)=\overline{\tau_{r}(M)}
\end{equation}
for any $3$--manifold $M$, where $\overline{\cdot}$ means
complex conjugation. This formula also follows directly from
\cite{[16]} and the remarks concerning normalization
following (\ref{eq:inv1}).
\end{remark}

Since the figure $8$ knot is amphicheiral, $M_{-p/q}$ and
$M_{p/q}$ are orientation reversing homeomorphic. By
(\ref{eq:unitarity}) we therefore have
\begin{equation}\label{eq:mirror2}
\overline{\tau_{r}(M_{p/q})} = \tau_{r}(M_{-p/q}).
\end{equation}
This formula also follows directly by (\ref{eq:rational}) and
the facts that $J_{K}'(\lambda)$ is real and $S(-p/q)=-S(p/q)$.
That $J_{K}'(\lambda)$ is real follows by amphicheirality of
$K$ but can also be seen directly from (\ref{eq:figure8Le}) by
\begin{eqnarray*}
J_{K}'(\lambda) &=& \sum_{m=0}^{\lambda-1} \prod_{l=1}^{m}
   \xi^{-\lambda/2}\xi^{l/2}(1-\xi^{\lambda-l})
   \xi^{-\lambda/2}\xi^{-l/2}(1-\xi^{\lambda+l}) \\
 &=& \sum_{m=0}^{\lambda-1} \prod_{l=1}^{m}
       (\xi^{(\lambda-l)/2}-\xi^{-(\lambda-l)/2})
       (\xi^{(\lambda+l)/2}-\xi^{-(\lambda+l)/2}) \\
 &=& \sum_{m=0}^{\lambda-1} (-4)^{m} \prod_{l=1}^{m}
       \sin(\pi(\lambda-l)/r)\sin(\pi(\lambda+l)/r).
\end{eqnarray*}
By (\ref{eq:rational}) and (\ref{eq:figure8Le}) we get the
following preliminary formula for the RT--invariants of the
manifolds $M_{p/q}$:
\begin{eqnarray}\label{eq:inv2}
\tau_{r}(M_{p/q}) &=& \frac{ia(r)}{2\sin(\pi/r)}\qsum
   \exp\left( 2\pi ir \frac{dn^{2}}{q} \right) \nonumber \\*
&& \times \sum_{k=1}^{r-1} \exp\left( \frac{\pi ir}{2q}
\left[p\left(\frac{k}{r}\right)^{2}-4n\frac{k}{r}\right]
                   \right) \sin\left(\frac{\pi}{q}
         \left[2nd-\frac{k}{r}\right]\right) \nonumber \\*
&& \times \sum_{m=0}^{r-1}
         \frac{\xi^{-(m+1/2)k}}{(1-\xi^{k})}
         \prod_{l=0}^{m} (1-\xi^{k-l}) (1-\xi^{k+l}).
\end{eqnarray}

\section{A complex double contour integral formula for
$\tau_{r}(M_{p/q})$}

\noindent In this section we derive a complex double contour
integral formula for the RT--invariants $\tau_{r}(M_{p/q})$ by
using methods similar to Kashaev \cite{[14]}. When we consider
the expression for the summand in the multi sum
(\ref{eq:inv2}), we see that the expression as it stands only
makes sense for non-negative integers $m$. In order to make
sense of this expression for arbitrary complex values for $m$,
let us consider the quantum dilogarithm of Faddeev
$$
S_{\gamma}(\zeta) = \exp \left( \frac{1}{4} \int_{C_{R}}
\frac{e^{\zeta z}}{\sinh(\pi z)\sinh(\gamma z)z} \dte z\right)
$$
defined on
$\Delta_{\gamma} = \{ \zeta \in \C \; | \;
|\rea(\zeta)| < \pi +\gamma \; \}$,
where $\gamma \in ]0,1[$ and $C_{R}$ is the contour
$]-\infty,-R] + \Upsilon_{R} + [R,\infty[$, where
$\Upsilon_{R}(t)=R e^{\sqrt{-1}(\pi-t)}$, $t \in [0,\pi]$ and
$R \in ]0,1[$.

The function $S_{\gamma} : \Delta_{\gamma} \to \C$ is
holomorphic and it satisfies the following well-known
functional equation (see \cite{[7]} or \cite{[14]}).

\begin{lemma}\label{lem:functional}
For $\zeta \in \C$ with $|\rea (\zeta) | < \pi$ we have
$$
(1+e^{\I \zeta})S_{\gamma}(\zeta+\gamma) =
S_{\gamma}(\zeta-\gamma).
$$
\end{lemma}

For the sake of completeness we have given a proof in Appendix
A. We use \reflem{lem:functional} to extend $S_{\gamma}$ to a
meromorphic function on the complex plane $\C$.

From now on we fix $\gamma=\pi/r$ (so $r>3$).
By \reflem{lem:functional} we get that
$$
S_{\gamma}(\zeta) = S_{\gamma}(\zeta+2\pi)
\prod_{j=0}^{r-1} \left( 1 + e^{\I(\zeta+(2j+1)\pi/r)} \right).
$$
If we write $\zeta =-\pi+2\pi x$ we get that
$$
\prod_{j=0}^{r-1} \left( 1+e^{\I(\zeta +(2j+1)\pi/r)} \right)
= \prod_{j=0}^{r-1} \left( 1-w^{j}e^{2\pi\I (x+\frac{1}{2r})}
\right),
$$
where $w=e^{2\pi\I/r}$. Using
$1 - z^{r} = \prod_{j=0}^{r-1} (1-w^{j}z)$ we get that
\begin{equation}\label{eq:functional}
S_{\gamma}(-\pi+2\pi x) = \left(1+e^{2\pi \I xr} \right)
S_{\gamma}(-\pi+2\pi(x+1))
\end{equation}
for $x \in \C$. Let
$$
x_{n} = \frac{n}{r} + \frac{1}{2r}, \hspace{.2in} n \in \Z.
$$
Then $x \mapsto S_{\gamma}(-\pi+2\pi x)$ is analytic on
$\C \sm \{ x_{n} | n =r,r+1,\ldots \,\}$. If $m$ is a positive
integer then $\{ x_{n} | n = mr, mr+1,\ldots, (m+1)r-1 \,\}$
are poles of order $m$, while the points
$\{ x_{n} | n = -mr, -mr+1,\ldots, -mr+r-1 \,\}$ are zeros of
order $m$.

Let us use the function $S_{\gamma}$ to give another expression
for $\tau_{r}(M_{p/q})$. By \reflem{lem:functional} we have
that
$$
\prod_{l=0}^{m} (1-\xi^{k \pm l}) = \prod_{l=0}^{m}
\frac{S_{\gamma}(-\pi+2\gamma(k\pm l) -\gamma)}
     {S_{\gamma}(-\pi+2\gamma(k\pm l) +\gamma)}.
$$
Therefore
\begin{eqnarray*}
\prod_{l=0}^{m} (1-\xi^{k - l}) &=&
\frac{S_{\gamma}(-\pi+2\gamma(k - m) - \gamma)}
     {S_{\gamma}(-\pi+2\gamma k +\gamma)}, \\
\prod_{l=0}^{m} (1-\xi^{k + l}) &=&
\frac{S_{\gamma}(-\pi+2\gamma k -\gamma)}
     {S_{\gamma}(-\pi+2\gamma(k + m) +\gamma)}.
\end{eqnarray*}
So
$$
\prod_{l=0}^{m} (1-\xi^{k - l})(1-\xi^{k + l}) =
(1-\xi^{k})\frac{S_{\gamma}(-\pi+2\gamma(k - m) -\gamma)}
                {S_{\gamma}(-\pi+2\gamma(k + m) +\gamma)},
$$
and then by (\ref{eq:inv2})
$$
\tau_{r}(M_{p/q}) = \beta(r)\qsum \sum_{k=1}^{r-1}
\sum_{m=0}^{r-1} f_{n,r}\left( \frac{k}{r},\frac{m+1/2}{r}
\right),
$$
where
$$
f_{n,r}(x,y) = \sin\left(\frac{\pi}{q}(x-2nd)\right)
e^{2\pi i r \left( \frac{dn^{2}}{q} + \frac{p}{4q} x^{2}
                    -\frac{n}{q}x-xy\right)}
\frac{S_{\gamma}(-\pi+2\pi(x - y))}
     {S_{\gamma}(-\pi+2\pi(x + y))}
$$
and
\begin{eqnarray}\label{eq:betaconstant}
\beta(r) &=& -\frac{ia(r)}{2\sin(\pi/r)} \nonumber \\*
&=& \frac{i\sign(q)}{r\sqrt{|q|}}
    e^{-\frac{3\pi i}{4}\sign(pq)}
\exp\left(\frac{\pi i}{2r}\left[3\sign(pq)
   -\frac{p}{q}+\dS \left(\frac{p}{q}\right)\right]\right).
\end{eqnarray}
Note that $d$ is equal to the inverse of $p \bmod{q}$ and that
the functions $f_{n,r}$ are independent of the choice of this
inverse. By the remarks following \reflem{lem:functional} the
functions $f_{n,r}$ are holomorphic on
$\Omega_{r} \times \Omega_{r}$, where
\begin{equation}\label{eq:typicaldomain}
\Omega_{s} = \left\{ w \in \C \mid -\frac{1}{4s}
< \rea(w) < 1 + \frac{1}{4s} \; \right\}
\end{equation}
for $s \in ]0,\infty]$. By the residue theorem we therefore
end up with

\begin{lemma}\label{lem:contour}
The quantum invariants of $M_{p/q}$ are given by
$$
\tau_{r}(M_{p/q}) = \frac{\beta(r) r^{2}}{4} \qsum
  \int_{C_{r}^{1}} \cot(\pi r x)
   \left( \int_{C_{r}^{2}} \tan(\pi ry)
            f_{n,r}( x,y) \dte y \right) \dte x,
$$
where $\beta(r)$ is given by (\ref{eq:betaconstant}) and
$$
f_{n,r}(x,y) = \sin\left(\frac{\pi}{q}(x-2nd)\right)
e^{2\pi i r \left( \frac{dn^{2}}{q} + \frac{p}{4q} x^{2}
            -\frac{n}{q}x-xy\right)}
\frac{S_{\gamma}(-\pi+2\pi(x - y))}
{S_{\gamma}(-\pi+2\pi(x + y))},
$$
and where $C_{r}^{1}$ is a closed curve in $\Omega_{r}$ such
that the poles $\{ k/r \; | \; k = 1,2,\ldots,r-1 \; \}$ for
$x \mapsto \cot(\pi r x)$ lies inside $C_{r}^{1}$ and the poles
$0$ and $1$ lies outside $C_{r}^{1}$, and $C_{r}^{2}$ is a
closed curve in $\Omega_{r}$ such that the poles
$\{ (m+1/2)/r \; | \; m=0,1,\ldots,r-1 \; \}$ for
$y \mapsto \tan(\pi r y)$ lies inside $C_{r}^{2}$. Both curves
are oriented in the anti-clockwise direction.
\end{lemma}

Using the function $S_{\gamma}$ we can also express $J_{K}'(r)$
as a contour integral. By \reflem{lem:functional} we get that
$$
J_{K}'(r) = r \sum_{m=0}^{r-1}
\frac{S_{\gamma}(\pi-(2m+1)\gamma)}
     {S_{\gamma}(-\pi+(2m+1)\gamma)}.
$$
We have here used that
$$
\frac{S_{\gamma}(-\pi+\gamma)}{S_{\gamma}(\pi-\gamma)} =
\prod_{j=1}^{r-1}\frac{S_{\gamma}(\pi-(2j+1)\gamma)}
                      {S_{\gamma}(-\pi+(2j+1)\gamma)} =
\prod_{j=1}^{r-1}\left(1-e^{-2\pi\I\frac{j}{r}}\right) = r.
$$
If we put
$$
g_{r}(x) = \frac{S_{\gamma}(\pi-2\pi x)}
                {S_{\gamma}(-\pi + 2\pi x)}
$$
for $x \in \Omega_{\frac{1}{2}r}$ (see
(\ref{eq:typicaldomain})) we get
$$
J_{K}'(r) = r \sum_{m=0}^{r-1} g_{r}
\left(\frac{m+1/2}{r}\right),
$$
and we can write this sum as the single contour integral
\begin{equation}\label{eq:contourfigure8}
J_{K}'(r) = \frac{\I r^{2}}{2} \int_{C_{r}^{2}}
                        \tan(\pi rx) g_{r}( x) \dte x,
\end{equation}
where $C_{r}^{2}$ is given as in \reflem{lem:contour}.

\section{The large $r$ asymptotics of $J'_K(r)$ and
$\tau_{r}(M_{p/q})$}

\noindent In this section we investigate the large $r$
asymptotics of $\tau_{r}(M_{p/q})$ or more precisely the
leading term of this asymptotics, using the saddle point
method. We begin by calculating the large $r$ asymptotics of
$J_{K}'(r)$ using the expression (\ref{eq:contourfigure8}).
This calculation will demonstrate the use of the saddle point
method and will serve as a warm up for the more difficult
considerations of the asymptotics of $\tau_{r}(M_{p/q})$ in the
final part of this section.

\subsection{Semiclassical asymptotics of the quantum
dilogarithm}

\noindent It is well-known that the semiclassical asymptotics,
i.e.\ the small $\gamma$ asymptotics of the quantum dilogarithm
$S_\gamma$, is given by Euler's dilogarithm
\begin{equation}\label{eq:dilog}
\Li_{2}(z) = -\int_{0}^{z} \frac{\Log (1-w)}{w} \dte w
\end{equation}
for $z \in \C \sm ]1,\infty[$. Here and elsewhere $\Log$
denotes the principal logarithm. For $\zeta \in \C$ satisfying
either $\rea(\zeta) = \pm \pi$ and $\ima(\zeta) \geq 0$ or
$|\rea(\zeta)| < \pi$ one can check (see Appendix A) that
$$
\frac{1}{2\I \gamma} \Li _{2} (-e^{\I \zeta} ) =
\frac{1}{4} \int_{C_{R}}
\frac{e^{\zeta z}}{\sinh(\pi z) \gamma z^{2}} \dte z,
$$
hence we have that
\begin{equation}\label{eq:Sgdilog}
S_\gamma(\zeta) = \exp\left( \frac{1}{2\I \gamma}
      \Li _{2} (-e^{\I \zeta} ) + I_\gamma(\zeta)\right)
\end{equation}
for such $\zeta$, where
$$
I_{\gamma}(\zeta) = \frac{1}{4} \int_{C_{R}}
      \frac{e^{\zeta z}}{z\sinh(\pi z)}
      \left( \frac{1}{\sinh(\gamma z)}
              - \frac{1}{\gamma z}\right) \dte z.
$$

\begin{lemma}\label{lem:Sapproximation}
 If $|\rea(\zeta)| < \pi$ then
$$
|I_{\gamma}(\zeta)| \leq A \left( \frac{1}{\pi -\rea(\zeta)}
          + \frac{1}{\pi + \rea(\zeta)} \right) \gamma
          + B\left(1+e^{-\ima(\zeta) R} \right)\gamma,
$$
and for $|\rea(\zeta)| \leq \pi$ we have
$$
|I_{\gamma}(\zeta)| \leq 2A +
          B \left(1+e^{-\ima(\zeta) R} \right)\gamma,
$$
where $A$ and $B$ are positive constants only depending on $R$.
\end{lemma}

A proof is given in Appendix A. On the unit circle the
imaginary part of the dilogarithm is given by Clausen's
function $\Cl_{2}$, i.e.\
\begin{equation}\label{eq:dilogimacircle}
\ima\left( \Li_{2}( e^{i\theta})\right) = \Cl_{2}(\theta)
   = \sum_{n=1}^{\infty} \frac{\sin(n\theta)}{n^{2}}
   = -\int_{0}^{\theta}
         \Log\left|2\sin\left(\frac{t}{2}\right)\right|\dte t
\end{equation}
for $\theta \in \R$. One sees that $\Cl_{2}$ is increasing on
$[0,\pi/3] \cup [5\pi/3,2\pi]$ and decreasing on
$[\pi/3,5\pi/3]$. In particular, $\Cl_{2}$ attains its maximum
value at $\pi/3$ and its minimum value at $5\pi/3$. Moreover
\begin{equation}\label{eq:volume}
-\Cl_{2} \left( \frac{5\pi}{3}\right) =
                    \Cl_{2}\left(\frac{\pi}{3}\right)
  = -2\int_{0}^{\pi/6} \Log|2\sin(\phi)|\dte \phi
  = 2\J \left(\frac{\pi}{6}\right) = \frac{1}{2}\Vol(K),
\end{equation}
where $\J$ is Lobachevsky's function and $\Vol(K)$ is the
hyperbolic volume of the complement of the figure $8$ knot,
see e.g.\ \cite[Sec.~10.4]{[25]}.

\subsection{The large $r$ asymptotics of $J_{K}'(r)$}

\noindent We calculate the leading term of the large $r$
asymptotics of $J_{K}'(r)$, using the saddle point
method like Kashaev \cite{[14]}. Our calculation supplements
the calculation of Kashaev with the required analytic error
estimates. Let
\begin{eqnarray*}
C_{r}^{2} = C(\vep) &=& [\I+\vep,-\I+\vep]
+ [-\I+\vep,1-\vep-\I] \\
&& + [1-\vep-\I,1-\vep+\I] + [1-\vep+\I,\vep+\I],
\end{eqnarray*}
where $\vep \in ]0, \frac{1}{4r}[$. We let $C_{+}(\vep)$ be the
part of the contour $C(\vep)$ above the real axes and
$C_{-}(\vep)$ the part below the real axes. By
(\ref{eq:contourfigure8}) we have
$$
J_{K}'(r) = \frac{\I r^{2}}{2}
\left( J_{+}(r,\vep) + J_{-}(r,\vep) \right),
$$
where
$$
J_{\pm}(r,\vep) =  \int_{C_{\pm}(\vep)} \tan(\pi r x)g_{r}(x)
\dte x.
$$
Away from the real axis, the factor $\tan(\pi r x)$ can be
approximated by $\pm \I$ depending on whether we are in the
upper or lower half-plane. In fact we have
\begin{equation}\label{eq:tan+}
|\tan(\pi r x) - \I| \leq \left\{ \begin{array}{ll}
   4e^{-2\pi r \ima(x)}, & \ima(x) \geq \frac{1}{\pi r}, \\
   2e^{-2\pi r \ima(x)}, & r\rea(x) \in \Z, \ima(x) \geq 0
     \end{array}\right.
\end{equation}
and
\begin{equation}\label{eq:tan-}
|\tan(\pi r x) + \I| \leq \left\{ \begin{array}{ll}
  4e^{2\pi r \ima(x)}, & \ima(x) \leq -\frac{1}{\pi r}, \\
  2e^{2\pi r \ima(x)}, & r\rea(x) \in \Z, \ima(x) \leq 0.
    \end{array}\right.
\end{equation}
Therefore we write
$$
J_\pm(r,\vep) = \pm\I \int_{C_\pm(\vep)}g_r(x) \dte x +
\int_{C_\pm(\vep)}(\tan(\pi r x) \mp \I)g_r(x) \dte x.
$$
The estimates on $\tan(\pi r x)\pm \I$ can be used (see
Appendix B) to prove that
\begin{equation}\label{eq:figure8estimate1}
\left| \sum_{\mu =\pm 1}
    \int_{C_{\mu}(\vep)}(\tan(\pi r x) -\mu \I)g_r(x) \dte
       x\right| \leq K_{1} \frac{1}{r},
\end{equation}
where $K_{1}$ is a constant independent of $r$ and $\vep$.
Let now
\begin{equation}\label{eq:figure8Phase}
\Phi(x) = \frac{1}{2\pi\I}
\left( \Li_{2}(e^{-2\pi\I x}) - \Li_{2}(e^{2\pi\I x})\right).
\end{equation}
Note that $\Phi$ is analytic on
$D=\C \sm \{ \, x \in \C \, |\, \rea(x) \in \Z \,\}$ but not in
the points $\Z$, so here we see the reason for using the small
deformation parameter $\vep$. We have
\begin{eqnarray*}
\int_{C_\pm(\vep)}g_r(x) \dte x &=&
\int_{C_\pm(\vep)} e^{r\Phi(x)} \dte x \\
&& + \int_{C_\pm(\vep)}
\left(\exp\left( I_{\gamma}(\pi-2\pi x)
  - I_{\gamma}(-\pi + 2\pi x)\right) -1 \right)
e^{r\Phi(x)} \dte x.
\end{eqnarray*}
However, as we will see in Appendix B, the estimate in
\reflem{lem:Sapproximation} implies that
\begin{equation}\label{eq:figure8estimate2}
\left| \int_{C_{\mu}(\vep)}
   \left(\exp\left( I_{\gamma}(\pi-2\pi x)
           - I_{\gamma}(-\pi + 2\pi x)\right) -1\right)
   e^{r\Phi(x)} \dte x\right| \leq \frac{K_{2}\Log(r)}{r}
   e^{\frac{r}{2\pi}\Vol (K)}
\end{equation}
for $\mu =\pm 1$, where $K_{2}$ is a constant independent of
$r$ and $\vep$. We will see below that the estimates
(\ref{eq:figure8estimate1}) and (\ref{eq:figure8estimate2})
imply that the leading order large $r$ asymptotics of $J'_K(r)$
is given by
\begin{equation}\label{eq:Asymf8}
J'_K(r) \sim \frac{r^{2}}{2}
  \left( \int_{C_-(\vep)} e^{r\Phi(x)} \dte x
       - \int_{C_+(\vep)} e^{r\Phi(x)} \dte x \right),
\end{equation}
to which we can apply the saddle point method, see e.g.\
\cite[Chap.~5]{[6]}. First we determine the stationary points
of the phase function $\Phi$. On $D$ we have
$$
\Phi'(x) = \Log\left(1-e^{2\pi\I x}\right) +
           \Log\left(1-e^{-2\pi\I x}\right).
$$
If we put $z=e^{2\pi\I x}$, then $\Phi'(x)=0$ implies that
\begin{equation}\label{eq:hyp}
z^{2}-z+1=0.
\end{equation}
The equation (\ref{eq:hyp}) has the solutions
$z_{\pm}=e^{\pm \I\pi/3}$. We have
$1 - z_{\pm} = 1/2 \mp i\sqrt{3}/2$
which both have norm $1$ and are each others conjugate, so
$$
\Log\left(1-z_{+} \right) + \Log\left(1-z_{-}\right) = 0.
$$
We note that $z_{\pm}=e^{\pm2\pi \I \frac{1}{6}}$ correspond to
the $x$--points $\pm 1/6 +\Z$. These points are non-degenerate
critical points. In fact,
$$
\Phi''(x) = 2\pi\I\frac{e^{2\pi\I x} + 1}{e^{2\pi\I x} - 1}
$$
on $D$, so in particular $\Phi''(x_{\pm})=\pm 2\pi \sqrt{3}$
for $x_{\pm} \in \pm 1/6 +\Z$. The imaginary part of $\Phi(x)$
is zero for $x \in \R$ and
$$
\Phi(x) = \frac{1}{2\pi}
\ima \left( \Li_{2}(e^{-2\pi\I x}) - \Li_{2}(e^{2\pi\I x})
\right) = -\frac{1}{\pi}\Cl_{2}(2\pi x),
$$
by (\ref{eq:dilogimacircle}). Let $x_{\pm} \in \pm 1/6 +\Z$.
By (\ref{eq:volume}) we have
$\Cl_{2}(2\pi x_{-})=-\Cl_{2}(2\pi x_{+})
= -\Cl_{2}(\pi/3)=-\Vol(K)/2$, i.e.\
$$
\Phi(x_{\pm}) = \mp \frac{1}{2\pi}\Vol(K).
$$
By Cauchy's theorem we have
$$
\int_{C_-(\vep)} e^{r\Phi(x)} \dte x
          - \int_{C_+(\vep)} e^{r\Phi(x)} \dte x
= 2\int_{C_-(\vep)} e^{r\Phi(x)} \dte x
= -2 \int_{C_+(\vep)} e^{r\Phi(x)} \dte x.
$$
Deform $C_{-}(\vep)$ to $[\vep,1-\vep]$ keeping the end points
fixed. This does not change the integral
$\int_{C_{-}(\vep)} e^{r\Phi(x)} \dte x$. Let $x_{0}=5/6$. By
terminology borrowed from \cite[Sec.~5.4]{[6]} the axis of the
saddle point $x_{0}$ is the real axis (i.e.\ the directions of
steepest descent are along the real axis). From the analysis of
\cite[Sec.~5.7]{[6]} it follows that we can find a $\delta>0$
(independent of $r$ and $\vep$) such that
$[x_{0}-\delta,x_{0}+\delta] \subseteq [1/(4r),1-1/(4r)]$ and
such that we have an asymptotic expansion
$$
\int_{-\delta}^{\delta} e^{r\Phi(x_{0} + t)} \dte t
   \sim \frac{1}{3^{1/4}\sqrt{r}} e^{\frac{r}{2\pi}\Vol(K)}
   \left(1+\sum_{n=1}^{\infty} d_{n} r^{-n} \right)
$$
in the limit $r \ria \infty$, where the $d_{n}$'s are certain
complex numbers. Finally we note that
$$
\left|\int_{\vep}^{x_{0}-\delta} e^{r\Phi(t)}\dte t +
\int_{x_{0}+\delta}^{1-\vep} e^{r\Phi(t)}\dte t \right|
\leq e^{rc},
$$
where
$c = \max \{ -\Cl_{2}(2\pi(x_{0}-\delta))/\pi,
                 -Cl_{2}(2\pi(x_{0}+\delta))/\pi \}
   < \Vol(K)/2\pi$,
see above (\ref{eq:volume}). We have shown

\begin{lemma}\label{volume}
The leading order large $r$ asymptotics of $J_{K}'(r)$ is given
by
\begin{equation}\label{eq:asymptoticsK}
J_{K}'(r) \sim 3^{-1/4} r^{3/2}
       \exp\left(\frac{r}{2\pi}\Vol(K)\right).
\end{equation}
In fact
$$
J_{K}'(r) = 3^{-1/4} r^{3/2}
    \exp\left(\frac{r}{2\pi}\Vol(K)\right)
 + \Os \left( r\Log(r)
         \exp\left(\frac{r}{2\pi}\Vol(K)\right)\right)
$$
in the limit $r \ria \infty$.\HS
\end{lemma}

In particular
$$
\lim_{r \ria \infty} \frac{2\pi\Log(J_{K}'(r))}{r}
= \Vol(K)
$$
as predicted by the volume conjecture of Kashaev \cite{[14]}
and Murakami, Murakami \cite{[24]} and as proven by Ekholm and
others, see \cite{[22]}. However, the arguments of Ekholm and
others can't see the finer details of the asymptotic behaviour
(\ref{eq:asymptoticsK}), namely the polynomial part
$3^{-1/4} r^{3/2}$.

\begin{remark}\label{rem:hyperbolicity}
The complement $S^{3} \sm K$ of the figure $8$ knot can be
decomposed into two ideal hyperbolic tetrahedra each
parametrized by a certain complex number. This decomposition
then defines a hyperbolic  structure on this complement if and
only if a certain set of conditions is satisfied. In fact, if
the complex parameters for the two tetrahedra are respectively
$a$ and $b$, then these conditions are equivalent to $a=b$ and
$b^{2}-b+1=0$, which is equation (\ref{eq:hyp}) after
substituting $b$ for $z$. We refer to \cite[Sec.~3]{[23]} for
more details. This phenomenon, that one finds the hyperbolicity
equation for the knot complement as the equation for the
stationary points of the phase function, seems to be a
general principal for hyperbolic knots as argued by Thurston
and Yokota, cf.\ \cite{[32]}, \cite{[35]}. However, there are
major unsolved analytic difficulties in their approach.
Basically they conjecture that one can carry out an asymptotic
analysis similar to the one we carried out above for the figure
$8$ knot. To prove their conjecture one first has to show how
to give an exact (multi-dimensional) contour integral formula
for the Jones polynomial of a hyperbolic knot like our
(\ref{eq:contourfigure8}). A main part consists of choosing a
correct (multi-dimensional) contour. Secondly, one has to carry
out an asymptotic analysis similar to the one leading to
(\ref{eq:Asymf8}). This analysis is relatively simple for the
figure $8$ knot due to the fact that we have a single
(one-dimensional) contour in this case. In general one gets a
contour of dimension $>1$ and the asymptotic analysis is
expected to be harder (as also illustrated by the asymptotic
analysis of the double-contour integral expression for the
invariant $\bar{\tau}_{r}(M_{p/q})$ in
\reflem{lem:contour}, see next section.)
\end{remark}

\subsection{The large $r$ asymptotics of $\tau_{r}(M_{p/q})$}

\noindent In Sec.~5 we will see that the signs of the phases
in the asymptotics of
$\bar{\tau}_{r}(M_{p/q}) = \overline{\tau_{r}(M_{p/q})}$
agrees with the Chern--Simons values, hence we work with
this conjugate invariant. Because of (\ref{eq:mirror2})
we can always obtain the asymptotic expansion of
$\tau_{r}(M_{p/q})$ by complex conjugation or by replacing
either $p$ by $-p$ or $q$ by $-q$. By \reflem{lem:contour} we
have
\begin{equation}\label{eq:barinvariant}
\bar{\tau}_{r}(M_{p/q}) = \beta_{1}(r) \qsum
  \int_{C_{r}^{1}\times C_r^2} \cot(\pi r x) \tan(\pi ry)
        \bar{f}_{n,r}(x,y) \dte x  \dte y,
\end{equation}
where
\begin{equation}\label{eq:beta1constant}
\beta_{1}(r) = \frac{i\sign(q)r}{4\sqrt{|q|}}
         e^{\frac{3\pi i}{4}\sign(pq)}
         \exp\left(-\frac{\pi i}{2r}\left[3\sign(pq)
      -\frac{p}{q}+\dS\left(\frac{p}{q}\right)\right]\right)
\end{equation}
and
$$
\bar{f}_{n,r}(x,y)=\sin\left(\frac{\pi}{q}(x-2nd)\right)
      e^{2\pi i r \left( -\frac{dn^{2}}{q}-\frac{p}{4q} x^{2}
             + \frac{n}{q}x-xy\right)}
\frac{S_{\gamma}(-\pi+2\pi(x - y))}
     {S_{\gamma}(-\pi+2\pi(x + y))},
$$
where $\gamma=\pi/r$ as usual. Let for $k,l \in \Z$
$$
\Omega_{k,l} = \{ \; (x,y) \in \C^{2} \mid
   \rea(x)+\rea(y) \in [k,k+1],
   \hspace{.1in}\rea(x)-\rea(y) \in [-l,-l+1]\;\}.
$$
For $(x,y)\in \Omega_{k,l}$, we have by (\ref{eq:functional})
and (\ref{eq:Sgdilog}) that
\begin{eqnarray*}
\bar{f}_{n,r}(x,y) &=& \sin\left(\frac{\pi}{q}( x-2nd)\right)
  \left( 1+e^{2\pi i(x-y)r} \right)^{l}
  \left( 1+e^{2\pi i(x+y)r} \right)^{k}
e^{2\pi ir\Phi_{n}(x,y)} \\*
&& \hspace{.2in} \times \exp\left(I_{\gamma}(-\pi+2\pi(x-y+l))
                          -I_{\gamma}(-\pi+2\pi(x+y-k))\right),
\end{eqnarray*}
where
$$
\Phi_{n}(x,y) = -\frac{dn^{2}}{q} - \frac{p}{4q} x^{2}
+ \frac{n}{q}x - xy + \frac{1}{4\pi^{2}}
\left( \Li_{2}(e^{2\pi i(x+y)}) - \Li_{2}(e^{2\pi i (x-y)})
\right).
$$
We note that this expression and the above expression for
$\bar{f}_{n,r}$ are only valid for $(x,y) \in \C^{2}$
satisfying the two conditions
\begin{eqnarray*}
\mbox{\rm (i)}&& \rea(x) + \rea(y) \notin \Z \vee
       \ima(x)+\ima(y) \geq 0, \\*
\mbox{\rm (ii)}&& \rea(x) - \rea(y) \notin \Z \vee
       \ima(x)-\ima(y) \geq 0.
\end{eqnarray*}

Observe that for $k,l \in \{0,1\},$ which corresponds to the
four different $\Omega_{k,l}$ intersecting
$C_r^1\times C_r^2$, we have that
$$
\left( 1+e^{2\pi i(x-y)r} \right)^{l}
\left( 1+e^{2\pi i(x+y)r} \right)^{k}
 = \sum_{(a,b)\in F_{k,l}} e^{2\pi i(a(x+y)+b(x-y))r},
$$
where
$F_{k,l}= \{(a,b)\in \{0,1\}^{2} \mid a \leq k,
\quad b \leq l\}$.
Hence, for such $k,l$ we have
\begin{eqnarray*}
\bar{f}_{n,r}(x,y) &=& \sum_{(a,b)\in F_{k,l}}
           \sin\left(\frac{\pi}{q}( x-2nd)\right)
           e^{2\pi ir\Phi^{a,b}_{n}(x,y)} \\*
&& \times \exp\left(I_{\gamma}(-\pi+2\pi(x-y+l)) -
                    I_{\gamma}(-\pi+2\pi(x+y-k))\right)
\end{eqnarray*}
for $(x,y) \in \Omega_{k,l}$, where
\begin{equation}\label{eq:Phasefunctiontypical}
\Phi_{n}^{a,b}(x,y)= a(x+y) + b(x-y) + \Phi_{n}(x,y).
\end{equation}
Let
$$
\Omega_{k,l}^{\mu,\nu} = \{ (x,y) \in \Omega_{k,l} \mid
\mu \ima(x) \geq 0, \quad \nu \ima(y) \geq 0\}.
$$

\begin{conjecture}\label{conj:mainconj}
There exists surfaces
$\Sigma_{k,l,a,b}^{\mu,\nu,n} \subset \Omega_{k,l}^{\mu,\nu}$
for $(k,l) \in \{0,1\}^{2}$, $(a,b) \in F_{k,l}$,
$(\mu,\nu) \in \{ \pm 1 \}^{2}$ and $n \in \Z / |q|\Z$ such
that the leading order large $r$ asymptotics of the quantum
invariant is given by
\begin{eqnarray}\label{asympconjf}
\bar{\tau}_{r}(M_{p/q}) &\sim&
\frac{i\sign(q)r}{4\sqrt{|q|}}e^{\frac{3\pi i}{4}\sign(pq)}
\qsum \sum_{(k,l) \in \{0,1\}^{2}} \sum_{(a,b)\in F_{k,l}}
\sum_{(\mu,\nu)\in \{\pm 1\}^{2}} \mu\nu \nonumber \\*
&& \times \int_{\Sigma_{k,l,a,b}^{\mu,\nu,n}}
          \sin\left(\frac{\pi}{q}( x-2nd)\right)
          e^{2\pi ir\Phi_{n}^{a,b}(x,y)} \dte x \dte y.
\end{eqnarray}
Moreover, if $p/q \neq 0$ the surfaces
$\Sigma_{k,l,a,b}^{\mu,\nu,n} \subset \Omega_{k,l}^{\mu,\nu}$
can be chosen such that any critical point of $\Phi_{n}^{a,b}$
belonging to $\Sigma_{k,l,a,b}^{\mu,\nu,n}$ also belongs to the
set $\mS$ in (\ref{eq:SUcriticalpoint}).
\end{conjecture}

The rational behind this conjecture is that we anticipate an
analysis of the expression (\ref{eq:barinvariant}) paralleling
our analysis of $J'_K(r)$ should be applicable. I.e.\ $\tan$
and $\cot$ should be approximated by $\pm \I$ depending on the
signs of $\ima(y)$ and $\ima(x)$ and $\bar{f}_{n,r}(x,y)$ by
$\sum_{(a,b) \in F_{k,l}}
\sin\left(\frac{\pi}{q}( x-2nd)\right)
e^{2\pi ir\Phi_{n}^{a,b}(x,y)}$
for some deformation of the part of $C_r^1\times C_r^2$ which
is contained in $\Omega_{k,l}$. We have partial analytic
results supporting this conjecture. The factor in front of the
sum is simply the leading term in the large $r$ asymptotics of
$\beta_{1}(r)$, see (\ref{eq:beta1constant}). In case $p/q = 0$
our results in Appendix C show that we have to include critical
points not belonging to $\mS$ in our surfaces
$\Sigma_{k,l,a,b}^{\mu,\nu,n}$.

In the remaining part of this section we will compute the
large $r$ asymptotics of the right-hand side of
(\ref{asympconjf}). The only task left is to calculate the
leading term in the large $r$ asymptotic of the integrals
\begin{equation}\label{eq:integraltypical}
I = \int_{ \Sigma_{k,l,a,b}^{\mu,\nu,n} }
   \sin\left(\frac{\pi}{q}( x-2nd)\right)
       e^{2\pi ir\Phi_{n}^{a,b}(x,y)} \dte x \dte y,
\end{equation}
where $k,l \in \{0,1\}$, $(a,b) \in F_{k,l}$,
$\mu,\nu \in \{ \pm 1 \}$ and $n \in \Z$. By the properties of
the surfaces $\Sigma_{k,l,a,b}^{\mu,\nu,n}$ postulated in
\refconj{conj:mainconj} this leading asymptotics should be
calculable by the saddle point method. Let us give some
details. Assume that $(x_0,y_0)$ is a critical point of
$\Phi=\Phi_{n}^{a,b}$ belonging to the surface
$\Sigma=\Sigma_{k,l,a,b}^{\mu,\nu,n}$ and assume that
$$
\ima\left( (\Phi(x,y) - \Phi(x_{0},y_{0}))\right) \geq 0
$$
for all $(x,y)$ in a small neighborhood of $(x_0,y_0)$ in
$\Sigma$ with equality only in $(x_0,y_0)$. In this case we say
that $\Sigma$ pass through the critical point $(x_0,y_0)$ in
the directions of steepest descent. Moreover, assume this is
satisfied for all critical points of $\Phi$ inside $\Sigma$.
Then the main contribution to the integral $I$ in the large $r$
limit is given by integrating over small neighborhoods of the
critical points in $\Sigma$.

Let us therefore begin by computing the stationary points of
$\Phi_{n}^{a,b}$. To this end, it is more convenient to work
with the functions
\begin{equation}\label{eq:Phasefunctiontypical1}
\Psi^{a,b}_{n}(x,y) = ax + by + \Phi_{n}(x,y),
\end{equation}
$a,b,n \in \Z$, so $\Phi_{n}^{a,b}=\Psi_{n}^{a+b,a-b}$.

Let $a,b,n \in \Z$ and put $\Psi=\Psi_{n}^{a,b}$. Let
$z=e^{2\pi i x}$ and $w=e^{2\pi iy}$. Then
\begin{eqnarray}\label{eq:critical}
&&2\pi i\frac{\partial \Psi}{\partial x}(x,y)
= 2\pi i(a-y) - \frac{p}{2q} 2\pi i x
+ \frac{2\pi i n}{q} + \Log(1-zw)- \Log(1-zw^{-1}),
\nonumber \\*
&&2\pi i\frac{\partial \Psi}{\partial y}(x,y)
= 2\pi i (b-x)  + \Log(1-zw) + \Log(1-zw^{-1}),
\end{eqnarray}
where we have to assume (which we will also assume in what
follows) that both $zw$ and $zw^{-1}$ are different from $1$.

We will need to specify a certain square root of $z$, namely
let $v=e^{\pi ix}$. Then
$\frac{\partial \Psi}{\partial x}(x,y)=0$ implies that
\begin{equation}\label{eq:stationary1}
v^{-p} = \left(\frac{w-v^{2}}{1-v^{2}w}\right)^{q},
\end{equation}
and
$\frac{\partial \Psi}{\partial y}(x,y) = 0$ implies that
\begin{equation}\label{eq:stationary2}
(1-v^{2}w)(w-v^{2}) = v^{2}w.
\end{equation}
Both equations (\ref{eq:stationary1}) and
(\ref{eq:stationary2}) are independent of $a$, $b$, and $n$.
We note that $(v,w)=(0,0)$ is the only solution to
(\ref{eq:stationary2}) with $v$ or $w$ equal to zero. Note,
moreover, that $(v,w)$ is a common solution to
(\ref{eq:stationary1}) and (\ref{eq:stationary2}) if and only
if $(\bar{v},\bar{w})$ is a common solution to these two
equations. By writing (\ref{eq:stationary1}) as
$$
(-1)^{q}v^{p+2q} =
\left( \frac{1-v^{2}w}{1-v^{-2}w} \right)^{q}
$$
and (\ref{eq:stationary2}) as
$$
(1-v^{2}w)(1-v^{-2}w) = -w
$$
we see that $(v,w)$ is a nonzero solution to
(\ref{eq:stationary1}) and (\ref{eq:stationary2}) if and only
if $(v^{-1},w)$ is such a solution. If $p$ is even, then
$(v,w)$ is a solution to (\ref{eq:stationary1}) and
(\ref{eq:stationary2}) if and only if $(-v,w)$ is such a
solution.

Let us oppositely begin with a common solution
$(v,w) \in \C^{*} \times \C^{*}$ to (\ref{eq:stationary1}) and
(\ref{eq:stationary2}), where, as usual,
$\C^{*} = \C \sm \{0\}$. Write $v = e^{\pi i x}$ and
$w = e^{2\pi iy}$ with $\rea(x) \in ]-1,1]$ and
$\rea(y) \in ]-1/2,1/2]$, i.e.\ $x = \frac{1}{\pi i}\Log(v)$
and $y = \frac{1}{2\pi i}\Log(w)$. One now easily deduce from
(\ref{eq:stationary2}) that there exists a unique $b \in \Z$
such that
\begin{equation}\label{eq:critical2}
0 = 2\pi i (b-x)  + \Log\left(1-e^{2\pi i(x+y)}\right)
                  + \Log\left(1-e^{2\pi i(x-y)}\right),
\end{equation}
and from (\ref{eq:stationary1}) we deduce that there exists a
unique $n \in \Z$ such that
\begin{equation}\label{eq:critical1}
0 = -2\pi iy - \frac{p}{q} \pi i x + \frac{2\pi i n}{q}
    + \Log\left(1-e^{2\pi i(x+y)}\right)
    - \Log\left(1-e^{2\pi i(x-y)}\right).
\end{equation}
That is, there exists a unique pair of integers $b,n$ such
that $(x,y)$ is a stationary point of $\Psi_{n}^{0,b}$.

\begin{remark}
Let us make a slight digression by giving some general remarks
about the set of solutions to (\ref{eq:stationary2}). Assume
that $v,w \in \C^{*}$ and let $z=v^{2}$. Then
(\ref{eq:stationary2}) can be written in the following two
ways
\begin{eqnarray}\label{eq:stationary2no2}
&& z^{2} - \left( w + \frac{1}{w} + 1 \right)z + 1 = 0,
\nonumber \\*
&& w^{2} - \left( z + \frac{1}{z} - 1 \right)w + 1 = 0.
\end{eqnarray}
It is straightforward to see that if $(z,w)$ is a solution to
(\ref{eq:stationary2no2}) with $w \in \R \sm \{0\}$, then
$z$ is real and positive if $w>0$, $z$ is real and negative if
$w \in ]-\infty,-(3+\sqrt{5})/2] \cup [-(3-\sqrt{5})/2,0[$,
and $z \in S^{1}$ if
$$
-(3 + \sqrt{5})/2 \leq  w \leq - (3 - \sqrt{5})/2.
$$
If $(z,w)$ is a solution to (\ref{eq:stationary2no2}) with
$z \in \R \sm \{0\}$, then $w$ is real and negative if $z<0$,
$w$ is real and positive if
$z \in ]0,(3-\sqrt{5})/2] \cup [(3+\sqrt{5})/2,\infty[$,
and $w \in S^{1}$ if $z \in [(3-\sqrt{5})/2,(3+\sqrt{5})/2]$.

Later on we will be particularly interested in common
solutions to (\ref{eq:stationary1}) and (\ref{eq:stationary2})
with $v \in S^{1}$ and $w \in \R \sm \{0\}$. Assume that
$(v,w)$ is such a solution, and write $v=e^{\pi i x}$,
$x \in ]-1,1]$. In that case the 2nd of the equations in
(\ref{eq:stationary2no2}) is
\begin{equation}\label{eq:wrealzcircle}
w^{2} + (1-2\rea(z))w + 1 = 0.
\end{equation}
Since $z \in S^{1}$ we have $1-2\rea(z) \in [-1,3]$ and since
$w$ is real we also have $|1-2\rea(z)| \geq 2$, so
$1-2\rea(z) \in [2,3]$ or $\cos(2\pi x)=\rea(z) \in [-1,-1/2]$.
By (\ref{eq:wrealzcircle}) we find that
\begin{equation}\label{eq:wrealcircle2}
w = w_{\pm}(x) = \cos(2\pi x) - \frac{1}{2} \pm
             \sqrt{\cos^{2}(2\pi x)-\cos(2\pi x)-\frac{3}{4}}.
\end{equation}
If $w>0$ we have already seen that $z \in ]0,\infty[$. But
then $z=1$ so $w^{2} - w + 1 = 0$ by (\ref{eq:stationary2no2})
contradicting the fact that $w$ is real.

Taking absolute values we get from
(\ref{eq:stationary1}) that $|1 - zw|=|w - z|$ and then from
(\ref{eq:stationary2}) that
$|w|=|w-z|^{2}=(w-z)(w-\bar{z})=w^{2}-2\rea(z)w+1$ (still
assuming that $(v,w) \in S^1 \times \R$). Since $w<0$
this equation is equivalent to (\ref{eq:wrealzcircle}). In
other words we have to consider the argument equation following
from (\ref{eq:stationary1}) to obtain information from that
equation not obtainable from (\ref{eq:stationary2}).
\end{remark}

Let us now turn to the second derivative of $\Psi$ in a
critical point $(x_{0},y_{0})$, i.e.\ the Hessian
$H=H(x_{0},y_{0})$ of $\Psi$ in $(x_{0},y_{0})$ (which is
equal to the Hessian of $\Phi_{n}$ in $(x_{0},y_{0})$). Put
$(v_{0},w_{0})=(e^{\pi ix_{0}},e^{2\pi iy_{0}})$ and
$z_{0}=v_{0}^{2}$. By a small computation, using the fact that
$(v_{0},w_{0})$ is a solution to (\ref{eq:stationary2}), we
find that
\begin{eqnarray}\label{eq:Hessian}
H_{11} &=& \frac{\partial^{2}\Psi}{\partial x^{2}}(x_{0},y_{0})
= \frac{1}{w_{0}}-w_{0}-\frac{p}{2q}, \nonumber \\*
H_{12} = H_{21} &=&
\frac{\partial^{2}\Psi}{\partial x\partial y}(x_{0},y_{0})
= z_{0}-\frac{1}{z_{0}}, \nonumber \\*
H_{22} &=& \frac{\partial^{2}\Psi}{\partial y^{2}}(x_{0},y_{0})
= \frac{1}{w_{0}}-w_{0}.
\end{eqnarray}
Let us begin by establishing under which circumstances the
critical point $(x_{0},y_{0})$ is non-degenerate. By
(\ref{eq:Hessian})
$$
\det(H) = \left(w_{0} + \frac{1}{w_{0}}\right)^{2}
        - \left( z_{0}+\frac{1}{z_{0}} \right)^{2}
        + \frac{p}{2q}\left(w_{0} -\frac{1}{w_{0}}\right).
$$
By (\ref{eq:stationary2no2}) we have
$w_{0} + \frac{1}{w_{0}} = z_{0} + \frac{1}{z_{0}} - 1$ so
\begin{equation}\label{eq:detHessian3}
\det(H) = 1 - 2\left( z_{0} + \frac{1}{z_{0}} \right)
          + \frac{p}{2q}\left( w_{0} - \frac{1}{w_{0}} \right).
\end{equation}
If we introduce $y'_0 = iy_0$ we get the formula
\begin{equation}\label{eq:detHessian2}
\det(H) = 1 - 4\cos(2\pi x_0) + \frac{p}{2q} \sinh(2 \pi y_0').
\end{equation}
We are particularly interested in critical points where
$w_{0} < 0$ and $z_{0} \in S^{1}$ (or equivalently
$(x_0,y_0) \in \mS$, see (\ref{eq:SUcriticalpoint})). By
definition of $z_0$ we have
$z_{0} - 1/z_{0} = 2i\sin(2\pi x_0)$ and by (\ref{eq:Hessian})
we then get
\begin{equation}\label{eq:detHessian1}
\det(H) = \left( \frac{1-w_{0}^{2}}{w_{0}}
- \frac{p}{2q} \right)\left( \frac{1-w_{0}^{2}}{w_{0}}\right)
+ 4\sin^{2}(2\pi x_0).
\end{equation}
Thus if $w_{0} \in ]-1,0[$ and
$p/q \geq 0$ or if $w_{0} \in ]-\infty,-1[$ and
$p/q \leq 0$ we have $\det(H)>0$. If $w_{0}=-1$ we
also have $\det(H) >0$ by (\ref{eq:wrealcircle2}).
Since $\rea(y_0) \in \frac{1}{2} + \Z$ we get by
(\ref{eq:detHessian3}) the alternative formula
\begin{equation}\label{eq:detHessian5}
\det(H) = 1 - 4\cos(2\pi x_0)
+ \frac{p}{2q} \sinh\left(2 \pi \ima(y_0)\right).
\end{equation}
Finally we get by (\ref{eq:wrealcircle2}) and
(\ref{eq:detHessian3}) that
\begin{equation}\label{eq:detHessian4}
\det(H) = 1 - 4\cos(2\pi x_0) \pm \frac{p}{q}
\sqrt{\cos^2(2\pi x_0)-\cos(2\pi x_0)-\frac{3}{4}}.
\end{equation}
for $w=w_{\pm}(x_0)$.

\begin{proposition}\label{prop:Nondeg}
Let the situation be as above and assume that $(x_0,y_0)$ is
a critical point of $\Psi$ belonging to the
set $\mS$ in (\ref{eq:SUcriticalpoint}). If $|p/q| < \sqrt{20}$
then $(x_0,y_0)$ is non-degenerate. If $|p/q| > \sqrt{20}$ then
$(x_0,y_0)$ is non-degenerate except if
$$
\rea(z_0) = \cos(2\pi x_0) = \frac{1}{2}
+ \frac{4 - \left| \frac{p}{q} \right|
\sqrt{\frac{p^2}{q^2} - 15}}{\frac{p^2}{q^2} - 16}
$$
and $w_0 = w_+(x_0)$ if $p/q < 0$ and $w_0 = w_-(x_0)$ if
$p/q > 0$.
\end{proposition}

\begin{proof}
Let $z_0=e^{2\pi ix_0}$ and $w_0=e^{2\pi iy_0}$ as above and
put $c=\cos(2\pi x_0)$. Assume that $(x_0,y_0)$ is degenerate,
i.e.\ that $\det(H) = 0$. Then
$$
\left( 16 - \frac{p^2}{q^2} \right) c^{2}
+ \left( \frac{p^2}{q^2} - 8 \right) c
+ 1 + \frac{3}{4}\frac{p^2}{q^2} = 0
$$
by (\ref{eq:detHessian4}). If $p^2/q^2 < 15$ this equation does
not have real solutions, so $p^2/q^2 \geq 15$. If $p^2/q^2=16$
we find that $c=-13/8$ giving a contradiction.
Thus $p^2/q^2 \in [15,16[ \cup ]16,\infty[$ and
\begin{equation}\label{eq:c}
c = \frac{8 - \frac{p^2}{q^2} \pm 2\left| \frac{p}{q} \right|
                 \sqrt{ \frac{p^2}{q^2} - 15 } }
         {2 \left( 16 - \frac{p^2}{q^2} \right)}.
\end{equation}
By this and $c \in [-1,-1/2]$, see below
(\ref{eq:wrealzcircle}), we get that
\begin{equation}\label{eq:inequalities}
\frac{3}{2}\frac{p^2}{q^2} - 20 \leq
\pm \left| \frac{p}{q} \right|\sqrt{\frac{p^2}{q^2} - 15}
\leq \frac{p^2}{q^2} - 12
\end{equation}
for $p^2/q^2 \in [15,16[$ with the opposite inequalities
for $p^2/q^2 \in ]16,\infty[$.
For $p^2/q^2 \in [15,16[$ we have
$\frac{3}{2}\frac{p^2}{q^2} - 20 >0$ so we have a plus in
front of
$\left| \frac{p}{q} \right|\sqrt{\frac{p^2}{q^2} - 15}$
in (\ref{eq:inequalities}) and in that case the first
inequality in (\ref{eq:inequalities}) is equivalent to
\begin{equation}\label{eq:inequalities2}
\frac{5}{4} \left( \frac{p}{q} \right)^{4}
- 45\left( \frac{p}{q} \right)^{2} + 400 \leq 0,
\end{equation}
which forces $p^2/q^2 \in [16,20]$ giving a contradiction.

Therefore $p^2/q^2 >16$ and (\ref{eq:c}) and $c \in [-1,-1/2]$
leads, as already stated, to (\ref{eq:inequalities}) with the
opposite inequalities. Since $p^2/q^2 - 12 >0$ we again have
a plus in front of
$\left| \frac{p}{q} \right|\sqrt{\frac{p^2}{q^2} - 15}$, and
the second inequality is automatically satisfied for all
$p^2/q^2 \geq 16$. This time the first inequality leads to
(\ref{eq:inequalities2}) with the opposite inequalities, so
we can conclude that $p^2/q^2 \geq 20$ and in fact therefore
$|p/q| > \sqrt{20}$. The formula for $\cos(2\pi x_0)$
is just (\ref{eq:c}) with the plus sign in front of
$\left| \frac{p}{q} \right|\sqrt{ \frac{p^2}{q^2} - 15 }$.
\end{proof}

We note that the only other solutions
$(z,w) \in S^{1} \times ]-\infty,0[$ to (\ref{eq:wrealzcircle})
satisfying that $\rea(z)$ is equal to the right-hand side
of (\ref{eq:c}) are $(z_{0},1/w_{0})$, $(\bar{z_{0}},w_{0})$
and $(\bar{z_{0}},1/w_{0})$. By the remarks following
(\ref{eq:detHessian1}) we get that among these three points
only the point $(\bar{z_{0}},w_{0})$ can actually satisfy, that
the right-hand side of (\ref{eq:detHessian1}) is zero.
We conjecture that also in the case $|p/q|>\sqrt{20}$ all
critical points of $\Psi$ which belong to the set $\mS$ are
non-degenerate. To confirm this one has probably to use the
argument equation following from (\ref{eq:stationary1}). By a
direct check we have confirmed this conjecture in
the cases
$|p/q| \in \{14/3,5,16/3,6,20/3,22/3,8,26/3,28/3,10,32/3\}$.

We are now ready to calculate the contribution to the leading
order large $r$ asymptotics of the integral
(\ref{eq:integraltypical}) coming from a critical point $(x,y)$
of the phase function $\Phi_{n}^{a,b}$ belonging to
$\mS \cap \Sigma_{k,l,a,b}^{\mu,\nu,n}$.

The integral $I$ in (\ref{eq:integraltypical}) is a double
contour integral and as such is calcalculated by choosing a
contour for each variable $x$ and $y$. We thus have an
``inner'' integral and an ``outer'' integral and in general the
contour for the inner integral can depend on where we are on
the outer contour. Here we will, however, only consider the
very simple case where these two contours are independent of
each other.

Therefore let $\alpha,\beta \in S^{1}$ and let
$\gamma_{\alpha}, \gamma_{\beta} : I_{\delta} \to \C$ be given
by $\gamma_{\alpha}(s) = x_{0} + \alpha s$,
$\gamma_{\beta}(s) = y_{0}+\beta s$, where
$I_{\delta} = [-\delta,\delta]$ for a sufficiently small
$\delta>0$. According to the saddle point method the main
contributions to the integral (\ref{eq:integraltypical}) comes
from integrating over a small neighborhood of
$\mS \cap \Sigma$ in $\Sigma=\Sigma_{k,l,a,b}^{\mu,\nu,n}$.
We are therefore lead to consider an integral of the form
\begin{eqnarray*}
K(x_{0},y_{0}) &=& \int_{\gamma_{\alpha}}
    \left( \int_{\gamma_{\beta}}
    \sin\left(\frac{\pi}{q}(x-2nd)\right)
    e^{2\pi ir\Psi(x,y)} \dte y \right)\dte x \\*
&=& \alpha\beta\int_{-\delta}^{\delta}
    \left( \int_{-\delta}^{\delta}
    \sin\left(\frac{\pi}{q}(x_{0}+\alpha s-2nd)\right)
    e^{2\pi ir\Psi(x_{0}+\alpha s,y_{0}+\beta t)}
    \dte t \right)\dte s,
\end{eqnarray*}
where $\Psi=\Psi_{n}^{a,b}$ as above. (If $(x_{0},y_{0})$ is a
non-degenerate critical point on the boundary of $\Sigma$, then
one or both or the integrals $\int_{-\delta}^{\delta}$ should
be replaced by $\int_{0}^{\delta}$ (or $\int_{-\delta}^{0}$)
and the contribution coming from that point in
\refthm{thm:asymptotics} should be multiplied by $\frac{1}{2}$
or $\frac{1}{4}$.) By a Taylor expansion we find that
$$
\Psi(x_{0} + \alpha s,y_{0} + \beta t) = \Psi(x_{0},y_{0})
+ \frac{1}{2}\la A\left(\begin{array}{c} s\\t\end{array}
\right),
\left(\begin{array}{c} s\\t\end{array}\right) \ra + h(s,t),
$$
where $A = \diag(\alpha,\beta) H \diag(\alpha,\beta)$ and
$$
\la \left(\begin{array}{c} x_{1}\\y_{1}\end{array}\right),
    \left(\begin{array}{c} x_{2}\\y_{2}\end{array}\right)\ra
= x_{1}x_{2} + y_{1}y_{2}
$$
for $(x_{1},y_{1}),(x_{2},y_{2}) \in \C^{2}$, and where
$h(s,t)$ is a remainder term being a sum of terms of the
form $c_{m,n}s^{m}t^{n}$, $m,n \in \{0,1,2,\ldots,...\}$,
$m+n \geq 3$, $c_{m,n} \in \C$. We search for $\alpha$ and
$\beta$ such that there exists a $\delta > 0$ satisfying
$$
\rea\left( 2\pi i (\Psi(x_{0} + \alpha s,y_{0} + \beta t)
             - \Psi(x_{0},y_{0}))\right) < 0
$$
for all
$(s,t) \in I_{\delta} \times I_{\delta} \sm \{ (0,0) \}$.
This amounts to finding $\alpha$ and $\beta$ such that
$$
\ima\left(\la A\left(\begin{array}{c} s\\t\end{array}\right),
\left(\begin{array}{c} s\\t\end{array}\right)\ra \right) > 0
$$
for all $(s,t) \in \R^{2} \sm \{ (0,0)\}$. Since
$$
\ima\left(\la A\left(\begin{array}{c} s\\t\end{array}\right),
    \left(\begin{array}{c} s\\t\end{array}\right)\ra \right)
= \la \ima(A) \left(\begin{array}{c} s\\t\end{array}\right),
              \left(\begin{array}{c} s\\t\end{array}\right) \ra
$$
this corresponds to finding $\alpha$ and $\beta$ such that
$\ima(A)$ is positive definite. It is essential for the method
that such $\alpha$ and $\beta$ exist. (If not, the analysis
becomes more difficult since higher order terms in the Taylor
expansion of $\Psi(x_0 +\alpha s,y_0 +\beta t)$ need to be
involved, and it is likely that the analysis can not be carried
through in this case.)
For that reason we make the following definition.

\begin{definition}\label{defi:posdefinite}
A critical point $(x,y)$ of the phase function $\Psi$ is called
positive definite if $(x,y)$ is non-degenerate and there exists
$(\alpha,\beta) \in S^1 \times S^1$ such that $\ima(A)$ is
positive definite, where
$$
A = \diag(\alpha,\beta) H \diag(\alpha,\beta),
$$
where $H$ is the Hessian of $\Psi$ in $(x,y)$.
\end{definition}

Assume that $(x_0,y_0)$ is positive definite and let
$(\alpha,\beta) \in S^1 \times S^1$ such that $\ima(A)$ is
positive definite. Then the main contribution to the integral
$K(x_{0},y_{0})$ in the limit of large $r$ is given by
$$
K_{\main}(x_{0},y_{0}) = \alpha\beta
 \sin\left(\frac{\pi}{q}(x_{0}-2nd)\right)
 e^{2\pi ir\Psi(x_{0},y_{0})} \int_{\R^{2}}
 e^{\pi ir\la A\left(\begin{array}{c} s\\t\end{array}\right),
 \left(\begin{array}{c} s\\t\end{array}\right)\ra}
 \dte s \dte t,
$$
and this integral can be evaluated using the results of
\cite[Sec.~3.4]{[12]}. In fact,
$$
\int_{\R^{2}} e^{\pi ir\la A
\left(\begin{array}{c} s\\t\end{array}\right),
\left(\begin{array}{c} s\\t\end{array}\right)\ra} \dte s \dte t
= \frac{1}{r} \left(\det(-iA)\right)^{-1/2},
$$
where $A$ is independent of $r$. Note here that the set $S$ of
complex symmetric $2 \times 2$--matrices $B$ with $\rea(B)$
positive definite is an open convex set in the $3$--dimensional
complex vector space of symmetric $2 \times 2$--matrices. It
follows that their is a unique analytic branch of
$B \mapsto \left( \det(B) \right)^{1/2}$ on $S$ such that
$\left( \det(B) \right)^{1/2} >0$ for $B$ real. We have used
that branch in the above result.
We note that $\det(-iA)=-(\alpha\beta)^2\det(H)$.
In conclusion we can state the
following result.

\begin{theorem}\label{thm:asymptotics}
Let $\mS_{k,l,a,b}^{\mu,\nu,n}$ be the set of critical points
of $\Phi_n^{a,b}$ in $\Sigma_{k,l,a,b}^{\mu,\nu,n}$. IF
\refconj{conj:mainconj} is true and if all the critical points
of $\Phi_n^{a,b}$ in $\Sigma_{k,l,a,b}^{\mu,\nu,n}$ are
positive definite, $k,l =0,1$, $(a,b) \in F_{k,l}$,
$\mu,\nu \in \{\pm 1 \}$, $n \in \Z/|q|\Z$, then the leading
order large $r$ asymptotics of the quantum invariant is given
by
\begin{eqnarray*}
\bar{\tau}_{r}(M_{p/q}) &\sim&
\frac{i\sign(q)}{4\sqrt{|q|}} e^{\frac{3\pi i}{4}\sign(pq)}
\qsum \sum_{(k,l)\in \{0,1\}^{2}} \sum_{(a,b)\in F_{k,l}}
\sum_{(\mu,\nu)\in \{\pm 1\}^{2}} \mu\nu \\
&& \times
\sum_{(x,y) \in \mS_{k,l,a,b}^{\mu,\nu,n}}
\left( -\det(H_{(x,y)}\Phi_{n}^{a,b})\right)^{-1/2} \\
&& \times \sin\left(\frac{\pi}{q}( x-2nd)\right)
         e^{2\pi ir\Phi_{n}^{a,b}(x,y)},
\end{eqnarray*}
where $H_{(x,y)}\Phi_{n}^{a,b}$ is the Hessian of
$\Phi_{n}^{a,b}$ in $(x,y)$.
\end{theorem}

According to the above theorem the growth rate of the quantum
invariants of $M_{p/q}$ is $r^{0}$ if we can prove that the
values of the phase functions $\Phi_n^{a,b}$ in the points
$S_{k,l,a,b}^{\mu,\nu,n}$ are real. This we do in Sec.~5.3,
see \refcor{cor:imaPsi}. This is in agreement with our computer
studies of the quantum invariants of $M_{p/q}$. To obtain
\refthm{thm:MT} from \refthm{thm:asymptotics}, hence a leading
asymptotics as predicted by the AEC, we have to prove i) that
the union of the sets of critical points
$\mS_{k,l,a,b}^{\mu,\nu,n}$ corresponds to the flat
(irreducible) $\SU(2)$--connections on $M_{p/q}$ and ii) that
the values of the relevant phase functions $\Phi_{n}^{a,b}$ in
these critical points are equal to the Chern--Simons invariants
of these flat connections. This we also do in Sec.~5.3. Note
that we have used the expression (\ref{eq:detHessian5}) for the
determinant of the Hessian $H_{(x,y)}\Phi_{n}^{a,b}$ in
\refthm{thm:MT}.

\begin{remark}\label{rem:Remark4}
The manifolds $M_{p/q}$, $|p/q| \in \{1,2,3\}$, are Seifert
manifolds and by \cite{[10]} (\ref{eq:asympqi}) holds for
these manifolds if we put $m_{\bar \rho}=4$ for all
$\bar \rho \in \mM_{p/q}'$.
Moreover, the results in
Appendix C show that the part of the leading order large $r$
asymptotics of $\tau_r(M_0)$ associated to the irreducible
$\SU(2)$--connections on $M_0$ is equal to the right-hand side
of (\ref{eq:asympqi}) if we put $m_{\bar \rho}=4$ for both
points $\bar \rho$ in $\mM_{0}'$. We actually expect that
$m_{\bar \rho}=4$ for all $\bar \rho \in \mM_{p/q}'$ for all
rational $p/q$. Recall here that when we calculated the leading
order large $r$ asymptotics of $J_{K}'(r)$ in Sec.~4.2 we had
to part the contour $C(\vep)$ into two parts $C_{\pm}(\vep)$
because of the tangent factor in the contour integral
(\ref{eq:contourfigure8}). The relevant stationary point was in
that case on the real axes, namely it was $x_0=5/6$, and to
obtain a contour passing through that point we had to deform
$C_{-}(\vep)$ to $[\vep,1-\vep]$ and $C_{+}(\vep)$ to
$[1-\vep,\vep]$. In that way the stationary point $x_0$
gave two equal contributions. (In fact we used Cauchy's
theorem so that we only had to work with either $C_{-}(\vep)$
or $C_{+}(\vep)$, but that of course also produced the factor
$2$.) In the case with the quantum invariants
$\bar \tau _r(M_{p/q})$ we expect a similar phenomenon with
the difference that each relevant stationary point contribute
with $4$ instead of $2$ equal contributions.
We expect that the relevant critical points $(x,y)$
belong to $[1/3,2/3] \times \{ y \in \C | \rea(y)=1/2\}$.
Because of the cotangent and tangent factors we had to
separate our double contour integral expression
for $\bar \tau _r(M_{p/q})$ in (\ref{eq:barinvariant}) into $4$
parts causing the sum
$\sum_{(\mu,\nu) \in \{\pm 1\}^2} \mu\nu$
in the asymptotic formula for $\bar{\tau}_r(M_{p/q})$. The
contour $C_{r}^1$ was thus separated into the part with
$\ima(x) \leq 0$ and the part with $\ima(x) \geq 0$. The first
of these parts can be deformed to $[0,1]$ while the second can
be deformed to $[1,0]$. The contour $C_{r}^2$ was separated
into the part with $\ima(y) \geq 0$ and the part with
$\ima(y) \leq 0$. Both
parts can be deformed to a contour containing the part of the
line $\rea(y)=1/2$ containing the $y$'s in the relevant
stationary points. In this way we obtain as claimed that each
relevant critical point contribute $4$ times. Like in the case
of the the knot invariant $J_{K}'(r)$ we find that the signs
$\mu$ and $\nu$ are cancelled by taking care of the
orientations of the different involved contours.
\end{remark}

Let us end this section by examine to what extend the critical
points $S_{k,l,a,b}^{\mu,\nu,n}$ are positive definite, see
\refdefi{defi:posdefinite}. Since $A$ is symmetric we have
$\ima(A)_{ij}=\ima(A_{ij})$.
Here
\begin{eqnarray*}
A_{11} &=& \alpha^{2}H_{11}, \\
A_{12} = A_{21} &=& \alpha\beta H_{12}, \\
A_{22} &=& \beta^{2}H_{22},
\end{eqnarray*}
where $H$ is the Hessian in (\ref{eq:Hessian}).
We have $z_0=e^{2\pi ix_0}$ and therefore
$H_{12}=H_{21}=2i\sin(2\pi x_0)$. Moreover,
$w_0 \in \{w_{\pm}(x_0)\}$ using notation from
(\ref{eq:wrealcircle2}) and $H_{11}=H_{22}-\frac{p}{2q}$, where
$$
H_{22}=\frac{1}{w_{\pm}(x_0)} - w_{\pm}(x_0) = \mp
2\sqrt{\cos^2(2\pi x_0) - \cos(2\pi x_0) - \frac{3}{4}}.
$$
If
$M=\left( \begin{array}{cc} a & c \\ c & b \end{array}\right)$
is a real symmetric matrix, then $M$ has real eigenvalues, and
$M$ is positive definite if and only if both of these
eigenvalues are positive, i.e.\ if and only if $a$ and $b$ are
both positive and $ab>c^{2}$. Now assume that $M=\ima(A)$.
Since $\rea(z_0) \in [-1,-1/2]$ we have that $H_{12}$ is zero
if and only if $z_0=-1$. In that case $H_{22} = \mp \sqrt{5}$
which is irrational. Therefore $H_{11} \neq 0$ and we can
always choose $\alpha$ and $\beta$ so $a>0$ and $b>0$ and hence
$ab>0=c^{2}$. Therefore assume in what follows that
$H_{12} \neq 0$. If $H_{11}$ and $H_{22}$ are both positive
(meaning that $w_0=w_{-}(x_0)$) we can let
$\alpha=\beta=e^{i\pi/4}$ giving $a=H_{11}$, $b=H_{22}$ and
$c=0$. If $H_{11}$ and $H_{22}$ are both negative (meaning that
$w_0=w_{+}(x_0)$) then we can let $\alpha=\beta=e^{-i\pi/4}$
giving $a=-H_{11}$, $b=-H_{22}$ and $c=0$.
If $\mu=\sign(H_{11}) \neq \sign(H_{22})$ we can let
$\alpha=e^{\mu i\pi/4}$ and $\beta=e^{-\mu i\pi/4}$ giving
$a=|H_{11}|$, $b=|H_{22}|$ and $c=D$, writing $H_{12}=iD$,
$D \in \R$. Thus $\det(M)=-\det(H)$, so if $\det(H)<0$ this
choice of $\alpha$ and $\beta$ works. The case
$\sign(H_{11}) \neq \sign(H_{22})$ and $\det(H)>0$ is more
problematic, and the above analysis is actually too simplified
as the examples $|p/q| \in \{4/3,10/3,16/3,22/3\}$ show.

One can in an even simpler way see that the above approach
is too simplified, namely by considering the cases where
$H_{11}=0$ or $H_{22}=0$. In that case $a=0$ or $b=0$ so
$\ima(A)$ is not positive definite. We have $H_{22}=0$ if and
only if $\cos(2\pi x_0)=-1/2$ and that case occurs if and only
if $p=6m+3$ for some integer $m$.

Since $\cos(2\pi x_0) \in [-1,-1/2]$ the entry $H_{11}$ can
not be zero if $|p/q|>\sqrt{20}$.
For $|p/q| < \sqrt{20}$, $H_{11}=0$ implies that
$\cos(2\pi x_0) = \frac{1}{2} - \sqrt{1+\frac{p^2}{16q^2}}$.
Finally note that if $H_{11}=H_{22}=0$ then
$p/q=0$, but in that case $H_{22} \neq 0$ for critical points
$(x_0,y_0) \in \mS$, see Appendix C.

The solution to these problems is to let the contour for the
``inner'' integral depend on the contour for the ``outer''
integral. We will not give more details here, but refer to
a later paper. Note that the expression for the matrix $A$ in
\refdefi{defi:posdefinite} in this more general approach is
different from the one stated in that definition.

\section{Classical Chern--Simons theory on $M_{p/q}$}

\noindent In this section we will describe the classical
theory, that is the classical Chern--Simons theory on the
manifolds $M_{p/q}$. The $\SU(2)$ Chern--Simons functional is
a map with values in $\R/\Z$ defined on the set $\mA$ of gauge
equivalence classes of connections in a principal $\SU(2)$
bundle on $M_{p/q}$ (all such bundles being trivializable).
Inside $\mA$ sits the moduli space $\mM_{p/q}$ of flat
$\SU(2)$--connections on $M_{p/q}$, that is the set of
classical solutions to the $\SU(2)$ Chern--Simons field theory.
Recall that
\[\mM_{p/q} = \Hom(\pi_{1}(M_{p/q}),\SU(2))/\SU(2) \]
from which it is clear that  $\mM_{p/q}$ is a compact space.
The Chern--Simons functional is constant on the connected
components of $\mM_{p/q}$, thus there are only finitely many
different values on flat connections. Let $\mM_{p/q}'$ be the
subset of $\mM_{p/q}$ consisting of nonabelian representations.
Recall that these representations correspond to the irreducible
connections, while the abelian representations correspond to
the reducible connections.

The main results in this section are \refthm{thm:comparison}
and \refthm{thm:main} which tie up the Chern--Simons theory
to the large $r$ asymptotics of $\bar{\tau}_{r}(M_{p/q})$ by
showing that a certain subset of the critical points of the
phase functions $\Phi_{n}^{a,b}$ corresponds to $\mM_{p/q}'$.
Under this correspondence, the Chern--Simons functional is
taken to the phase functions $\Phi_{n}^{a,b}$.

We begin by giving a description of $\mM_{p/q}$ following
Riley \cite{[28]}, \cite{[29]} and Kirk \& Klassen \cite{[18]}.

\subsection{The moduli space of flat $\SU(2)$--connections on
$M_{p/q}$}

\noindent In the following
$\pi=\pi_{1}\left( S^{3} \sm \nbd (K) \right)$ denotes the knot
group of the figure $8$ knot. We have a presentation
\begin{equation}\label{eq:pi}
\pi = \la \, x,y\, | wx=yw \,\ra,
\end{equation}
where $w=[x^{-1},y]$, and where $\mu=x$ and
$\lambda=yx^{-1}y^{-1}x^{2}y^{-1}x^{-1}y$ are the elements of
$\pi$ corresponding to the meridian and the preferred longitude
of $K$. The $\SL(2,\C)$ representation variety of $\pi$ was
analyzed by Riley \cite{[28]}, \cite{[29]} relevant to our
work. Consider a group $G$ given by a presentation
$$
G = \la \, x,y\, | wx=yw \,\ra,
$$
where
$w=x^{\vep_{1}}y^{\vep_{2}}x^{\vep_{3}}\cdots
y^{\vep_{\alpha-1}}$,
where $\alpha$ is odd and $\vep_{j} = \vep_{\alpha-j} = \pm 1$,
$j=1,2,\ldots,\alpha-1$. Such groups are denoted $2$--bridge
kmot groups by Riley since they generalize the $2$--bridge knot
groups. Following Riley we say that a representation
$\psi : G \to \SL(2,\C)$ is affine when the image of $\psi$
fixes exactly one point in $\C\Pro^{1}$ and not affine when
this image has no fixed points. We note that if $\psi$ is
nonabelian then $\psi$ is affine if and only if $\psi(x)$
and $\psi(y)$ have a common eigenvector, and $\psi$ is not
affine if these two matrices have no common eigenvector.
Let $H$ be some subgroup of $\SL(2,\C)$. Then we will say
that two representations $\psi_{1},\psi_{2} : G \to \SL(2,\C)$
are $H$--equivalent if they are conjugate to each other by a
matrix in $H$, i.e.\ if there exists a matrix $U \in H$ such
that $\psi_{2}(\gamma)=U\psi_{1}(\gamma)U^{-1}$ for all
$\gamma \in G$. In particular, we will say that $\psi_{1}$ and
$\psi_{2}$ are equivalent if they are  $\SL(2,\C)$--equivalent.
For $(t,u) \in \C^{*} \times \C$ we put
\begin{eqnarray*}
C_{0}(t) &=& \left( \begin{array}{cc} t & 1 \\
                0 & 1 \end{array}\right), \hspace{.2in}
D_{0}(t,u) = \left( \begin{array}{cc} t & 0 \\
                         -tu & 1 \end{array}\right), \\
C_{1}(t) &=& \left( \begin{array}{cc} t & 1 \\
                0 & t^{-1} \end{array}\right), \hspace{.2in}
D_{1}(t,u) = \left( \begin{array}{cc} t & 0 \\
                     -u & t^{-1} \end{array}\right), \\
C_{2}(t) &=& \left( \begin{array}{cc} t & t^{-1} \\
                0 & t^{-1} \end{array}\right), \hspace{.2in}
D_{2}(t,u) = \left( \begin{array}{cc} t & 0 \\
                     -tu & t^{-1} \end{array}\right).
\end{eqnarray*}
We note that $C_{\nu}(t)$ and $D_{\nu}(t,u)$ are elements of
$\SL(2,\C)$, $\nu=1,2$. If $s$ is a square root of $t$ and
$V(s)=\diag(s,s^{-1})$ then
\begin{equation}\label{eq:konjugering}
V(s)C_{2}(t)V(s)^{-1} = C_{1}(t), \hspace{.2in}
V(s)D_{2}(t,u)V(s)^{-1} = D_{1}(t,u),
\end{equation}
and
\begin{equation}\label{eq:C0C2}
C_{0}(t) = sC_{2}(s), \hspace{.2in} D_{0}(t,u) = sD_{2}(s,u).
\end{equation}
Let $W_{\nu}(t,u)$ denote the matrix obtained by replacing
$x$ and $y$ by respectively $C_{\nu}(t)$ and $D_{\nu}(t,u)$
in the expression for $w$, and let
$$
\phi(t,u) = W_{11} + (1-t)W_{12},
$$
where $W=W(t,u)=W_{0}(t,u)$. We let $\rho_{(t,u)}$ be the
assignment $x \mapsto C_{2}(t)$, $y \mapsto D_{2}(t,u)$. We
have the following $\SL(2,\C)$ version of
\cite[Theorem 1]{[28]}.

\begin{theorem}\label{thm:SL(2,C)}
Let $(s,u) \in \C^{*} \times \C$. None of the assignments
$\rho_{(s,u)}$ extend to an abelian
$\SL(2,\C)$--representation of $G$. The assignment
$\rho_{(s,u)}$ extends to a nonabelian representation
$\rho_{(s,u)} : G \to \SL(2,\C)$ if and only if
\begin{equation}\label{eq:phizero}
\phi(s^{2},u) = 0.
\end{equation}
Conversely, if $\psi : G \to \SL(2,\C)$ is a nonabelian
representation, then there exists a pair
$(s,u) \in \C^{*} \times \C$ satisfying
(\ref{eq:phizero}) such that $\psi$ and $\rho_{(s,u)}$
are equivalent. When $\psi$ is affine this pair is unique,
and when $\psi$ is not affine the pair $(s,u)$ can only be
replaced by $(s^{-1},u)$.
\end{theorem}

\begin{proof}
The theorem follows by results of \cite{[28]}, \cite{[29]}.
The assignment $\rho_{(s,u)}$ extends to a
$\SL(2,\C)$--representation of $G$ if and only if
$W_{2}(s,u)C_{2}(s,u) = D_{2}(s,u)W_{2}(s,u)$.
The matrices $C_{2}(s)$ and $D_{2}(s,u)$ commute if and only
if $s = \pm 1$ and $u=0$, and since $C_{2}(\pm 1)=\pm C_{2}(1)$
is different from $D_{2}(\pm 1,0)=\pm D_{2}(1,0)$ we have that
the assignment $\rho_{(\pm 1,0)}$ does not extend to a
$\SL(2,\C)$--representation of $G$. (We note that if $W=W(1,0)$
then $W_{21}=0$ since the matrices $C_{0}(1), D_{0}(1,0)=I$ and
their inverses are upper triangular. But we also have that
$\phi(1,0)=0$ would imply that $W_{11}=0$ contradicting the
fact that $W$ is invertible. Therefore $\phi(1,0) \neq 0$.)

Let $\sigma=\sum_{j=1}^{\alpha-1} \vep_{j}$. Then
$s^{\sigma}W_{2}(s,u)=W(s^{2},u)$. By (the proof of)
\cite[Theorem 1]{[28]} we have
$W(s^{2},u)C_{0}(s^{2}) = D_{0}(s^{2},u)W(s^{2},u)$
if and only if $\phi(s^{2},u)=0$, so by (\ref{eq:C0C2}) we
find that $\rho_{(s,u)}$ extends to a (necessarily nonabelian)
$\SL(2,\C)$--representation of $G$ if and only if
$\phi(s^{2},u)=0$.

If $\psi : G \to \SL(2,\C)$ is an arbitrary nonabelian
representation it follows by \cite[Lemma 7]{[29]} and
(\ref{eq:konjugering}) that there exists a pair
$(s,u) \in \C^{*} \times \C$ (necessarily satisfying
(\ref{eq:phizero})) such that $\psi$ and $\rho_{(s,u)}$
are equivalent. By the above any such pair is different from
$(\pm 1,0)$ and by \cite[Lemma 8]{[29]} and
(\ref{eq:konjugering}) we then get the final statement of
the theorem.
\end{proof}

If $(s,u) \in \C^{*} \times \C$ with $\phi(s^{2},u)=0$ then
$\rho_{(s,u)}$ is affine if and only if $C_{2}(s)$ and
$D_{2}(s,u)$ have a common eigenvector. But this happens
exactly when $u=0$ or $u=(s-s^{-1})^{2}$.

Let us now restrict to the case where $G$ is the figure $8$
knot group $\pi$. Then
\begin{equation}\label{eq:phifunction}
\phi(t,u)=u^{2}+\left(3-(t+t^{-1})\right)(u+1)
\end{equation}
so in particular $\phi(s^{2},0)=3-s^{2}-s^{-2}=0$ if and only
if $s^{4}-3s^{2}+1=0$ i.e.\ if and only if
$s=\mu_{1}\sqrt{(3+\mu_{2}\sqrt{5})/2}$ for some
$\mu_{1},\mu_{2} \in \{ \pm 1 \}$. If
$u=(s-s^{-1})^{2}=s^{2}+s^{-2}-2$ then
$\phi(s^{2},u) = u^{2} + (3-u-2)(u+1) = u^{2} + 1 - u^{2} = 1$,
so we conclude that $\rho_{(s,u)}$ is affine if and only if
$u=0$ and $s^{4}-3s^{2}+1=0$. Let
$$
\mN = \Hom(\pi,\SL(2,\C))/\SL(2,\C)
$$
be the space of conjugacy classes of
$\SL(2,\C)$--representations of $\pi$ and let $\mN_{\nab}$ be
the subset consisting of classes represented by nonabelian
$\SL(2,\C)$--representations. Moreover, let
\begin{equation}\label{eq:Ntilde}
\tilde{\mN} = \{\, (s,u) \in
        \C^{*} \times \C \, | \, \phi(s^{2},u)=0 \,\},
\end{equation}
and let $\Phi : \tilde{\mN} \to \mN$ be the map which maps
$(s,u)$ to the class represented by $\rho_{(s,u)}$. We have
shown

\begin{corollary}\label{cor:SL(2,C)}
The image of $\Phi$ is $\mN_{\nab}$. If we let
$$
\tilde{\mN}_{0} = \{ \,(\mu_{1}\sqrt{(3+\mu_{2}\sqrt{5})/2},0)
         \,|\,\mu_{1},\mu_{2} \in \{\pm 1\}\,\}
$$
then $\Phi|_{\tilde{\mN}_{0}} : \tilde{\mN}_{0} \to \mN$ is
injective and $\Phi^{-1}(\Phi(s,u))=\{(s,u),(s^{-1},u)\}$
for any $(s,u) \in \tilde{\mN} \sm \tilde{\mN}_{0}$.\HS
\end{corollary}

Recall that $M_{p/q}$ denotes the closed oriented $3$--manifold
obtained by surgery on $S^{3}$ along the figure $8$ knot with
rational surgery coefficient $p/q$. The representation
$\rho_{(s,u)}$, $(s,u) \in \tilde{\mN}$, extends to a
$\SL(2,\C)$--representation of $\pi_{1}(M_{p/q})$ if and only
if
$$
\rho_{(s,u)}(\mu)^{p}\rho_{(s,u)}(\lambda)^{q} = 1.
$$
To analyze this criterion it is an advantage to diagonalize
$\rho_{(s,u)}(\lambda)$. Assume in the following that
$(s,u) \in \tilde{\mN}$. By a rather long but completely
elementary and straightforward calculation we find that
$$
\rho_{(s,u)}(\lambda) = \left(
  \begin{array}{cc} \lambda_{11}(s,u) & \lambda_{12}(s,u) \\
            0 & \lambda_{11}(s^{-1},u) \end{array} \right),
$$
where
$$
\lambda_{11}(s,u) =
-1 + s^{-2} - 2s^{2} + s^{4} + u(s^{-2}-s^{2})
$$
and $\lambda_{12}(s,u) = 2u(1-u)$ if $s^{2}=1$ and
$$
\lambda_{12}(s,u) = \frac{\lambda_{11}(s,u)
              - \lambda_{11}(s^{-1},u)}{s^{2}-1}
$$
for $s^{2} \neq 1$. We note that
$\lambda_{11}(s^{-1},u) = \lambda_{11}(s,u)^{-1}$.

In general, if
$A = \left(\begin{array}{cc} \alpha & \beta \\
       0 & \alpha^{-1} \end{array}\right) \in \SL(2,\C)$,
then $A$ can be diagonalized if and only if $A$ is not
parabolic, i.e.\ if and only if $\tr(A) \neq \pm 2$ or
equivalently if and only if $\alpha \neq \pm 1$ (except, of
course, if $\beta=0$). If $\alpha \neq \pm 1$ then $(1,0)$ is
an eigenvector with eigenvalue $\alpha$ and
$(-\beta/(\alpha-\alpha^{-1}),1)$ is an eigenvector with
eigenvalue $\alpha^{-1}$. If $s^{2}=1$ then
$\lambda_{11}(s,u)=-1$ and $\lambda_{12}(s,u)=\pm i2\sqrt{3}$.
If $s^{2} \neq 1$ then $\lambda_{11}(s,u) = \pm 1$ if and only
if $\lambda_{12}(s,u)=0$. Therefore $\rho_{(s,u)}(\lambda)$ is
diagonalizable (or diagonal) if and only if $s^{2} \neq 1$. In
case $s^{2} \neq 1$ and $\lambda_{12}(s,u) \neq 0$ we have
$$
\frac{\lambda_{12}(s,u)}{\lambda_{11}(s,u)
         -\lambda_{11}(s,u)^{-1}}
= \frac{1}{s^{2}-1} = \frac{s^{-1}}{s-s^{-1}}.
$$
We conclude that if $s^{2} \neq 1$, then $\C^{2}$ has a basis
consisting of a set of common eigenvectors for the matrices
$\rho_{(s,u)}(\mu)$ and $\rho_{(s,u)}(\lambda)$, namely
$u_{1}=(1,0)$ and $u_{2}=(-1/(s^{2}-1),1)$. If we let
$\tilde{\rho}_{(s,u)} : \pi \to \SL(2,\C)$ be the
representation
$\tilde{\rho}_{(s,u)}(\gamma)=U^{-1}\rho_{(s,u)}(\gamma)U$,
where $U \in \SL(2,\C)$ with $j$th column $u_{j}$, we
therefore have
\begin{equation}\label{eq:rhoa1}
\tilde{\rho}_{(s,u)}(x) = \diag\left(s,s^{-1}\right), \qquad
\tilde{\rho}_{(s,u)}(\lambda) =
   \diag\left(\lambda_{11}(s,u),\lambda_{11}(s,u)^{-1}\right).
\end{equation}
In particular, $\rho_{(s,u)} : \pi \to \SL(2,\C)$,
$s^{2} \neq 1$, extends to a representation of
$\pi_{1}(M_{p/q})$ if and only if
\begin{equation}\label{eq:rhosurgery}
s^{-p} = \lambda_{11}(s,u)^{q}.
\end{equation}
Recall here that $\rho_{(s,u)}$ and $\rho_{(s^{-1},u)}$
are equivalent for $(s,u) \in \tilde{\mN} \sm \tilde{\mN}_{0}$,
cf.\ \refcor{cor:SL(2,C)}. But as noted above
$\lambda_{11}(s^{-1},u) = \lambda_{11}(s,u)^{-1}$ in accordance
with (\ref{eq:rhosurgery}). For $(s,u) \in \tilde{\mN}_{0}$
we have $\lambda_{11}(s,u)=1$ and $|s| \neq 1$, so
$\rho_{(s,u)}$ extends to a representation of
$\pi_{1}(M_{p/q})$ if and only if $p=0$.

A direct check shows that if $s^{2}=1$ then $\rho_{(s,u)}$
does not extend to a representation of $\pi_{1}(M_{p/q})$ for
any rational number $p/q$. In fact, if $s = \pm 1$, then
$$
\rho_{(s,u)}(\mu)^{p} = (\pm 1)^{|p|}
  \left(\begin{array}{cc} 1 & p \\ 0 & 1 \end{array}\right)
$$
and
$$
\rho_{(s,u)}(\lambda)^{q} =
       \left(-\left(\begin{array}{cc} 1 & \vep i2\sqrt{3} \\
            0 & 1 \end{array}\right)\right)^{q}
= (-1)^{|q|}\left(\begin{array}{cc} 1 & \vep i2q\sqrt{3} \\
          0 & 1 \end{array}\right)
$$
for a $\vep \in \{-1,1\}$. On the other hand, since
$\lambda_{11}(s,u) = -1$, we have that (\ref{eq:rhosurgery})
is satisfied if $s=1$ and $q$ is even or $s=-1$ and both $p$
and $q$ are odd. For the following we note that if $s^{2}=1$
then $u^{2} + u + 1 = \phi(1,u) = 0$ so $u$ is not real.

We are mostly interested in the $\SU(2)$--representations of
$\pi$. Let in the following $\mM_{\nab}$ be the set of
conjugacy classes of nonabelian $\SU(2)$--representations of
$\pi$.

\begin{proposition}\label{prop:SU(2)}
Let $(s,u) \in \tilde{\mN}$. The representation
$\rho_{(s,u)} : \pi \to \SL(2,\C)$ is $\SL(2,\C)$--equivalent
to a representation  $\pi \to \SU(2)$ if and only if $|s|=1$
and $u$ is real. If we write $s=e^{2\pi i\theta}$,
$\theta \in ]-1/2,1/2]$, then $u \in \{ u_{\pm} \}$, where
$u_{\pm} = u_{\pm}(\theta)$ are the two solutions to
$\phi(e^{4\pi i\theta},u) = 0$, i.e.\
$$
u_{\pm}(\theta) = \cos(4\pi\theta)-\frac{3}{2}
\pm \sqrt{\cos^{2}(4\pi\theta)-\cos(4\pi\theta)-\frac{3}{4}}.
$$
Since $u$ is real we have
$\theta \in [-1/3,-1/6] \cup [1/6,1/3]$.
The representation $\rho_{(e^{2\pi i \theta},u_{\vep})}$,
$\theta \in [-1/3,-1/6] \cup [1/6,1/3]$ and
$\vep \in \{\pm \}$, is $\SL(2,\C)$--equivalent to a
$\SU(2)$--representation $\bar{\rho}_{\theta,\vep}$ which
satisfies
$$
\bar{\rho}_{\theta,\vep}(\mu) = \diag\left(e^{2\pi i \theta},
         e^{-2\pi i\theta}\right),\hspace{.2in}
\bar{\rho}_{\theta,\vep}(\lambda) =
\diag(L_{\vep},L_{\vep}^{-1}),
$$
where $\mu$ and $\lambda$ are the elements of $\pi$
corresponding to the meridian and the preferred longitude of
$K$, and
$$
L_{\pm} = L_{\pm}(\theta) =
   \lambda_{11}\left(e^{2\pi i\theta},u_{\pm}\right)
= - 1 + e^{-4\pi i\theta} - 2e^{4\pi i\theta}
+ e^{8\pi i \theta}
+ u_{\pm}\left(e^{-4\pi i\theta} - e^{4\pi i\theta}\right).
$$
We note that $\bar{\rho}_{-\theta,\vep}$ and
$\bar{\rho}_{\theta,\vep}$ are $\SU(2)$--equivalent. In
particular, the space $\mM_{\nab}$ can be parametrized by the
two arcs
$(e^{2\pi i\theta},u_{+}(\theta))$, $\theta \in [1/6,1/3]$, and
$(e^{2\pi i\theta},u_{-}(\theta))$, $\theta \in [1/6,1/3]$.
These two arcs only coincide at the endpoints, so
topologically $\mM_{\nab}$ is a circle.
\end{proposition}

This proposition follows from \cite[Proposition 4]{[28]}, see
also \cite[Proposition 5.4]{[18]}. For a more geometric
argument determining the topological type of $\mM_{\nab}$,
see \cite{[19]}.

For $\theta \in [-1/3,-1/6] \cup [1/6,1/3]$ and
$\vep \in \{\pm\}$, the representation
$\bar{\rho}_{\theta,\vep}$ extends to a representation of
$\pi_{1}(M_{p/q})$ if and only if
$$
\bar{\rho}_{\theta,\vep}(\mu)^{p}
\bar{\rho}_{\theta,\vep}(\lambda)^{q} = 1,
$$
i.e.\ if and only if
\begin{equation}\label{eq:p/q}
e^{-2\pi ip\theta} = L_{\vep}(\theta)^{q}.
\end{equation}
From this (use e.g.\ (\ref{eq:Lreaima})) we see that

\begin{corollary}
Let $p/q \in \Q$ be arbitrary. The moduli space of irreducible
flat $\SU(2)$--connections on $M_{p/q}$ is a finite set.\HS
\end{corollary}

Let us end this section by finding the abelian
$\SU(2)$--representations of $\pi_{1}(M_{p/q})$ (up to
equivalence). Therefore, let $\theta \in ]-1/2,1/2]$
and let $\rho_{\theta}$ be the assignment
$$
\rho_{\theta}(\mu) =
\diag (e^{2\pi i \theta},e^{-2\pi i \theta}).
$$
By (\ref{eq:pi}) this assignment extends to an abelian
$\SU(2)$--representation of $\pi$ for any
$\theta \in ]-1/2,1/2]$ by letting
$\rho_{\theta}(y) = \rho_{\theta}(x)$. Moreover, any abelian
$\SU(2)$--representation of $\pi$ is $\SU(2)$--equivalent to
$\rho_{\theta}$ for some $\theta \in ]-1/2,1/2]$. For any
$\theta \in ]-1/2,1/2]$ we have $\rho_{\theta}(\lambda)=1$, and
$\rho_{\theta}$ extends to a representation of
$\pi_{1}(M_{p/q})$ if and only if $\rho_{\theta}(\mu)^{p} = 1$,
i.e.\ if and only if $p\theta \in \Z$. If
$A = \left(\begin{array}{cc} 0 & -1 \\
1 & 0 \end{array} \right)$
then $A\rho_{\theta}(\mu)A^{-1} = \rho_{-\theta}(\mu)$ so we
can assume that $\theta \in [0,1/2]$. Note, moreover, that two
matrices $\diag (e^{i \phi},e^{-i \phi})$ and
$\diag (e^{i \psi},e^{-i \psi})$, $\phi,\psi \in ]-\pi,\pi]$,
are conjugate in $\SU(2)$ if and only if $\phi=\psi$ or
$\phi =-\psi$. We conclude

\begin{proposition}\label{prop:SU(2)abelian}
For $p \neq 0$ the set of conjugacy classes of abelian
$\SU(2)$--represen-tations of $\pi_{1}(M_{p/q})$ is given by
$$
\left\{ \; [\rho_{j/|p|}] \; \left| \; j=0,1,\ldots,
    \left[\frac{|p|}{2}\right] \; \right. \right\},
$$
where $[|p|/2]$ is the integer part of $|p|/2$. For $p=0$
the set of conjugacy classes of abelian $\SU(2)$--representations
of $\pi_{1}(M_{p/q})$ is given by
$$
\left\{ \; [\rho_{\theta}] \; \left| \; \theta \in
           [0,\frac{1}{2}] \; \right. \right\},
$$
so topologically this set is a closed interval.\HS
\end{proposition}

\subsection{Chern--Simons invariants}

\noindent We begin by recalling formulas from \cite{[18]} for
the Chern--Simons invariants of the flat $\SU(2)$--connections
on $M_{p/q}$. The basic tool will be \refthm{thm:KK} below due
to P.\ A.\ Kirk and E.\ P.\ Klassen.

Let $M$ be a closed oriented $3$--manifold with a knot $K$ in
its interior and let $X$ be the complement of a tubular
neighborhood of $K$ in $M$. Moreover, let $\mu$ be a meridian
of $K$ and $\lambda$ a longitude, both in $\partial X$. Let $G$
be $\SU(2)$ or $\SL(2,\C)$. If $\rho : \pi_{1}(X) \to G$ is a
representation, then $\rho$ extends to a $G$--representation of
$\pi_{1}(M)$ if and only if $\rho(\mu)=1$. Now assume that
$\rho_{t} : \pi_{1}(X) \to G$ is a piecewise smooth path of
representations, $t \in I=[0,1]$. Choose a piecewise smooth
path $g : I \to G$ such that
\begin{eqnarray}\label{eq:pathkonjugering}
g_{t} \rho_{t}(\mu) g_{t}^{-1} &=&
\diag \left( e^{2\pi i \alpha(t)},e^{-2\pi i \alpha(t)}
\right), \nonumber \\*
g_{t} \rho_{t}(\lambda) g_{t}^{-1} &=&
\diag \left( e^{2\pi i \beta(t)},e^{-2\pi i \beta(t)} \right)
\end{eqnarray}
for some piecewise smooth curves $\alpha,\beta$. If $G=\SU(2)$
this is always possible by \cite[Lemma 3.1]{[18]}, and in that
case $\alpha$ and $\beta$ are real-valued. If $G=\SL(2,\C)$ the
above is possible if the path $\rho$ avoids the parabolic
representations (i.e.\ upper triangular with $1$s or $-1$s on
the diagonal), cf.\ \cite[Remark p.~354]{[18]}. (See the text
in connection to (\ref{eq:rhoa1}) for the case of the figure
$8$ knot.) In that case the curves $\alpha$ and $\beta$ are
complex-valued. We then have

\begin{theorem}[ {\cite[Theorem 4.2]{[18]}} ]\label{thm:KK}
Assume that $\rho_{0}(\mu)=\rho_{1}(\mu)=1$. Thinking of
$\rho_{0}$ and $\rho_{1}$ as flat $G$--connections on $M$, we
have
$$
\CS(\rho_{1})-\CS(\rho_{0})
= -2\int_{0}^{1} \beta(t)\alpha'(t) \dte t \pmod{\Z},
$$
where $\CS$ is the Chern--Simons functional associated to
$G$.\HS
\end{theorem}

We note that in case $G=\SL(2,\C)$ the Chern--Simons functional
takes values in $\C/\Z$.

Next consider a knot $K$ in $S^{3}$ and let $X$ be the knot
complement. Let $p/q$ be a rational number and let $N_{p/q}$ be
the closed oriented $3$--manifold obtained by $p/q$ surgery on
$S^{3}$ along $K$. Let $\mu$ and $\lambda$ be classes in
$\pi_{1}(\partial X)$ represented by respectively a meridian
and the preferred longitude of $K$. Choose integers
$c,d \in \Z$ such that $pd-qc=1$. Let $V$ be a tubular
neighborhood of $K$ considered as a subspace of $N_{p/q}$. We
note that $\mu'=p\mu+q\lambda$ and $\lambda'=c\mu+d\lambda$ are
represented by respectively a meridian of $V$ and a longitude
of $V$. Assume that $\rho_t : \pi_{1}(X) \to G$, $t \in I$, is
a piecewise smooth curve of representations from the trivial
representation to a representation, which extends to a
representation of $\pi_{1}(N_{p/q})$, i.e.\ $\rho_1(\mu')=1$.
Assume, moreover, that $\rho_t$ avoids the parabolic
representations in case $G=\SL(2,\C)$, and choose curves
$\alpha,\beta$ as in (\ref{eq:pathkonjugering}) with
$\alpha(0)=\beta(0)=0$. By \refthm{thm:KK} we have
\begin{eqnarray}\label{eq:CSsurgery}
\CS(\rho_1) &=& -2\int_{0}^{1}
(c\alpha(t)+d\beta(t))(p\alpha'(t)+q\beta'(t))\dte t \\*
&=& -2\int_{0}^{1} \beta(t)\alpha'(t)\dte t
    -cp\alpha^{2}(1) - dq\beta^{2}(1) - 2cq\alpha(1)\beta(1)
            \pmod{\Z}.\nonumber
\end{eqnarray}
(We have corrected a sign error in
\cite[Formula (*) p.~361]{[18]}.) We note that this expression
is independent of the choice of $c,d$. The condition
$\rho_1(\mu')=1$ is equivalent to
\begin{equation}\label{eq:p/q1}
p\alpha(1) + q\beta(1) \in \Z.
\end{equation}
Now let $K$ be the figure $8$ knot and let
$\bar{\rho}_{\theta,\vep}$ be a $\SU(2)$--representation of
$\pi_{1}(M_{p/q})$, i.e.\ (\ref{eq:p/q}) is satisfied.
Following \cite{[18]} we determine a formula for
$\CS(\bar{\rho}_{\theta,\vep})$. For later we will here pay
special attention to the branches of the logarithm. By
\refprop{prop:SU(2)} we have
\begin{eqnarray}\label{eq:Lreaima}
\rea\left(L_{-}(\theta)\right) &=&
\rea\left(L_{+}(\theta)\right)
= 2\cos^{2}(4\pi\theta) - \cos(4\pi\theta)-2, \nonumber \\*
\ima\left(L_{\pm}(\theta)\right) &=&
\mp 2\sin(4\pi\theta)\sqrt{\cos^{2}(4\pi\theta)
         - \cos(4\pi\theta) - \frac{3}{4}}
\end{eqnarray}
for all $\theta \in [-1/3,-1/6] \cup [1/6,1/3]$. From these
identities we see that
\begin{equation}\label{eq:Lspecialvalues}
L_{\pm}(1/3) = L_{\pm}(1/6) = -1, \hspace{.2in}
L_{\pm}(1/4) = 1,
\end{equation}
and that $\ima(L_{+})<0$ and $\ima(L_{-})>0$ on $]1/6,1/4[$
with the opposite signs on $]1/4,1/3[$. We conclude that
$L_{+}(\theta)$ and $L_{-}(\theta)$ run through $S^{1}$ both
beginning and ending in $-1$, $L_{+}$ in the anti-clockwise and
$L_{-}$ in the clockwise direction, as $\theta$ runs through
$[1/6,1/3]$. We can therefore use the principal logarithm
$\Log$ to define continuous curves
$\beta_{\pm} : [1/6,1/3] \to \R$ by
\begin{equation}\label{eq:beta}
\beta_{\pm}(\theta) = \frac{1}{2\pi i} \Log(L_{\pm}(\theta))
                      + f_{\pm}(\theta),
\end{equation}
where
\begin{equation}\label{eq:f+}
f_{+}(\theta) = \left\{
\begin{array}{cl} 0,& \theta = \frac{1}{6}, \\
1,& \theta \in ]\frac{1}{6},\frac{1}{3}],\end{array}\right.
\end{equation}
and
\begin{equation}\label{eq:f-}
f_{-}(\theta) = \left\{
\begin{array}{cl} 0,& \theta \in [\frac{1}{6},\frac{1}{3}[, \\
-1,& \theta = \frac{1}{3}.\end{array}\right.
\end{equation}
By (\ref{eq:p/q}) and the definition of $\beta_{\pm}$ the
representation $\bar \rho _{\theta,\vep}$ extends to a
representation of $\pi_1(M_{p/q})$ if and only if
\begin{equation}\label{eq:p/q'}
p\theta+q\beta_{\vep}(\theta) \in \Z
\end{equation}
which is nothing but (\ref{eq:p/q1}) (reparametrized).
We note that $\beta_{\pm}$ are smooth on
$]\frac{1}{6},\frac{1}{3}[$ but not in the end points $1/6$ and
$1/3$. The terms $f_{\pm}$ have been chosen so that
$\beta_{\pm}(1/6)= 1/2$. This is needed for the proof of

\begin{proposition}[ {\cite[p.~362]{[18]}} ]\label{prop:CS1}
Let $\theta \in [1/6,1/3]$ and let $\bar{\rho}_{\theta,\pm}$
be as in \refprop{prop:SU(2)}. If $\bar{\rho}_{\theta,\vep}$
extends to a representation of $\pi_{1}(M_{p/q})$ for a
$\vep \in \{\pm 1\}$ then
$$
\CS(\bar{\rho}_{\theta,\vep}) = -\frac{1}{6} - cp\theta^{2}
- dq\beta_{\vep}^{2}(\theta) - 2cq\theta\beta_{\vep}(\theta)
- 2\int_{1/6}^{\theta}\beta_{\vep}(t)\dte t \pmod{\Z},
$$
where $\beta_{\pm}$ are the curves defined by
(\ref{eq:beta}).\HS
\end{proposition}

Kirk and Klassen prove the above result by explicitly
constructing a piecewise smooth path
$\rho : [0,1] \to \Hom(\pi,\SL(2,\C))$ from the trivial
representation to
$\bar{\rho}_{\frac{1}{6},+} = \bar{\rho}_{\frac{1}{6},-}$ with
piecewise smooth curves $\alpha,\beta:[0,1] \to \C$ as in
(\ref{eq:pathkonjugering}) satisfying $\alpha(1)=\frac{1}{6}$
and $\beta(1)=\frac{1}{2}$. Moreover, they use that
$\int_{0}^{1} \beta(t)\alpha'(t) \dte t=\frac{1}{12}$ and the
fact that $y \mapsto \bar{\rho}_{y,\vep}$,
$[1/6,\theta] \to \Hom(\pi,\SU(2))$ is a path from
$\bar{\rho}_{\frac{1}{6},\vep}$ to $\bar{\rho}_{\theta,\vep}$
with associated functions $\alpha(y)=y$ and
$\beta(y)=\beta_{\vep}(y)$. By the choice of the functions
$f_{\pm}$, these $\alpha$-- and $\beta$--functions are
continuations of the ones used for the path $\rho$ from the
trivial representation to $\bar{\rho}_{\frac{1}{6},\pm}$.

Kirk and Klassen argue that
$\int_{0}^{1} \beta(t)\alpha'(t) \dte t = \frac{1}{12}$ using a
comparison between a computer calculation and the Chern--Simons
invariants of flat $\SU(2)$--connection on the Seifert manifold
$M_{-3}$. By following Kirk and Klassen's argument
\cite[pp.~361--362]{[18]} it is actually not hard to give an
explicit calculation of this integral. Let us give some
details. The path $\rho$ from the trivial connection to
$\bar{\rho}_{\frac{1}{6},\pm}$ consists of three parts, where
$\beta$ is identically zero along the first two parts. We can
therefore concentrate on the third part. Let
$a=(3+\sqrt{5})/2$, and let $\alpha : [0,1] \to \C$ be a
piecewise smooth curve from $\frac{1}{2\pi i} \Log(\sqrt{a})$
to $1/6$. Let $t(s)=e^{4\pi i \alpha(s)}$ and choose a
piecewise smooth solution $u(s)$ to $\phi(t(s),u(s))=0$, where
$\phi$ is given by (\ref{eq:phifunction}). Then
$s \mapsto \rho_{(e^{2\pi i\alpha(s)},u(s))} =: \eta_{s}$
is the third piece of our curve $\rho$ (reparametrized).
Assuming $\alpha(s) \notin \frac{1}{2}\Z$ (so as to avoid the
parabolic representations) we can diagonalize exactly as
demonstrated after \refcor{cor:SL(2,C)} and get
$$
\tilde{\eta}_{s}(\mu) = \diag(T,T^{-1}),\qquad
\tilde{\eta}_{s}(\lambda) =
\diag\left( \lambda_{11}(T,u),\lambda_{11}(T,u)^{-1}\right),
$$
where $T=T(s)=e^{2\pi i\alpha(s)}$, $u=u(s)$, $\lambda_{11}$
is as below \refcor{cor:SL(2,C)}, and where
$\tilde{\eta}_{s}(\gamma) = U^{-1}\eta_{s}(\gamma)U$ for
$\gamma \in \pi$, where
$U=U(s)=\left( \begin{array}{cc}
                1 & -1/(T^{2}-1) \\
                0 & 1
                \end{array}
\right)$.
This shows that $\alpha(s)$ indeed plays the role as the
$\alpha$--curve for our path $\eta_{s}$. The $\beta$--curve
should be a piecewise smooth curve $\beta(s)$ starting at zero
such that
\begin{equation}\label{eq:betany}
e^{2\pi i\beta(s)} = \lambda_{11}(e^{2\pi i\alpha(s)},u(s))
\end{equation}
for $s \in [0,1]$. We note that
$\lambda_{11}(e^{2\pi i\alpha(1)},u(1))=-1$ so we must have
$\beta(1) \in \frac{1}{2} + \Z$. Since $\bar{\rho}_{1/6,\pm}$
extends to a $\SU(2)$--representation of $\pi_{1}(M_{-3})$ we
get from (\ref{eq:CSsurgery}) that
$$
- 2\int_{0}^{1} \beta(s) \alpha'(s) \dte s - \frac{1}{12}
+ \frac{1}{3}\beta(1) \pmod{\Z}
$$
is a Chern--Simons value of a flat $\SU(2)$--connection on
$M_{-3}$, so $\int_{0}^{1} \beta(s) \alpha'(s) \dte s$ is real.

We now only have to choose a nice $\alpha$--curve. Let $\delta$
be a small positive parameter less than $|\alpha(0)| < 1/6$,
and let $\alpha = \alpha_{1} + \alpha_{2} + \alpha_{3}$, where
$\alpha_{1}$ is the line segment $[\alpha(0), -i\delta]$,
$\alpha_{2}$ is the part of the circle with centre zero and
radius $\delta$ running from $-i\delta$ to $\delta$, and
$\alpha_{3} = [\delta,1/6]$. Here the parameter $\delta$ is
introduced to avoid passing through zero (so as to avoid
parabolic representations). With this choice it is not hard to
show that there is a unique continuous solution
$(u(s),\beta(s))$ to $\phi(t(s),u(s))=0$ and (\ref{eq:betany})
with $\beta(0)=0$ and $\beta(1)=\frac{1}{2}$ and to show that
$\int_{0}^{1} \beta(s) \alpha'(s)\dte s = \frac{1}{12}$ for
this solution. (Use that the integral is real and independent
of $\delta$ and calculate its $\delta \ria 0_{+}$ limit).

\refprop{prop:CS1} gives a formula for the Chern--Simons
invariants of the irreducible flat $\SU(2)$--connections on
$M_{p/q}$. Let us also determine the Chern--Simons invariants
of the reducible flat $\SU(2)$--connections on $M_{p/q}$.
Therefore let $\rho_{j/|p|}$ be as in
\refprop{prop:SU(2)abelian}, $p \neq 0$, and let
$\alpha(t)=jt/|p|$ and $\beta(t)=0$, $t \in [0,1]$, and get by
(\ref{eq:CSsurgery}) that
\begin{equation}\label{eq:CSabelianpnot0}
\CS(\rho_{j/|p|}) = -cpj^{2}/p^{2} = -cj^{2}/p \pmod{\Z},
\end{equation}
where $c$ is the inverse of $-q \bmod{p}$.

In case $p=0$ the moduli space of flat reducible
$\SU(2)$--connections on $M_{0}$ can be identified with a
closed interval, cf.\ \refprop{prop:SU(2)abelian}. Since the
Chern--Simons functional is constant on each of the connected
components of the moduli space of flat connections, we conclude
that the Chern--Simons invariant of any of the reducible
$\SU(2)$--connections is equal to the Chern--Simons invariant
of the trivial connection, i.e.\ it is equal to zero. This of
course also follows from Kirk and Klassen's result. In fact,
if $\rho_{\theta}$ denotes the abelian representation from
\refprop{prop:SU(2)abelian}, then we can put
$\alpha(t)=\theta t$ and $\beta(t) =0$, $t \in [0,1]$, and get
by (\ref{eq:CSsurgery}) that
\begin{equation}\label{eq:CSabelianp0}
\CS(\rho_{\theta}) = 0 \pmod{\Z}.
\end{equation}

\newpage

\subsection{A comparison between Chern--Simons invariants and
critical values of the phase functions $\Phi_{n}^{a,b}$}

\noindent The purpose of this section is to combine the
Chern--Simons theory described in the previous two sections
with the asymptotic analysis in Sec.~4.3. First we will
describe a correspondence between the critical points of the
phase functions $\Phi_{n}^{a,b}$ in
(\ref{eq:Phasefunctiontypical}) and the nonabelian
$\SL(2,\C)$--representations of $\pi_{1}(M_{p/q})$. Thereafter
we will show that the set of Chern--Simons invariants of flat
irreducible $\SU(2)$--connections on $M_{p/q}$ is a certain
subset of the critical values of the functions
$\Phi_{n}^{a,b}$. Like in Sec.~4.3 we will most of the time
work with the shifted phase functions $\Psi_{n}^{a,b}$ in
(\ref{eq:Phasefunctiontypical1}). We refer to
\refrem{rem:PhiPsi} below for a comparison between the phase
functions $\Phi_{n}^{a,b}$ and $\Psi_{n}^{a,b}$. Recall the set
$\tilde{\mN}$ given by (\ref{eq:Ntilde}).

\begin{theorem}\label{thm:comparison}
The map $(x,y) \mapsto \rho_{(e^{\pi ix},e^{2\pi iy} -1)}$
gives a surjection $\varphi$ from the set of critical points
$(x,y)$ of the functions $\Psi^{a,b}_{n}$, $a,b,n \in \Z$, with
$x \notin \Z$ onto the set of representations
$\rho_{(s,u)} : \pi \to \SL(2,\C)$, $(s,u) \in \tilde{\mN}$,
which extend to $\SL(2,\C)$--representations of
$\pi_{1}(M_{p/q})$.

If $(x,y)$ is a critical point, let us say of $\Psi_{n}^{a,b}$,
then $(x+2k,y+l)$ is a critical point of
$\Psi_{n+pk}^{a+l,b+2k}$ for any $k,l \in \Z$, so
$\varphi^{-1}(\varphi(x,y))=\{\;(x+2k,y+l)\;|\;k,l \in \Z\;\}$
if $x \notin \Z$.
\end{theorem}

\begin{proof}
The last statement follows immediately from
(\ref{eq:critical}), so let us concentrate on the first part.
Let $\psi : \C^{2} \to \C$ be given by
$$
\psi(z,w) = (1-zw)(w-z) - zw,
$$
so $\psi(v^{2},w)=0$ if and only if $(v,w)$ is a solution to
(\ref{eq:stationary2}). Then
$$
\psi(z,u+1) = -z\phi(z,u)
$$
for $(z,u) \in \C^{*} \times \C$, where $\phi$ is the function
(\ref{eq:phifunction}). It follows that $(v,w)$ is a solution
to (\ref{eq:stationary2}) if and only if either $(v,w)=(0,0)$
or $(v,w-1) \in \tilde{\mN}$.

Next assume that $(v,u) \in \tilde{\mN}$ and put
$(z,w)=(v^{2},u+1)$. Since $\psi(z,w)=0$ and $z \neq 0$ we have
that $w \neq 0$ and $w-z \neq 0$. Therefore $1-zw = wz/(w-z)$
and $(v,w)$ is a solution to (\ref{eq:stationary1}) if and only
if
\begin{eqnarray*}
v^{-p} &=& \left(\frac{(w-z)^{2}}{wz}\right)^{q}
    = \left( wz^{-1} - 2 + zw^{-1} \right)^{q} \\*
&=& \left( \lambda_{11}(v,w-1) - 1 - z^{2}
               + zw^{-1} + wz + z \right)^{q}.
\end{eqnarray*}
But $\psi(z,w)=0$ and $w \neq 0$ implies that
$-1-z^{2}+zw^{-1}+wz+z=0$. Thus $(v,w)$ is a solution to
(\ref{eq:stationary1}) if and only if $(s,u)=(v,u)$ is a
solution to (\ref{eq:rhosurgery}).

The above shows together with the discussion around
(\ref{eq:rhosurgery}) that there is a one to one correspondence
between common nonzero solutions $(v,w)$ to
(\ref{eq:stationary1}), (\ref{eq:stationary2}) with
$v^{2} \neq 1$ and representations
$\rho_{(v,w-1)} : \pi \to \SL(2,\C)$,
$(v,w-1) \in \tilde{\mN}$, which extend to
$\SL(2,\C)$--representations of $\pi_{1}(M_{p/q})$. By the
remarks following (\ref{eq:critical}) this proves the theorem.
\end{proof}

Recall the set $\mS$ in (\ref{eq:SUcriticalpoint}). Moreover,
recall that if $(x,y)$ is a critical point of $\Psi_n^{a,b}$,
$a,b,n \in \Z$, with $(x,y) \in \mS$ then $x \notin \Z$ is
automatically satisfied. By \refprop{prop:SU(2)} and the above
theorem we therefore have

\begin{corollary}
The surjection $\varphi$ in \refthm{thm:comparison} restricts
to a surjection from the set of critical points $(x,y)$ of
$\Psi_n^{a,b}$, $a,b,n \in \Z$, with $(x,y) \in \mS$ onto the
set of representations
$\rho_{(s,u)} : \pi_1(M_{p/q}) \to \SL(2,\C)$,
$(s,u) \in \tilde{\mN}$, which are equivalent to
$\SU(2)$--representations of $\pi_{1}(M_{p/q})$.
\end{corollary}

\begin{theorem}\label{thm:main}
Let $(x,y) \in \C^{2}$ such that
$\rho_{(e^{\pi ix},e^{2\pi iy} -1)}$
is equivalent to a nonabelian $\SU(2)$--representation of
$\pi_{1}(M_{p/q})$ and choose in accordance with
\refthm{thm:comparison} integers $a,b,n$ such that $(x,y)$ is a
critical point of $\Psi_{n}^{a,b}$. If $a',b',n'$ is another
such set of integers, then $b'=b$ and $a'+n'/q=a+n/q$. In fact
we have
\begin{eqnarray*}
a + \frac{n}{q} &=& y + \frac{p}{2q}x
+ \frac{i}{2\pi} \left(\Log\left(1-e^{2\pi i(x+y)}\right)
- \Log\left(1-e^{2\pi i(x-y)}\right)\right), \\*
b &=& x
+ \frac{i}{2\pi} \left(\Log\left(1-e^{2\pi i(x+y)}\right)
+ \Log\left(1-e^{2\pi i(x-y)}\right)\right).
\end{eqnarray*}
Moreover, for any such set of integers $a,b,n$ we have
$$
\CS(\bar{\rho}_{\theta,\vep}) = \Psi_{n}^{a,b}(x,y) \pmod{\Z},
$$
where $\theta \in [-1/3,-1/6] \cup [1/6,1/3]$ and
$\vep \in \{ \pm \}$ with $e^{2\pi i\theta}=e^{\pi ix}$ and
$1+u_{\vep}(\theta) = e^{2\pi iy}$, and
$\bar{\rho}_{\theta,\vep}$ is a $\SU(2)$--representation of
$\pi_{1}(M_{p/q})$ equivalent to
$\rho_{(e^{\pi ix},e^{2\pi iy} -1)}$ and given by
\refprop{prop:SU(2)}.
\end{theorem}

\begin{remark}\label{rem:PhiPsi}
i) In general we have $\Phi_{n}^{a,b} = \Psi_{n}^{a+b,a-b}$,
so there is a ono to one correspondence between the phase
functions $\Phi_{n}^{a,b}$ and the phase functions
$\Psi_{n}^{a',b'}$ with $a'+b'$ even.
If $(x,y)$ is a critical point
of $\Psi_{n}^{a',b'}$ with $a'+b'$ odd, then $(x,y)$ is also a
critical point of $\Psi_{n-q}^{a'+1,b'}$ by \refthm{thm:main}.
In that way we see that any critical point of one of the
functions $\Psi_{n}^{a',b'}$ is also a critical point of one
of the functions $\Phi_{n}^{a,b}$.

\noindent ii) Assume that $|p/q| > \sqrt{20}$ and assume that
$(x_j,y_j)$, $j=1,2$, are degenerate critical points of phase
functions $\Phi_n^{a,b}$ belonging to $\mS$. Let
$v_j=e^{\pi ix_j}$ and $w_j=e^{2\pi iy_j}$. By
\refprop{prop:Nondeg} we have
$\cos(2\pi x_1) = \cos(2\pi x_2)$ and $w_j=w_{-\vep}(x_j)$,
where $\vep=\sign(p/q)$, so $w_1=w_2$. Moreover,
$x_1, x_2 \in [-2/3,-1/3] \cup [1/3,2/3] + 2\Z$. We are
interested in
comparing the equivalence classes of the representations
$\rho_j=\rho_{(v_j,w_j - 1)}$. For that purpose we can without
loss of generality assume that $x_1,x_2 \in [1/3,2/3]$. Then
$x_2=x_1$ or $x_2=1-x_1$. In the second case $v_2 = - v_1$.
Note here that $-v_1 = v_1^{-1}$ if and only if
$v_1 \in \{ \pm i \}$. But in that case $\cos(2\pi x_1)=-1$ and
by (\ref{eq:detHessian4}) we get that
$p/q \in \{ \pm 2/\sqrt{5} \}$ which is impossible. Thus
$\rho_1$ and $\rho_2$ are not equivalent. If $p$ is odd we
must have $x_2=x_1$ since only one of the points
$(v_1,w_1)$, $(-v_1,w_1)$ can be a solution to
(\ref{eq:stationary1}). If $p$ is even both $(v_1,w_1)$ and
$(-v_1,w_1)$ are solutions to both of the equations
(\ref{eq:stationary1}) and (\ref{eq:stationary2}). Thus
\refprop{prop:Nondeg} shows that the set of degenerate critical
points in $\mS$ of the phase functions $\Phi_n^{a,b}$,
$a,b,n \in \Z$, is either empty or corresponds under the
surjection in \refthm{thm:comparison} to one point in
$\mM_{p/q}$ in case $p$ is odd and to two points in the moduli
space in case $p$ is even.
\end{remark}

The phase functions $\Psi_{n}^{0,0}$, $\Psi_{n}^{1,-1}$,
$\Psi_{n}^{1,1}$ and $\Psi_{n}^{2,0}$ are the ones entering
\refconj{conj:mainconj}. By \refthm{thm:main} the parameter
$b$ is fixed by a critical point of $\Psi_{n}^{a,b}$, while
$a$ and $n$ can be varied as long as we keep fixed $a+n/q$.
Because of this (see also \reflem{lem:Phivaluesym} below) we
only need to consider $\Psi_{n}^{0,0}$, $\Psi_{n}^{0,-1}$, and
$\Psi_{n}^{0,1}$. By \refthm{thm:main} (and also the proof of
this theorem and \refprop{prop:SU(2)}) we get

\begin{corollary}
Let $(x,y) \in \C^{2}$ such that
$\rho_{(e^{\pi ix},e^{2\pi iy} -1)}$
is equivalent to a nonabelian $\SU(2)$--representation of
$\pi_{1}(M_{p/q})$. Let $(x',y') \in \C^{2}$ such that
$(e^{\pi ix'},e^{2\pi iy'})=(e^{\pi ix},e^{2\pi iy})$ and
$\rea(x') \in ]-1,1]$ and $\rea(y') \in ]-1/2,1/2]$. Moreover,
let
\begin{eqnarray*}
n &=& q\left(y'+\frac{p}{2q}x'
+ \frac{i}{2\pi} \left(\Log\left(1-e^{2\pi i(x+y)}\right)
- \Log\left(1-e^{2\pi i(x-y)}\right)\right)\right), \\
b &=& x'
+ \frac{i}{2\pi} \left(\Log\left(1-e^{2\pi i(x+y)}\right)
+ \Log\left(1-e^{2\pi i(x-y)}\right)\right).
\end{eqnarray*}
Then $\theta:= x'/2 \in [-1/3,-1/6] \cup [1/6,1/3]$ and
$$
b = \left\{ \begin{array}{rl} -1,& \theta \in [-1/3,-1/4] \\
0, & \theta \in ]-1/4,-1/6] \cup [1/6,1/4] \\
1, & \theta \in ]1/4,1/3].
\end{array}\right.
$$
Moreover, $n \in \Z$ and $(x',y')$ is a critical point of
$\Psi_{n}^{0,b}$ and
$$
\CS(\bar{\rho}_{\theta,\vep}) =
\Psi_{n}^{0,b}(x',y') \pmod{\Z},
$$
where $\vep \in \{ \pm \}$ with
$1+u_{\vep}(\theta) = e^{2\pi iy}$, and
$\bar{\rho}_{\theta,\vep}$ is a $\SU(2)$--representation
of $\pi_{1}(M_{p/q})$ equivalent to
$\rho_{(e^{\pi ix},e^{2\pi iy} -1)}$ and given by
\refprop{prop:SU(2)}.
\end{corollary}

We note that the first part of \refthm{thm:main} is an
immediate consequence of (\ref{eq:critical}). To prove the
second part we start by observing that if $a,a',n,n',b$ are
integers such that $a'+n'/q=a+n/q$ then
$\Psi_{n}^{a,b}(x,y)-\Psi_{n'}^{a',b}(x,y)$ is an integer
independent of $(x,y) \in \C^{2}$. The remaining part of the
proof will consists of a series of lemmas. We start by

\begin{lemma}\label{lem:Phivaluesym}
Let $a,b,n \in \Z$ and assume that $(x,y)$ is a critical point
of $\Psi_{n}^{a,b}$. Let $l,k \in \Z$ and put $a'=a+l$,
$b'=b+2k$, and $n'=n+pk$. Then $(x+2k,y+l)$ is a critical point
of $\Psi_{n'}^{a',b'}$ and
$$
\Psi_{n'}^{a',b'} (x+2k,y+l) - \Psi_{n}^{a,b}(x,y) \equiv
0 \pmod{\Z}.
$$
\end{lemma}

\begin{proof}
That $(x+2k,y+l)$ is a critical point of $\Psi_{n'}^{a',b'}$ is
an immediate consequence of (\ref{eq:critical}) and was already
observed in \refthm{thm:comparison}. Put $x'=x+2k$ and
$y'=y+l$. Since $e^{2\pi i x'}=e^{2\pi ix}$ and
$e^{2\pi iy'}=e^{2\pi iy}$ we get
\begin{eqnarray*}
\Psi_{n'}^{a',b'} (x',y') - \Psi_{n}^{a,b}(x,y) &=&
a'x' + b'y' -x'y' + \frac{n'}{q}x' -\frac{p}{4q}{x'}^{2}
     - \frac{d}{q}{n'}^{2} \\
&& - ax - by + xy - \frac{n}{q}x + \frac{p}{4q}x^{2}
   + \frac{d}{q}n^{2} \\
&=& 2(a+l)k  + 2\frac{n}{q}(1-pd)k + \frac{p}{q} (1-pd) k^{2}.
\end{eqnarray*}
Here $pd-cq =1$ for an integer $c$ so
$$
\Psi_{n'}^{a',b'} (x',y') - \Psi_{n}^{a,b}(x,y)
= 2(a+l-nc)k - pck^{2}.
$$
\end{proof}

By \refprop{prop:SU(2)} and the above lemma we are thus left
with the following to prove: Let
$\theta \in [-1/3,-1/6] \cup [1/6,1/3]$ and assume that
$\bar{\rho}_{\theta,\vep}$ extends to a representation of
$\pi_{1}(M_{p/q})$ for $\vep \in \{\pm\}$. Let $b=b(\theta)$
and $n=n(\theta)$ be the unique integers such that
$(2\theta,y)$ is a stationary point of $\Psi_{n}^{0,b}$, where
$y=\frac{1}{2\pi i} \Log(1+u_{\vep}(\theta))$ (so
$\rea(y) = \frac{1}{2}$) (see around (\ref{eq:critical2}) and
(\ref{eq:critical1})). Then
\begin{equation}\label{eq:mainstatement}
\CS(\bar{\rho}_{\theta,\rho}) = \Psi_{n}^{0,b}(2\theta,y) \pmod{\Z}.
\end{equation}
In the course of the proof of this identity we will prove that
\begin{equation}\label{eq:CS2}
\CS(\bar{\rho}_{\theta,\vep}) = -\frac{1}{6}
- \frac{p}{q}\theta^{2} + \frac{m}{q}2\theta - \frac{d}{q}m^{2}
- 2\int_{1/6}^{\theta}\beta_{\vep}(t) \dte t \pmod{\Z},
\end{equation}
for $\theta \in [1/6,1/3]$, where
$m= p\theta + q\beta_{\vep}(\theta) \in \Z$ by (\ref{eq:p/q'}).
Unfortunately the proof of (\ref{eq:mainstatement}) is rather
technical, but we have tried to emphasize the main ideas of the
proof here and defer many technicalities to the appendices C
and D.

The Chern--Simons invariants of $\SU(2)$--connections are real,
so we begin by investigating to what extent a general phase
function $\Psi_{n}^{a,b}$, $a,b,n \in \Z$, is real in its
critical points. Assume that $(x,y)$ is such a critical point
and write $z=e^{2\pi i x}$ and $w=e^{2\pi i y}$ as usual. From
(\ref{eq:critical}) we find that
$$
a - \rea(y) = \frac{p}{2q}\rea(x) - \frac{n}{q}
- \frac{1}{2\pi} \ima \left(\Log(1-zw)-\Log(1-zw^{-1})\right),
$$
and
$$
b - \rea(x) =
- \frac{1}{2\pi} \ima\left(\Log(1-zw)+\Log(1-zw^{-1})\right).
$$
Moreover, we have $\ima(x) = -\frac{1}{2\pi}\Log|z|$ and
$\ima(y)=-\frac{1}{2\pi}\Log|w|$. By
(\ref{eq:Phasefunctiontypical1}) we therefore get
\begin{eqnarray*}
\ima\left(\Psi_{n}^{a,b}(x,y)\right) &=& (a-\rea(y))\ima(x)
+ (b-\rea(x))\ima(y) - \frac{p}{2q}\ima(x)\rea(x) \\
&& + \frac{n}{q}\ima(x) + \frac{1}{4\pi^{2}}
\ima\left( \Li_{2}(zw)-\Li_{2}(zw^{-1})\right).
\end{eqnarray*}
This together with the above expressions for $a-\rea(y)$ and
$b-\rea(x)$ leads to the formula
\begin{eqnarray*}
\ima\left(\Psi_{n}^{a,b}(x,y)\right) &=& \frac{1}{4\pi^{2}}
\left( \ima(\Log(1-zw)-\Log(1-zw^{-1}))\Log|z| \right. \\
&& + \ima(\Log(1-zw) + \Log(1-zw^{-1}))\Log|w| \\
&& + \left.
\ima\left( \Li_{2}(zw)-\Li_{2}(zw^{-1})\right) \right).
\end{eqnarray*}
By introducing the Bloch--Wigner dilogarithm function
$$
D(z) = \ima\left( \Li_{2}(z) \right) + \Arg(1-z) \Log|z|
$$
we obtain

\begin{lemma}\label{lem:imaPsi}
Let $a,b,n \in \Z$ and let $(x,y)$ be a critical point of
$\Psi_{n}^{a,b}$. Then
$$
\ima\left(\Psi_{n}^{a,b}(x,y)\right)
= \frac{1}{4\pi^{2}} \left( D\left(e^{2\pi i(x+y)}\right)
  - D\left(e^{2\pi i(x-y)}\right)\right).
$$\HS
\end{lemma}

We note that $D$ is analytic on $\C \sm \{ 0,1\}$ and
continuous on $\C$. Moreover, $D$ satisfies the identities
\begin{equation}\label{eq:WienerBlochidentities}
D(z) + D(\bar{z}) = 0, \hspace{.2in} D(z) + D(1/z) = 0
\end{equation}
for all $z \in \C \sm \{0\}$. From this we have

\begin{corollary}\label{cor:imaPsi}
Let $a,b,n \in \Z$ and let $(x,y)$ be a critical point of
$\Psi_{n}^{a,b}$ with $e^{2\pi ix} \in S^{1}$ and
$e^{2\pi i y } \in \R$. Then
$$
\ima\left(\Psi_{n}^{a,b}(x,y)\right) = 0.
$$\HS
\end{corollary}

\begin{remark}
Let $(v,w)$ be a non-zero solution to (\ref{eq:stationary1})
and (\ref{eq:stationary2}). Then $(\bar{v},\bar{w})$ is also
a solution to these two equations as already observed. Let
$z=v^{2}$ and write $(z,w)=(e^{2\pi ix},e^{2\pi iy})$. Then
$(\bar{z},\bar{w})=(e^{-2\pi i\bar{x}},e^{-2\pi i\bar{y}})$.
If $(x,y)$ is a critical point of $\Psi_{n}^{a,b}$, then
$(-\bar{x},-\bar{y})$ is a critical point of
$\Psi_{-n}^{-a,-b}$. By \reflem{lem:imaPsi} and
(\ref{eq:WienerBlochidentities}) we have
$$
\ima\left(\Psi_{n}^{a,b}(x,y)\right) =
    -\ima\left(\Psi_{-n}^{-a,-b}(-\bar{x},-\bar{y})\right).
$$
Hence, if this value is different from zero, then either
$\exp(2\pi ir\Psi_{-n}^{-a,-b}(-\bar{x},-\bar{y}))$ or
$\exp(2\pi ir\Psi_{n}^{a,b}(x,y))$ grows exponentially. We
note that by \refconj{conj:mainconj} and
\refcor{cor:imaPsi} stationary points leading to such
exponential growth do not contribute to the large $r$
asymptotics of $\bar{\tau}_{r}(M_{p/q})$, see also
\refthm{thm:asymptotics}.
\end{remark}

We now embark upon proving (\ref{eq:mainstatement}). We start
by reducing to the case $\theta \in [1/6,1/3]$. By
\refprop{prop:SU(2)} we know that the representations
$\bar{\rho}_{\theta,\vep}$ and $\bar{\rho}_{-\theta,\vep}$ are
$\SU(2)$--equivalent, so in particular they have the same
Chern--Simons invariant.

\begin{lemma}\label{lem:symmetry}
Let the situation be as in (\ref{eq:mainstatement}). Then
$$
\Psi_{n(\theta)}^{0,b(\theta)}(2\theta,y)
- \Psi_{n(-\theta)}^{0,b(-\theta)}(-2\theta,y) \equiv
0 \pmod{\Z}.
$$
\end{lemma}

To prove this and (\ref{eq:mainstatement}) we need the
technical \reflem{lem:Qfunctions} below keeping track of
branches of the logarithm in certain expressions. Let
$I=[-1/3,-1/6] \cup [1/6,1/3]$ and put
\begin{eqnarray}\label{eq:Q}
Q_{1}^{\pm}(\theta) &=& 1
- \left( 1+u_{\pm}(\theta) \right)^{-1}
e^{4\pi i\theta}, \nonumber \\
Q_{2}^{\pm}(\theta) &=& 1
- \left(1+u_{\pm}(\theta)\right) e^{4\pi i\theta}
\end{eqnarray}
for $\theta \in I$, where $u_{\pm}$ are defined in
\refprop{prop:SU(2)}. We also put
$Q_{3}^{\pm}(\theta)=1+u_{\pm}(\theta)$.

\begin{lemma}\label{lem:Qfunctions}
We have
\begin{equation}\label{eq:Qlog1}
\Log\left(Q_{1}^{\pm}(\theta)\right)
+ \Log\left(Q_{2}^{\pm}(\theta)\right)
= \Log\left(Q_{1}^{\pm}(\theta)Q_{2}^{\pm}(\theta)\right)
\end{equation}
and
\begin{equation}\label{eq:Qidentity1}
Q_{1}^{\pm}(\theta)Q_{2}^{\pm}(\theta) = e^{4\pi i\theta}
\end{equation}
for all $\theta \in I$. Moreover
\begin{equation}\label{eq:Qlog2}
\Log\left(Q_{1}^{\pm}(\theta)\right)
+ \Log\left(Q_{3}^{\pm}(\theta)\right)
- \Log\left(Q_{2}^{\pm}(\theta)\right)
= \Log\left(\frac{Q_{1}^{\pm}(\theta)Q_{3}^{\pm}(\theta)}
                 {Q_{2}^{\pm}(\theta)}\right)
+ e_{\pm}(\theta)2\pi i
\end{equation}
for all $\theta \in I$, where
$$
e_{+}(\theta) = \left\{
\begin{array}{cc} 1, & \theta \in ]1/6,1/4],\\
0, & \theta \in ]1/4,1/3[,\end{array}\right.
$$
and
$e_{-}(\theta) = 1 - e_{+}(\theta)$ for $\theta \in ]1/6,1/3[$,
$e_{\pm}(-1/4) = e_{\pm}(1/4)$,
$e_{\pm}(\theta) = 1 - e_{\pm}(-\theta)$ for
$\theta \in I \sm \{ \pm 1/6,\pm 1/4,\pm 1/3\}$, and
$e_{\pm}(\theta)=0$ for $\theta \in \{ \pm 1/6,\pm 1/3\}$.
Finally
\begin{equation}\label{eq:Qidentity2}
\frac{Q_{1}^{\pm}(\theta)Q_{3}^{\pm}(\theta)}
     {Q_{2}^{\pm}(\theta)}
= L_{\pm}(\theta)
\end{equation}
for all $\theta \in I$, where $L_{\pm}$ are given by
\refprop{prop:SU(2)}.
\end{lemma}

The proof of this lemma is given in Appendix D.

\begin{proofa}{\it of\,} \reflem{lem:symmetry}\,\,
By (\ref{eq:critical2}) and (\ref{eq:critical1}) we have
$$
b(\theta) = 2\theta - \frac{1}{2\pi i}
\left( \Log\left(Q_{1}^{\vep}(\theta)\right)
+ \Log\left(Q_{2}^{\vep}(\theta)\right)\right)
$$
and
$$
n(\theta) = p\theta + qy + q\frac{1}{2\pi i}
\left( \Log\left(Q_{1}^{\vep}(\theta)\right)
- \Log\left(Q_{2}^{\vep}(\theta)\right)\right).
$$
By \reflem{lem:Qfunctions}
\begin{equation}\label{eq:b}
b(\theta) = 2\theta
- \frac{1}{2\pi i} \Log\left( e^{4\pi i\theta} \right).
\end{equation}
Taking the real part of the expression for $n(\theta)$ we get
\begin{equation}\label{eq:n1}
n(\theta) = p\theta + \frac{1}{2}q
+ \frac{q}{2\pi} \left( \Arg\left(Q_{1}^{\vep}(\theta)\right)
- \Arg\left(Q_{2}^{\vep}(\theta)\right)\right).
\end{equation}
By (\ref{eq:b}) we get that
$$
b(-\theta) = - b(\theta), \hspace{.2in}
\theta \in [-1/3,-1/6] \cup [1/6,1/3] \sm \{ \pm 1/4\},
$$
and
$$
b(-1/4)=-1, \hspace{.2in} b(1/4)=0.
$$
By (\ref{eq:n1}), (\ref{eq:Qsym1}), and (\ref{eq:Qsym2}) we get
$$
n(-\theta) = q - n(\theta), \hspace{.2in}
\theta \in [-1/3,-1/6] \cup [1/6,1/3]\sm \{ \pm 1/4\},
$$
and
$$
n(-1/4) = n(1/4) - \frac{p}{2}.
$$
Let $b=b(\theta)$ and $n=n(\theta)$. Since
$\Psi_{n(\theta)}^{0,b(\theta)}(2\theta,y)$ is real by
\refcor{cor:imaPsi} and since
$\rea\left( \Li_{2}(\bar{z})\right) =
\rea\left( \Li_{2}(z)\right)$
for $z \in \C \sm ]1,\infty[$ we find that
\begin{eqnarray*}
\Psi_{n(\theta)}^{0,b(\theta)}(2\theta,y)
- \Psi_{n(-\theta)}^{0,b(-\theta)}(-2\theta,y)
&=& (b-b(-\theta))\rea(y) -4\theta\rea(y) \\
&& + \frac{(n+n(-\theta)}{q}2\theta
   - \frac{d(n^{2}-n(-\theta)^{2})}{q}.
\end{eqnarray*}
(If $\theta \in \{\pm 1/4\}$, use that
$\Li_{2}(e^{2\pi i(2\theta +y)})
- \Li_{2}(e^{2\pi i(2\theta-y)})$
is the same for $\theta =1/4$ and $\theta=-1/4$.) Assume first
that $\theta \neq \pm 1/4$. Then, since $\rea(y)=1/2$, we get
$$
\Psi_{n(\theta)}^{0,b(\theta)}(2\theta,y)
- \Psi_{n(-\theta)}^{0,b(-\theta)}(-2\theta,y)
= b + dq -2dn \in \Z.
$$
Next assume that $\theta = 1/4$. Then one finds
$$
\Psi_{n(\theta)}^{0,b(\theta)}(2\theta,y)
- \Psi_{n(-\theta)}^{0,b(-\theta)}(-2\theta,y)
= \frac{n}{q}(1-pd) - \frac{p}{4q}(1 - pd).
$$
But $1-pd=-qc$ for an integer $c$ so
$$
\Psi_{n(\theta)}^{0,b(\theta)}(2\theta,y)
- \Psi_{n(-\theta)}^{0,b(-\theta)}(-2\theta,y)
= \left(\frac{p}{4} -n \right)c.
$$
By (\ref{eq:p/q'}) $p/4 + q\beta_{\vep}(1/4) \in \Z$. But
$L_{\pm}(1/4)=1$ so $\beta_{\vep}(1/4) = f_{\vep}(1/4) \in \Z$,
so $p$ is divisible by $4$.
\end{proofa}

By using that
$y=\frac{1}{2\pi i} \Log\left(Q_{3}^{\vep}(\theta)\right)$
together with \reflem{lem:Qfunctions} we obtain the alternative
formula
$$
n(\theta) = p\theta + qe_{\vep}(\theta)
+ \frac{q}{2\pi i} \Log\left(L_{\vep}(\theta)\right).
$$
This and (\ref{eq:beta}) immediately leads to
\begin{equation}\label{eq:n3}
n(\theta) = p\theta + q\beta_{\vep}(\theta)
+ q(e_{\vep}(\theta)-f_{\vep}(\theta))
\end{equation}
for $\theta \in [1/6,1/3]$.
In the proof of (\ref{eq:mainstatement}) we need certain
symmetries for the functions $L_{\pm}$.
First we note that
\begin{equation}\label{eq:Lsym2}
L_{\pm}\left(\frac{1}{2} -\theta\right) = L_{\mp}(\theta)
\end{equation}
for $\theta \in [1/6,1/3]$ by (\ref{eq:Lreaima}).
Second, by the next lemma and \refprop{prop:SU(2)}, we have
\begin{equation}\label{eq:Lsym1}
L_{\mp}(\theta) = L_{\pm}(\theta)^{-1}
\end{equation}
for $\theta \in [-1/3,-1/6] \cup [1/6,1/3]$. (In particular,
$\bar{\rho}_{\theta,\pm}$ extends to a representation of
$\pi_{1}(M_{p/q})$ if and only if $\bar{\rho}_{\theta,\mp}$
extends to a representation of $\pi_{1}(M_{-p/q})$.)

\begin{lemma}\label{lem:lambdasymmetry}
Let $s \in \C^{*}$ and let $u_{\pm}$ be the two solutions to
$\phi(s^{2},u)=0$. Then
$$
(1+u_{+})(1+u_{-}) = 1,
$$
and
$$
\lambda_{11}(s,u_{+})\lambda_{11}(s,u_{-}) = 1.
$$
\end{lemma}

\begin{proof}
We have
$$
1 + u_{\pm} = \frac{1}{2}(s^{2}+s^{-2}-1)
\pm \frac{1}{2}\sqrt{s^{4}+s^{-4}-2(s^{2}+s^{-2})-1}
$$
so
$$
(1 + u_{+})(1 + u_{-}) = \frac{1}{4}(s^{2} + s^{-2} -1)^{2}
- \frac{1}{4}(s^{4} + s^{-4} - 2(s^{2} + s^{-2}) - 1) = 1.
$$
Moreover,
\begin{eqnarray*}
\lambda_{11}(s,u_{\pm}) &=& \frac{1}{2}(s^{4} + s^{-4})
- \frac{1}{2}(s^{2} + s^{-2}) - 1 \\*
&& \pm \frac{1}{2}(s^{-2} - s^{2})\sqrt{s^{4} + s^{-4}
- 2(s^{2} + s^{-2}) -1}
\end{eqnarray*}
so
\begin{eqnarray*}
\lambda_{11}(s,u_{+})\lambda_{11}(s,u_{-}) &=&
\left(\frac{1}{2}(s^{4} + s^{-4})
- \frac{1}{2}(s^{2} + s^{-2}) - 1 \right)^{2} \\*
&& - \frac{1}{4}(s^{2} - s^{-2})^{2}
\left(s^{4} + s^{-4} - 2(s^{2} + s^{-2}) - 1 \right).
\end{eqnarray*}
By simple reductions one gets the result.
\end{proof}

We are now ready to finalize the proof of
(\ref{eq:mainstatement}) and thereby the proof of
\refthm{thm:main}.

\begin{proofa}{\it of\,} (\ref{eq:mainstatement})\,\,
By \reflem{lem:symmetry} we can assume that
$\theta \in [1/6,1/3]$. Write $\theta_{0}$ for $\theta$ in the
following. Let us first observe that formula (\ref{eq:CS2}) is
an immediate consequence of \refprop{prop:CS1}. In fact, by
letting $c,d$ be integers as in \refprop{prop:CS1} we get
\begin{eqnarray*}
&& -2cq\theta_{0}\beta_{\vep}(\theta_{0}) = 2cp\theta_{0}^{2}
- 2cm\theta_{0}, \\
&& -dq\beta_{\vep}^{2}(\theta_{0}) = -\frac{d}{q}\left(m^{2}
- 2pm\theta_{0}+p^{2}\theta_{0}^{2}\right),
\end{eqnarray*}
and therefore
$$
- cp\theta_{0}^{2} - dq\beta_{\vep}^{2}(\theta_{0})
- 2cq\theta_{0}\beta_{\vep}(\theta_{0})
= \left(c - \frac{dp}{q}\right)p\theta_{0}^{2}
- 2\left(c - \frac{dp}{q}\right)m\theta_{0} - \frac{d}{q}m^{2}.
$$
But $c-dp/q=-1/q$, hence (\ref{eq:CS2}) follows. On the other
hand we have
$$
\Psi_{n}^{0,b}(2\theta_{0},y) = (b-2\theta_{0})y
+ \frac{n}{q}2\theta_{0} -\frac{p}{q}\theta_{0}^{2}
- \frac{d}{q} n^{2} + \frac{1}{4\pi^{2}}
\left( \Li_{2}(z_{0}w_{\vep})
- \Li_{2}(z_{0}w_{\vep}^{-1}) \right),
$$
where $z_{0}=e^{4\pi i\theta_{0}}$ and
$w_{\vep}=1+u_{\vep}(\theta_{0})$. By (\ref{eq:b}) and
(\ref{eq:n3}) we have $b-2\theta_{0}=-\frac{1}{2\pi i}\Log(z_{0})$
and $n=m+q\left(e_{\vep}(\theta_{0})-f_{\vep}(\theta_{0})\right)$
so
\begin{eqnarray*}
\Psi_{n}^{0,b}(2\theta_{0},y) &=& -\frac{p}{q}\theta_{0}^{2}
+ \frac{m}{q}2\theta_{0} - \frac{d}{q} m^{2}
+ 2\left(e_{\vep}(\theta_{0})
- f_{\vep}(\theta_{0})\right)\theta_{0}
- l_{\vep}(\theta_{0}) \\*
&& + \frac{1}{4\pi^{2}}
\left(\Log(z_{0})\Log(w_{\vep})+ \Li_{2}(z_{0}w_{\vep})
- \Li_{2}(z_{0}w_{\vep}^{-1}) \right),
\end{eqnarray*}
where
$l_{\vep}(\theta_{0}) = 2dm\left(e_{\vep}(\theta_{0})
- f_{\vep}(\theta_{0})\right) + dq\left(e_{\vep}(\theta_{0})
- f_{\vep}(\theta_{0})\right)^{2} \in \Z$.
We therefore get that
\begin{eqnarray*}
\Psi_{n}^{0,b}(2\theta_{0},y)
- \CS(\bar{\rho}_{\theta_{0},\vep}) &=& \frac{1}{6}
+ 2\left(e_{\vep}(\theta_{0})
- f_{\vep}(\theta_{0})\right)\theta_{0} \\*
&& + \frac{1}{4\pi^{2}} \left(\Log(z_{0})\Log(w_{\vep})
+ \Li_{2}(z_{0}w_{\vep})
- \Li_{2}(z_{0}w_{\vep}^{-1}) \right) \\*
&& + 2\int_{1/6}^{\theta_{0}}\beta_{\vep}(t) \dte t \pmod{\Z}.
\end{eqnarray*}
For $\theta_{0} \in [1/6,1/4]$ we note that
$e_{\vep}(\theta_{0}) = f_{\vep}(\theta_{0})$. We will consider
the special cases $\theta_{0} \in \{1/6,1/4,1/3\}$ first and
then handle the other cases afterwards.\newline
{\bf The cases $\theta_{0} \in \{1/6,1/3\}$.} In these
cases we have $w_{\vep}=-1$ so
\begin{eqnarray*}
&& \Psi_{n}^{0,b}(2\theta_{0},y)
- \CS(\bar{\rho}_{\theta_{0},\vep}) \\*
&& \hspace{.3in} = \frac{1}{6}+\frac{i}{4\pi} \Log(z_{0})
+ 2\left(e_{\vep}(\theta_{0})
- f_{\vep}(\theta_{0})\right)\theta_{0}
+ 2\int_{1/6}^{\theta_{0}}\beta_{\vep}(t) \dte t \pmod{\Z}.
\end{eqnarray*}
If $\theta_{0}=1/6$ we immmediately get that this is zero. If
$\theta_{0}=1/3$ we get
\begin{eqnarray*}
\Psi_{n}^{0,b}(2\theta_{0},y)
- \CS(\bar{\rho}_{\theta_{0},\vep})
&=& \frac{1}{6}+\frac{i}{4\pi} (4\pi i/3-2\pi i)
+ \frac{2}{3}\left(e_{\vep}(\theta_{0})
- f_{\vep}(\theta_{0})\right) \\*
&& + 2f_{\vep}(1/4)(1/3-1/6) \\*
&& + \frac{1}{\pi i}
\int_{1/6}^{1/3}\Log(L_{\vep}(t)) \dte t \pmod{\Z}.
\end{eqnarray*}
But if $\theta_0 \in [1/4,1/3]$ then
$$
\int_{1/4}^{\theta_{0}} \Log(L_{\pm}(\theta)\dte \theta
= -\int_{1/4}^{1/2-\theta_{0}}\Log(L_{\mp}(t))\dte t
= - \int_{1/2-\theta_{0}}^{1/4}\Log(L_{\pm}(t))\dte t,
$$
by (\ref{eq:Lsym2}) and (\ref{eq:Lsym1}) so
\begin{equation}\label{eq:Lintegral}
\int_{1/2-\theta_{0}}^{\theta_{0}}
\Log(L_{\pm}(\theta)\dte \theta = 0.
\end{equation}
In particular,
$$
\Psi_{n}^{0,b}(2\theta_{0},y)
- \CS(\bar{\rho}_{\theta_{0},\vep})
= \frac{1}{3} - \frac{2}{3}f_{\vep}(1/3)
+ \frac{1}{3}f_{\vep}(1/4) \pmod{\Z},
$$
where we also use that $e_{\pm}(1/3)=0$ by
\reflem{lem:Qfunctions}. By (\ref{eq:f+}) and (\ref{eq:f-})
this is zero.\newline
{\bf The case $\theta_{0} = 1/4$.} In this case we have
$z_{0}=-1$ so
\begin{eqnarray*}
\Psi_{n}^{0,b}(2\theta_{0},y)
- \CS(\bar{\rho}_{\theta_{0},\vep})
&=& \frac{1}{6}
+ \frac{1}{4\pi^{2}} \left(i\pi\Log(w_{\vep})
+ \Li_{2}(|w_{\vep}|) - \Li_{2}(|w_{\vep}|^{-1}) \right) \\*
&& + 2\int_{1/6}^{1/4}\beta_{\vep}(t) \dte t \pmod{\Z}.
\end{eqnarray*}
Since $(v_{0},w_{\vep})=(e^{2\pi i\theta_{0}},w_{\vep})$ is a
solution to (\ref{eq:stationary2}) we have
$w_{\vep} \in \{ (-3-\sqrt{5})/2,(-3+\sqrt{5})/2\}$, and by
(\ref{eq:1+u}) we then conclude that
$w_{\vep} = \frac{-3-\sqrt{5}}{2}$ if $\vep=-$ and
$w_{\vep} = \frac{-3+\sqrt{5}}{2}$ if $\vep=+$. By
(\ref{eq:dilog1/4}) we then get
\begin{eqnarray*}
\Psi_{n}^{0,b}(2\theta_{0},y)
- \CS(\bar{\rho}_{\theta_{0},\vep})
&=& \frac{1}{6} + \frac{1}{4\pi^{2}}\left(-\pi^{2}
- \vep \frac{\pi^{2}}{5}\right)
+ 2\int_{1/6}^{1/4}\beta_{\vep}(t) \dte t \\*
&=& \frac{1}{6}
+ \frac{1}{4}\left(-1 - \vep\frac{1}{5}\right)
+ \frac{1}{6}f_{\vep}(1/4) \\*
&& + \frac{1}{\pi} \int_{1/6}^{1/4}\Arg(L_{\vep}(t)) \dte t
\pmod{\Z},
\end{eqnarray*}
and this is zero by (\ref{eq:f+}), (\ref{eq:f-}), and
(\ref{eq:Specialintegral1/4}).\newline
{\bf The case $\theta \in ]1/6,1/3[ \sm \{1/4\}$.}
We have
\begin{eqnarray*}
\Psi_{n}^{0,b}(2\theta_{0},y)
- \CS(\bar{\rho}_{\theta_{0},\vep})
&=& \frac{1}{6}+\frac{1}{4\pi^{2}} \Log(z_{0})\Log(w_{\vep})
+ 2\left(e_{\vep}(\theta_{0})
- f_{\vep}(\theta_{0})\right)\theta_{0} \\*
&& + 2\int_{1/6}^{\theta_{0}}\beta_{\vep}(t) \dte t
+ R(2\theta_{0},y) \pmod{\Z},
\end{eqnarray*}
where
$$
R(x,y) = \frac{1}{4\pi^{2}}\left( \Li_{2}(e^{2\pi i(x+y)})
- \Li_{2}(e^{2\pi i(x-y)}) \right).
$$
Let us write $u$ for $u_{\vep}$ in the following. By definition
of the dilogarithm we have
$$
4\pi^{2}R(2\theta_{0},y) =
\int_{0}^{(1+u)^{-1}e^{4\pi i \theta_{0}}}
\frac{\Log(1-t)}{t} \dte t
- \int_{0}^{(1+u)e^{4\pi i \theta_{0}}}
\frac{\Log(1-t)}{t} \dte t.
$$
By (\ref{eq:1+u}) we have
$1+u_{-}(\theta) \leq -1 \leq 1+u_{+}(\theta) <0$ for all
$\theta \in [1/6,1/3]$ and $1+u_{\pm}(\theta)=-1$ if and only
if $\theta \in \{1/6,1/3\}$. Let $\theta_{1}=1/6$ if
$\theta_{0} \in [1/6,1/4[$ and $\theta_{1}=1/3$ if
$\theta_{0} \in ]1/4,1/3[$. We note that
$z \mapsto \Log(1-z)/z$ is analytic on $\C \sm [1,\infty[$ so
by Cauchy's theorem
\begin{eqnarray*}
4\pi^{2}R(2\theta_{0},y) &=& \int_{0}^{-e^{4\pi i \theta_{1}}}
\frac{\Log(1-t)}{t} \dte t +
\int_{-e^{4\pi i \theta_{1}}}^{(1+u)^{-1}e^{4\pi i \theta_{0}}}
\frac{\Log(1-t)}{t} \dte t \\*
&& - \int_{0}^{-e^{4\pi i \theta_{1}}}
\frac{\Log(1-t)}{t} \dte t
- \int_{-e^{4\pi i \theta_{1}}}^{(1+u)e^{4\pi i \theta_{0}}}
\frac{\Log(1-t)}{t} \dte t \\*
&=&
\int_{-e^{4\pi i \theta_{1}}}^{(1+u)^{-1}e^{4\pi i \theta_{0}}}
\frac{\Log(1-t)}{t} \dte t
- \int_{-e^{4\pi i \theta_{1}}}^{(1+u)e^{4\pi i \theta_{0}}}
\frac{\Log(1-t)}{t} \dte t .
\end{eqnarray*}
The curves
$\gamma_{\pm}(\theta)=(1+u(\theta))^{\pm 1} e^{4\pi i \theta}$
are smooth on $]1/6,1/3[$ so
\begin{eqnarray*}
4\pi^{2}R(2\theta_{0},y) &=&
\lieta \left( \int_{\theta_{1} + \mu\eta}^{\theta_{0}}
\Log(1-\gamma_{-}(\theta))\frac{\gamma_{-}'(\theta)}
{\gamma_{-}(\theta)} \dte \theta \right.\\*
&& \left. \hspace{.5in}
- \int_{\theta_{1} + \mu\eta}^{\theta_{0}}
\Log(1-\gamma_{+}(\theta))\frac{\gamma_{+}'(\theta)}
{\gamma_{+}(\theta)} \dte \theta \right),
\end{eqnarray*}
where $\mu =1$ if $\theta_{1}=1/6$ and $\mu=-1$ if
$\theta_{1}=1/3$. (The parameter $\eta$ is necessary because
$u$ is not differentiable in $1/6$ and $1/3$.) It follows that
$$
4\pi^{2}R(2\theta_{0},y) = \lieta \left(
4\pi i R_{1}(\theta_{0},\eta) - R_{2}(\theta_{0},\eta)\right),
$$
where
\begin{eqnarray*}
R_{1}(\theta_{0},\eta) &=&
\int_{\theta_{1} + \mu\eta}^{\theta_{0}}
\Log\left(Q_{1}^{\vep}(\theta)\right)
- \Log\left(Q_{2}^{\vep}(\theta)\right) \dte \theta, \\
R_{2}(\theta_{0},\eta) &=&
\int_{\theta_{1} + \mu\eta}^{\theta_{0}}
\left\{ \Log\left(Q_{1}^{\vep}(\theta)\right)
+ \Log\left(Q_{2}^{\vep}(\theta)\right) \right\}
\frac{u'(\theta)}{1+u(\theta)} \dte \theta,
\end{eqnarray*}
where the functions $Q_{i}^{\pm}$ are defined above
\reflem{lem:Qfunctions}. By \reflem{lem:Qfunctions}
\begin{eqnarray*}
R_{2}(\theta_{0},\eta) &=&
\int_{\theta_{1} + \mu\eta}^{\theta_{0}}
\Log\left(e^{4\pi i\theta}\right)
\frac{u'}{1+u} \dte \theta \\*
&=& 4\pi i \int_{\theta_{1}+\mu\eta}^{\theta_{0}} \theta
\frac{u'(\theta)}{1+u(\theta)} \dte \theta
- 2\pi ib(\theta_{0}) \int_{\theta_{1} + \mu\eta}^{\theta_{0}}
\frac{u'(\theta)}{1+u(\theta)} \dte \theta.
\end{eqnarray*}
Since
$u'(\theta)/(1+u(\theta))
= \frac{\dte}{\dte \theta} \log(1+u(\theta))$
for any branch $\log$ of the logarithm defined on an open
section of $\C^{*}$ containing $]-\infty,0[$ we have
$$
\int_{\theta_{1}+\mu\eta}^{\theta_{0}} \theta
\frac{u'(\theta)}{1+u(\theta)} \dte \theta =
\left[ \theta\Log(1+u(\theta))
\right]_{\theta_{1}+\mu\eta}^{\theta_{0}}
- \int_{\theta_{1}+\mu\eta}^{\theta_{0}}
\Log(1+u(\theta)) \dte \theta,
$$
and
$$
\int_{\theta_{1}+\mu\eta}^{\theta_{0}}
\frac{u'(\theta)}{1+u(\theta)} \dte \theta
= \left[ \Log(1+u(\theta))
\right]_{\theta_{1}+\mu\eta}^{\theta_{0}}.
$$
We therefore get
\begin{eqnarray*}
R(2\theta_{0},y) &=& \frac{1}{4\pi^{2}}\lieta \left(
    4\pi i R_{1}(\theta_{0},\eta) - R_{2}(\theta_{0},\eta)
\right) \\*
&=& \frac{b(\theta_{0})}{2} - \theta_{1}
-\frac{1}{4\pi^{2}} \Log(z_{0})\Log(w_{\vep}) \\*
&& + \frac{i}{\pi}\lieta \left( R_{1}(\theta_{0},\eta)
+ \int_{\theta_{1} + \mu\eta}^{\theta_{0}}
\Log\left(Q_{3}^{\vep}(\theta)\right)\dte \theta \right),
\end{eqnarray*}
where we use that
$1 + u_{\vep}(\theta) = Q_{3}^{\vep}(\theta) = w_{\vep}$ and
(\ref{eq:b}). By \reflem{lem:Qfunctions} we get
$$
R_{1}(\theta_{0},\eta) +
\int_{\theta_{1} + \mu\eta}^{\theta_{0}}
\Log\left(Q_{3}^{\vep}(\theta)\right)\dte \theta
=  \int_{\theta_{1} + \mu\eta}^{\theta_{0}}
\Log\left(L_{\vep}(\theta)\right)
+ 2\pi i e_{\vep}(\theta)\dte \theta,
$$
so
$$
\lieta\left(R_{1}(\theta_{0},\eta) +
\int_{\theta_{1} + \mu\eta}^{\theta_{0}}
\Log\left(Q_{3}^{\vep}(\theta)\right)\dte \theta\right)
= e_{\vep}(\theta_{0}) 2\pi i(\theta_{0} - \theta_{1})
+ \int_{\theta_{1}}^{\theta_{0}}
\Log(L_{\vep}(\theta))\dte \theta.
$$
But then
\begin{eqnarray*}
\Psi_{n}^{0,b}(2\theta_{0},y)
- \CS(\bar{\rho}_{\theta_{0},\vep})
&=& \frac{1}{6} + 2e_{\vep}(\theta_{0})\theta_{1}
- 2f_{\vep}(\theta_{0})\theta_{0}
+ \frac{b(\theta_{0})}{2} - \theta_{1} \\*
&& + 2\int_{1/6}^{\theta_{0}}\beta_{\vep}(t) \dte t
- \frac{1}{\pi i} \int_{\theta_{1}}^{\theta_{0}}
\Log(L_{\vep}(\theta))\dte \theta \pmod{\Z}.
\end{eqnarray*}
Here
$$
2\int_{1/6}^{\theta_{0}}\beta_{\vep}(t) \dte t =
2f_{\vep}(\theta_{0})\left(\theta_{0}-\frac{1}{6}\right)
+ \frac{1}{\pi i} \int_{\theta_{1}}^{\theta_{0}}
\Log(L_{\vep}(\theta))\dte \theta,
$$
so
\begin{eqnarray*}
\Psi_{n}^{0,b}(2\theta_{0},y)
- \CS(\bar{\rho}_{\theta_{0},\vep})
&=& \frac{1}{\pi i} \left(\int_{1/6}^{\theta_{0}}
\Log(L_{\vep}(\theta))\dte \theta
- \int_{\theta_{1}}^{\theta_{0}}
\Log(L_{\vep}(\theta))\dte \theta\right) \\*
&& + \frac{1}{6} + 2e_{\vep}(\theta_{0})\theta_{1}
- \frac{1}{3}f_{\vep}(\theta_{0})
+ \frac{b(\theta_{0})}{2} - \theta_{1} \pmod{\Z}.
\end{eqnarray*}
{\bf The subcase $\theta_{0} \in ]1/6,1/4[$.} Here we have
$\theta_{1}=1/6$ and $b(\theta_{0})=0$ so the result follows by
the fact that $e_{\vep}(\theta)=f_{\vep}(\theta)$ for
$\theta \in ]1/6,1/4[$.\newline
{\bf The subcase $\theta_{0} \in ]1/4,1/3[$.} In this case we
have $b(\theta_{0})=1$ and $\theta_{1}=1/3$ so
\begin{eqnarray*}
\Psi_{n}^{0,b}(2\theta_{0},y)
- \CS(\bar{\rho}_{\theta_{0},\vep})
&=& \frac{1}{\pi i} \left( \int_{1/6}^{\theta_{0}}
\Log(L_{\vep}(\theta))\dte \theta - \int_{1/3}^{\theta_{0}}
\Log(L_{\vep}(\theta))\dte \theta \right) \\*
&& + \frac{1}{3} + \frac{2}{3}e_{\vep}(\theta_{0})
- \frac{1}{3}f_{\vep}(\theta_{0})  \pmod{\Z}.
\end{eqnarray*}
By (\ref{eq:Lintegral}) we get
$$
\int_{1/3}^{\theta_{0}} \Log(L_{\pm}(t))\dte t =
\int_{1/3}^{\frac{1}{2}-\theta_{0}} \Log(L_{\pm}(t))\dte t
= -\int_{1/6}^{\theta_{0}}
\Log(L_{\pm}\left(\frac{1}{2}-t\right))\dte t.
$$
By (\ref{eq:Lsym2}) and (\ref{eq:Lsym1}) we then have
$$
\int_{1/3}^{\theta_{0}} \Log(L_{\pm}(t))\dte t =
\int_{1/6}^{\theta_{0}} \Log(L_{\pm}(t))\dte t,
$$
so
$$
\Psi_{n}^{0,b}(2\theta_{0},y)
- \CS(\bar{\rho}_{\theta_{0},\vep})
= \frac{1}{3} + \frac{2}{3}e_{\vep}(\theta_{0})
- \frac{1}{3}f_{\vep}(\theta_{0}) \pmod{\Z},
$$
and this is zero by \reflem{lem:Qfunctions}, (\ref{eq:f+}) and
(\ref{eq:f-}).
\end{proofa}

\begin{remark}\label{rem:whybar}
If $u_{\pm}=u_{\pm}(v)$ are the two solutions to
$\phi(v^{2},u)=0$ for $v \in \C^{*}$ fixed, then
$\lambda_{11}(v,u_{-})=\lambda_{11}(v,u_{+})^{-1}$ by
\reflem{lem:lambdasymmetry}. By the proof of
\refthm{thm:comparison} we therefore conclude that
$(v,1+u_{\vep})$, $\vep \in \{\pm\}$, is a solution to
(\ref{eq:stationary1}) and (\ref{eq:stationary2}) if and only
if $(v,1+u_{-\vep})$ is a solution to (\ref{eq:stationary2})
and
\begin{equation}\label{eq:stationary1-}
v^{p} = \left(\frac{w-v^{2}}{1-v^{2}w} \right)^{q}.
\end{equation}
If we work with the invariants $\tau_{r}$ instead of the
invariants $\bar{\tau}_{r}$, then by (\ref{eq:mirror2}) we have
to change $p/q$ to $-p/q$ everywhere in the above. The equation
(\ref{eq:stationary2}) will be the same but
(\ref{eq:stationary1}) will change to (\ref{eq:stationary1-}).
If $(v,u) \in \tilde{\mN}$ then we find as in the proof of
\refthm{thm:comparison} that $(v,w)=(v,u+1)$ is a solution to
(\ref{eq:stationary1-}) if and only if $(v,u)$ is a solution to
$$
v^{p}=\lambda_{11}(v,u)^{q},
$$
which is the ``wrong'' equation. This is one of the main
reasons for working with $\bar{\tau}_{r}$ instead of
$\tau_{r}$. Another reason is that one would get a minus sign
in (\ref{eq:mainstatement}).
\end{remark}

Let us end this section by some futher results on the
Chern--Simons invariants of flat irreducible
$\SU(2)$--connections on $M_{p/q}$.
By (\ref{eq:Lsym1}) $\bar \rho _{\theta,\vep}$ extends to
a $\SU(2)$--representation of $\pi_1(M_{p/q})$ if and only if
$\bar \rho _{\theta,-\vep}$ extends to
a $\SU(2)$--representation of $\pi_1(M_{-p/q})$. In that case
we find that
$$
\CS(\bar \rho _{\theta,-\vep}) =
- \CS(\bar \rho _{\theta,\vep})
$$
by using (\ref{eq:p/q'}).

If $p$ is even then
$e^{2\pi ip(1/2-\theta)}=e^{2\pi ip\theta}$. Since
$L_{\vep}(1/2-\theta)=L_{\vep}(\theta)^{-1}$ by
(\ref{eq:Lsym2}) and (\ref{eq:Lsym1}) and since $q$ is odd it
follows that $\bar \rho _{\theta,\vep}$ extends to a
$\SU(2)$--representation of $\pi_1(M_{p/q})$ if and only if
$\bar \rho _{1/2-\theta,\vep}$ extends to such a
representation, see (\ref{eq:p/q}). Assuming that
$\bar \rho _{\theta,\vep}$ extends to such a representation
we get
$$
\CS(\bar \rho _{\frac{1}{2}-\theta,\vep})
= \CS(\bar \rho _{\theta,\vep})
- \frac{1}{2} \frac{p}{2} \pmod{\Z}
$$
by using (\ref{eq:Lintegral}) and (\ref{eq:p/q'}).

By (\ref{eq:Lspecialvalues}) and (\ref{eq:p/q}) we find that
$\bar \rho _{1/6,-}=\bar \rho _{1/6,+}$ is a
$\SU(2)$--representation of $\pi_1(M_{p/q})$ if and only if
$p=6m+3$, $m \in \Z$, and q is odd. In that case
$$
\CS(\bar \rho _{1/6,\vep}) = - \frac{cp}{36} - \frac{dq}{4}
- \frac{1}{2}d \pmod{\Z}.
$$
By (\ref{eq:Lspecialvalues}) and (\ref{eq:p/q}) we also have
that $\bar \rho _{1/3,-}=\bar \rho _{1/3,+}$ is a
$\SU(2)$--representation of $\pi_1(M_{p/q})$ if and only if
$p=6m+3$, $m \in \Z$, and q is even. In that case
$$
\CS(\bar \rho _{1/3,\vep}) = \frac{1}{2} - \frac{cp}{9}
- \frac{dq}{4} \pmod{\Z}.
$$

Finally, by (\ref{eq:Lspecialvalues}) and (\ref{eq:p/q}),
$\rho_{1/4,-}$ and $\rho_{1/4,+}$ are
$\SU(2)$--representations of $\pi_1(M_{p/q})$ if and only if
$4$ divides $p$. In that case we find that
\begin{eqnarray*}
\CS(\bar \rho _{1/4,+}) &=& \frac{1}{5}
- \frac{cp}{16} \pmod{\Z}, \\
\CS(\bar \rho _{1/4,-}) &=& - \frac{1}{5}
- \frac{cp}{16} \pmod{\Z}.
\end{eqnarray*}
This follows by using (\ref{eq:Specialintegral1/4}).

\renewcommand{\theequation}{\Alph{appendixc}.\arabic{equation}}

\noindent We have in the following appendices collected
material of a technical nature.

\appendix{\,\,\,Proofs of \reflem{lem:functional},
\reflem{lem:Sapproximation} and (\ref{eq:Sgdilog})}

\tenrm\baselineskip=13pt

\begin{proofa}{\it of\,} \reflem{lem:functional}\,\,Let $a>0$.
Let $\vep =1$ if $\ima(\zeta) \geq 0$ and let $\vep=-1$
otherwise. Put $\delta_{a}^{-}=[-a,\I \vep a]$ and
$\delta_{a}^{+}=[\I\vep a,a]$. (Here, as usual, $[z_{1},z_{2}]$
denotes the line segment in $\C$ beginning at $z_{1}$ and
ending at $z_{2}$.) We have
$$
\frac{S_{\gamma}(\zeta-\gamma)}{S_{\gamma}(\zeta+\gamma)} =
\exp \left( - \frac{1}{2} \int_{C_{R}}
\frac{e^{\zeta z}}{\sinh(\pi z) z} \dte z \right).
$$
By an elementary argument one finds that the integrals
$\int_{\delta_{a}^{\pm}}
\frac{e^{\zeta z}}{\sinh(\pi z) z} \dte z$
converge to zero as $a \ria \infty$. Therefore
$$
\int_{C_{R}} \frac{e^{\zeta z}}{\sinh(\pi z) z} \dte z =
\vep 2\pi \I \left( b_{\vep}+ \sum_{n=1}^{\infty}
\Res_{z=\vep\I n} \left\{ \frac{e^{\zeta z}}{\sinh(\pi z) z}
\right\} \right),
$$
where $b_{1}=0$ and
$b_{-1}=\Res_{z=0} \left\{
\frac{e^{\zeta z}}{\sinh(\pi z) z} \right\}=\frac{\zeta}{\pi}$.
For $n \in \Z \sm \{0\}$ we have
$$
\Res_{z=\I n}
\left\{ \frac{e^{\zeta z}}{\sinh(\pi z) z} \right\}
=\frac{(-1)^{n}e^{\I \zeta n}}{\pi\I n},
$$
so
$$
\int_{C_{R}} \frac{e^{\zeta z}}{\sinh(\pi z) z} \dte z
= -(1-\vep)\I \zeta -2\Log\left( 1+e^{\vep\I \zeta } \right)
$$
giving the result.
\end{proofa}

To prove the identity (\ref{eq:Sgdilog}) we use the power
series expansion
\begin{equation}\label{eq:dilogseries}
\Li_{2}(z) = \sum_{n=1}^{\infty} \frac{z^{n}}{n^{2}}
\end{equation}
for the dilogarithm, valid for $|z| \leq 1$. In the course of
the following proof we will establish the identity
\begin{equation}\label{eq:dilogid1}
\frac{\zeta^{2}}{2} - \frac{\pi^{2}}{6}
- \Li _{2} (-e^{-\I \zeta} )
= \Li _{2} (-e^{\I \zeta} )
\end{equation}
valid for $\zeta = \pm \pi$ and all $\zeta \in \C$ with
$|\rea(\zeta)| < \pi$.

\begin{proofa}{\it of\,} (\ref{eq:Sgdilog})\,\,
Note first that the integral
$A_{\gamma}(\zeta) := \frac{1}{4} \int_{C_{R}}
\frac{e^{\zeta z}}{\sinh(\pi z) \gamma z^{2}} \dte z$
is convergent for all $\zeta \in \C$ with
$|\rea(\zeta)| \leq \pi$ since
$$
\int_{-\infty}^{-R}
\frac{e^{\zeta t}}{\sinh(\pi t) t^{2}} \dte t
= - \int_{R}^{\infty}
\frac{e^{-\zeta t}}{\sinh(\pi t) t^{2}} \dte t
$$
and
$$
\left| \int_{R}^{\infty}
\frac{e^{\zeta t}}{\sinh(\pi t) t^{2}} \dte t \right|
\leq \frac{2}{1-e^{-2\pi R}} \int_{R}^{\infty}
e^{(\rea(\zeta)-\pi)t} \frac{1}{t^{2}} \dte t.
$$
Let $b=\sign(\ima(\zeta))$, where $\sign(0)$ can be put to both
$1$ and $-1$ in the following. Moreover, let $h$ be a positive
parameter and let
$\delta_{h}^{-}(t) = (1 + ib)ht + ibh$ for $t \in [-1,0]$ and
let $\delta_{h}^{+}(t) = (1 - ib)ht + ibh$ for $t \in [0,1]$.
It is elementary to show that the integrals
$\int_{\delta_{h}^{\pm}}
\frac{e^{\zeta z}}{\sinh(\pi z) z^{2}} \dte z$
converge to zero as $h$ converges to infinity for
$|\rea(\zeta)| \leq \pi$. By the residue theorem we conclude
that
$$
A_{\gamma}(\zeta) = \frac{1}{4\gamma}
2\pi i \sum_{n=1}^{\infty} \Res_{z=in}
\left\{ \frac{e^{\zeta z}}{\sinh(\pi z) z^{2}}\right\}
$$
for $\ima(\zeta) \geq 0$ and $|\rea(\zeta)| \leq \pi$ and
$$
A_{\gamma}(\zeta) = -\frac{1}{4\gamma}
2\pi i \sum_{n=0}^{\infty} \Res_{z=-in}
\left\{ \frac{e^{\zeta z}}{\sinh(\pi z) z^{2}}\right\}
$$
for $\ima(\zeta) \leq 0$ and $|\rea(\zeta)| \leq \pi$. Using
(\ref{eq:dilogseries}) this leads directly to
$$
A_{\gamma}(\zeta) = \left\{ \begin{array}{ll}
\frac{1}{2\I \gamma}
\Li _{2} (-e^{\I \zeta} ) & ,\ima(\zeta) \geq 0, \\
\frac{1}{2\I \gamma}
\left[ \frac{\zeta^{2}}{2} - \frac{\pi^{2}}{6}
- \Li _{2} (-e^{-\I \zeta} ) \right] & ,\ima(\zeta)
\leq 0. \end{array} \right.
$$
for $|\rea(\zeta)| \leq \pi$. Left is to prove the identity
(\ref{eq:dilogid1}). To this end, let
$$
g(\zeta) = \frac{\zeta^{2}}{2} - \frac{\pi^{2}}{6}
- \Li_{2}(-e^{-i\zeta}) - \Li_{2}(-e^{i\zeta})
$$
for
$\zeta \in \Omega :=
\{ \zeta \in \C \; | \; |\rea(\zeta)| < \pi\}$.
By (\ref{eq:dilog}) we have
$$
g'(\zeta) = \zeta + i\left( \Log(1+e^{i\zeta})
- \Log(1+e^{-i\zeta})\right)
$$
and therefore $e^{ig'(\zeta)}=1$. Since $\Omega$ is connected,
$g$ is $C^{1}$ and $g'(0)=0$ we get that $g'$ is identically
zero on $\Omega$ so $g$ is constant on $\Omega$. Now
$g(0) = -\frac{\pi^{2}}{6} - 2\Li_{2}(-1)$ and
$$
\Li_{2}(-1) = \sum_{n=1}^{\infty} \frac{(-1)^{n}}{n^{2}}
= -\frac{\pi^{2}}{12}
$$
so $g(0)=0$. Finally we note that $g$ is well-defined and
continuous on $\Omega \cup \{ \pm \pi \}$ so $g(\pm \pi) =0$ by
continuity.
\end{proofa}

Note that the function $g$ in the above proof is a well-defined
analytic function on
$W = \{ \; \zeta \in \C \mid \rea(\zeta)
\notin \pi + 2\pi \Z \; \}$
and that $g$ is continuous on $W \cup \R$. As in the proof
above we find that $g$ is constant on each connected component
of $W$. Moreover, we can on each of these connected components
choose a branch of the dilogarithm such that $g$ extends to a
continuous (and hence a constant) function on the connected set
$W \cup \R$.

\begin{proofa}{\it of\,} \reflem{lem:Sapproximation}\,\,
The function $1/\sinh(w)$ has a simple pole at $w=0$ with
principal part $1/w$, i.e.\
$$
\phi(w)= \frac{1}{\sinh(w)} - \frac{1}{w}
$$
is holomorphic in a neighborhood of zero, in fact on the open
disk $D(0,\pi)$ with centre $0$ and radius $\pi$. Let
$a \geq R$, let $C_{R;a}=[-a,-R] \cup \Upsilon_{R} \cup [a,R]$,
and write
$I_{\gamma}(\zeta) = J_{\gamma}(\zeta) + K_{\gamma}(\zeta)$,
where
\begin{eqnarray*}
4J_{\gamma}(\zeta) &=& \int_{C_{R;a}}
\frac{e^{\zeta z}}{\sinh(\pi z) z} \phi(\gamma z) \dte z, \\
4K_{\gamma}(\zeta) &=& \int_{-\infty}^{-a}
\frac{e^{\zeta z}}{\sinh(\pi z) z} \phi(\gamma z) \dte z
+ \int_{a}^{\infty}
\frac{e^{\zeta z}}{\sinh(\pi z) z} \phi(\gamma z) \dte z.
\end{eqnarray*}
To estimate $K_{\gamma}(\zeta)$ we simply use that
$$
|\phi(\gamma t)| \leq \frac{1}{\sinh(\gamma t)}
+ \frac{1}{\gamma t} \leq \frac{2}{\gamma t}
$$
for $t>0$ leading to the bound
\begin{eqnarray*}
|K_{\gamma}(\zeta)| &\leq& \frac{1}{2\gamma} \int_{a}^{\infty}
\frac{e^{\rea(\zeta)t} + e^{-\rea(\zeta)t}}{\sinh(\pi t)}
\frac{1}{t^{2}} \dte t \\*
&\leq& \frac{1}{\gamma (1-e^{-2\pi a})} \int_{a}^{\infty}
\left( e^{-(\pi-\rea(\zeta))t} + e^{-(\pi + \rea(\zeta))t}
\right) \frac{1}{t^{2}} \dte t.
\end{eqnarray*}
For $|\rea(\zeta)| \leq \pi$ we therefore get
$$
|K_{\gamma}(\zeta)| \leq \frac{2}{\gamma (1-e^{-2\pi a})}
\int_{a}^{\infty} \frac{1}{t^{2}} \dte t =
\frac{2}{a\gamma (1-e^{-2\pi a})}.
$$
If $|\rea(\zeta)| < \pi$ we find that
\begin{eqnarray*}
|K_{\gamma}(\zeta)| &\leq&
\frac{1}{\gamma a^{2} (1-e^{-2\pi a})}
\int_{a}^{\infty} (e^{-(\pi-\rea(\zeta)) t}
+ e^{-(\pi+\rea(\zeta))t}) \dte t \\*
&=& \frac{1}{\gamma a^{2} (1-e^{-2\pi a})}
\left( \frac{e^{-(\pi-\rea(\zeta))a}}{\pi - \rea(\zeta)}
+ \frac{e^{-(\pi+\rea(\zeta))a}}{\pi + \rea(\zeta)} \right).
\end{eqnarray*}
Next let us estimate $J_{\gamma}(\zeta)$. First we use the
standard estimate
$$
\left| \int_{\Upsilon_{R}} \frac{e^{\zeta z}}{\sinh(\pi z) z}
\phi(\gamma z) \dte z \right| \leq \pi R M(\zeta,R),
$$
where
$M(\zeta,R) = \max_{z \in \Upsilon_{R}}
\left| \frac{e^{\zeta z}}{\sinh(\pi z) z} \phi(\gamma z)
\right|.$
We have
$$
\left| \phi(w) \right| = \frac{\sinh(w)-w}{w\sinh(w)}.
$$
Here $\sinh(w) - w = w^{3}h(w)$ and $w\sinh(w) = w^{2}k(w)$,
where $h$ and $k$ are entire functions. Note that $k$ is
different from zero on $D(0,\pi)$. Since $\gamma \in ]0,1[$ we
get
$$
M(\zeta,R) = 2\gamma L(R) N(\zeta,R),
$$
where $L(R) = \max_{|z| \leq R} |h(z)/k(z)|$ and
$N(\zeta,R)=\max_{z \in \Upsilon_{R}}\left|
\frac{e^{\zeta z}}{e^{\pi z}-e^{-\pi z}} \right|.$
We note that
$$
N(\zeta,R) \leq Q_{\pm}(R) \max_{z \in \Upsilon_{R}}
\left| e^{(\zeta \pm \pi)z} \right|,
$$
where
$Q_{\pm}(R) = \max_{z \in \Upsilon_{R}}
\frac{1}{\left| 1 - e^{\pm 2\pi z} \right|}$.
Put $Q(R)=Q_{-}(R)+Q_{+}(R)$ and get
$$
N(\zeta,R) \leq Q(R) \min_{\mu = \pm 1}
\left( \max_{z \in \Upsilon_{R}}
\left| e^{(\zeta +\mu \pi)z} \right|\right).
$$
Since $\rea(z) \in [-R,R]$ and $\ima(z) \in [0,R]$ for
$z \in \Upsilon_{R}$ we finally get
$$
N(\zeta,R) \leq Q(R) e^{2\pi R}
\left(1 + e^{-\ima(\zeta)R} \right).
$$
We have thus obtained the estimate
$$
\left| \frac{1}{4}\int_{\Upsilon_{R}}
\frac{e^{\zeta z}}{\sinh(\pi z) z} \phi(\gamma z) \dte z
\right| \leq \gamma B \left(1 + e^{-\ima(\zeta)R} \right),
$$
where $B=\frac{\pi}{2} RL(R)Q(R)e^{2\pi R}$.

Finally we have to estimate
$\int_{R}^{a} \frac{e^{\zeta t}-e^{-\zeta t}}{\sinh(\pi t)t}
\phi(\gamma t) \dte t$.
First observe that
$$
h(y) = \sum_{n=0}^{\infty} \frac{y^{2n}}{(2n+3)!}
$$
and
$$
k(y) = \sum_{n=0}^{\infty} \frac{y^{2n}}{(2n+1)!}
$$
so $h(y) \leq k(y)/6$ for $y \in \R$. Therefore
\begin{eqnarray*}
\left| \int_{R}^{a}
\frac{e^{\zeta t}-e^{-\zeta t}}{\sinh(\pi t)t}
\phi(\gamma t) \dte t \right|
&\leq& \frac{\gamma}{3} \int_{R}^{a}
\frac{e^{\rea(\zeta)t}+e^{-\rea(\zeta)t}}{e^{\pi t}-e^{-\pi t}}
\dte t \\*
&\leq& \frac{\gamma}{3(1-e^{-2\pi R})} \int_{R}^{a}
\left(e^{-(\pi-\rea(\zeta))t}+e^{-(\pi+\rea(\zeta))t}\right)
\dte t .
\end{eqnarray*}
For $|\rea(\zeta)| \leq \pi$ we get
$$
\left| \int_{R}^{a}
\frac{e^{\zeta t}-e^{-\zeta t}}{\sinh(\pi t)t} \phi(\gamma(t)
\dte t \right| \leq \frac{2a\gamma}{3(1-e^{-2\pi R})}.
$$
If $|\rea(\zeta) <\pi$ we get
\begin{eqnarray*}
&& \left| \int_{R}^{a}
\frac{e^{\zeta t}-e^{-\zeta t}}{\sinh(\pi t)t}
\phi(\gamma(t) \dte t \right| \leq
\frac{\gamma}{3(1-e^{-2\pi R})} \int_{0}^{a}
\left(e^{-(\pi-\rea(\zeta))t}+e^{-(\pi+\rea(\zeta))t}\right)
\dte t \\*
 &&\hspace{.2in} = \frac{\gamma}{3(1-e^{-2\pi R})}
\left( \frac{1 - e^{-(\pi-\rea(\zeta))a}}{\pi-\rea(\zeta)}
+ \frac{1 - e^{-(\pi+\rea(\zeta))a}}{\pi+\rea(\zeta)} \right).
\end{eqnarray*}
The lemma now follows by putting $a=1/\gamma >1>R$ and
$A = \frac{13}{12(1-e^{-2\pi R})}$.
\end{proofa}

\appendix{\,\,\,Proofs of the estimates
(\ref{eq:figure8estimate1}) and (\ref{eq:figure8estimate2})}

\noindent Let
\begin{eqnarray*}
J_{\pm}^{'}(r,\vep) &=& \int_{C_\pm(\vep)}
\left(\exp\left( I_{\gamma}(\pi-2\pi x)
- I_{\gamma}(-\pi + 2\pi x)\right) -1 \right)
e^{r\Phi(x)} \dte x, \\
J_{\pm}^{''}(r,\vep) &=&
\int_{C_\pm(\vep)}(\tan(\pi r x) \mp \I)g_r(x) \dte x,
\end{eqnarray*}
where $\Phi$ is given by (\ref{eq:figure8Phase}). Note first
that we are free to deform the contour $C_{\pm}(\vep)$ as long
as we stay inside the domain of analyticity of the integrands.
For the integrals $J_{\pm}^{'}(r,\vep)$ we will deform
$C_{\pm}(\vep)$ to $\mp [\vep,1-\vep]$. Since the integrands of
the integrals $J_{\pm}^{''}(r,\vep)$ are analytic on
$\Omega_{\frac{1}{2}r} \sm \{ \, (m+1/2)/r \, |\,
m=0,1,\ldots,r-1\,\}$,
we can deform $C_{\pm}(\vep)$ to $C_{\pm}(0)$ in these
integrals without changing their sum, i.e.\
$$
J_{-}^{''}(r,\vep) + J_{+}^{''}(r,\vep) = J_{-}^{''}(r)
+ J_{+}^{''}(r),
$$
where $J_{\pm}^{''}(r) = J_{\pm}^{''}(r,0)$. In the following
calculations we will need the identity
\begin{equation}\label{eq:dilogrea1}
\rea\left( \Phi(x) \right) = 2\pi\rea(x)\ima(x) - \pi\ima(x)
- \frac{1}{\pi} \ima\left( \Li_{2}
\left( e^{-2\pi\ima(x)} e^{2\pi i\rea(x)} \right) \right)
\end{equation}
valid for $x \in \Omega_{\infty}$ (see
(\ref{eq:typicaldomain})). This identity is an immediate
consequence of (\ref{eq:dilogid1}).

\begin{proofa}{\it of\,} (\ref{eq:figure8estimate1})\,\,
Let us first estimate $J_{+}^{''}(r)$. We partition
$C_{+}(0)$ into the three pieces $C_{+}^{1} = [\I,0]$,
$C_{+}^{2} = [1+\I,\I]$, and $C_{+}^{3} = [1,1+\I]$. Put
$$
I_{+}^{i}(r) = \int_{C_{+}^{i}}
(\tan(\pi r x) - \I)g_{r}(x) \dte x.
$$
By (\ref{eq:tan+}) we immediately get
$$
|I_{+}^{1}(r)| \leq
2\int_{0}^{1} e^{-2\pi r t} |g_{r}(\I t)| \dte t.
$$
To be able to use (\ref{eq:Sgdilog}) we introduce the positive
parameter $\vep$ again. In fact, we have by (\ref{eq:Sgdilog}),
\reflem{lem:Sapproximation} and Lebesgue's dominated
convergence theorem that
$$
|I_{+}^{1}(r)| \leq 2\exp(2A + 2B\pi/r) \livep \int_{0}^{1}
e^{-2\pi r t} \left| e^{r\Phi(\vep+ \I t)} \right| \dte t.
$$
By (\ref{eq:dilogrea1}) we immediately get that
$$
\left| e^{r\Phi(\vep+ \I t)} \right| =
\exp\left( r \left[ 2\pi\vep t - \pi t -\frac{1}{\pi}
\ima\left(\Li_{2}
\left( e^{2\pi \sqrt{-1} \vep} e^{-2\pi t} \right) \right)
\right] \right).
$$
Now by continuity of
$(t,\vep) \mapsto \Li_{2} \left( e^{2\pi \sqrt{-1} \vep}
e^{-2\pi t} \right)$
on $[0,1] \times [0,1]$ we can remove the parameter $\vep$
again by using Lebusgue's dominated convergence theorem once
more. This gives us
$$
\livep \int_{0}^{1} e^{-2\pi r t} \left| e^{r\Phi(\vep+ \I t)}
\right| \dte t = \int_{0}^{1} e^{-3\pi rt} \dte t
$$
leading to the estimate
$$
|I_{+}^{1}(r)| \leq \frac{2}{3\pi r}
\exp(2A+2B\pi/r)\left(1-e^{-3\pi r}\right).
$$
Next we estimate $I_{+}^{2}$. By (\ref{eq:tan+}) we get
$$
|I_{+}^{2}(r)| \leq 4e^{-2\pi r}
\int_{0}^{1} |g_{r}(\I+t)| \dte t.
$$
Similarly to the analysis of $I_{+}^{1}$ we introduce the
parameter $\vep$ and get
$$
|I_{+}^{2}(r)| \leq 4\exp(2A+2B\pi/r) e^{-2\pi r}
\livep \int_{\vep}^{1-\vep} |e^{r\Phi(\I+t)}| \dte t.
$$
Here
$$
|e^{r\Phi(\I+t)}| = e^{-\pi r(1-2t)}
\exp\left(-\frac{r}{\pi}\ima
\left(\Li_{2}(e^{-2\pi}e^{2\pi\I t})\right) \right)
$$
by (\ref{eq:dilogrea1}), so by Lebesgue's dominated convergence
theorem we get
$$
|I_{+}^{2}(r)| \leq 4\exp(2A+2B\pi/r) e^{-3\pi r}
\int_{0}^{1} e^{2\pi rt}
\exp\left(-\frac{r}{\pi}\ima \left(
\Li_{2}(e^{-2\pi}e^{2\pi\I t})\right) \right) \dte t.
$$
By definition of the dilogarithm we have
$$
\ima \left(\Li_{2}(e^{-2\pi}e^{2\pi\I t})\right) =
-\int_{0}^{1} \frac{
\Arg\left(1-se^{-2\pi}e^{2\pi\I t}\right)}{s} \dte s,
$$
which is non-negative for $t \in [0,1/2]$. For $t \in [1/2,1]$
we use that
$$
\ima \left(\Li_{2}(e^{-2\pi}e^{2\pi\I t})\right) =
\ima \left(\Li_{2}(e^{2\pi\I t})\right)
+ \int_{1-e^{-2\pi}}^{1}
\frac{\Arg\left(1-se^{2\pi\I t}\right)}{s} \dte s,
$$
where the last integral is positive. The first term is bounded
from below by
$-\sum_{n=1}^{\infty} \frac{1}{n^{2}} = -\frac{\pi^{2}}{6}$ by
(\ref{eq:dilogimacircle}). We therefore end up with
\begin{eqnarray*}
|I_{+}^{2}(r)| &\leq& 4\exp(2A+2B\pi/r) e^{-(3-1/6)\pi r}
\int_{0}^{1} e^{2\pi rt}\dte t \\*
&=& \frac{2}{\pi r} \exp(2A+2B\pi/r)
e^{-(3-1/6)\pi r}\left(e^{2\pi r} -1 \right).
\end{eqnarray*}
Finally, we estimate $I_{+}^{3}$. Similarly to the other cases
we get
$$
|I_{+}^{3}(r)| \leq 2\exp(2A+2B\pi/r) \livep \int_{0}^{1}
e^{-2\pi rt} \left| e^{r\Phi(1-\vep + \I t)} \right| \dte t.
$$
By (\ref{eq:dilogrea1}) we have
$$
\left| e^{r\Phi(1-\vep + \I t)} \right| =
\exp\left( r \left[ \pi(1-2\vep)t
- \frac{1}{\pi} \ima \left(
\Li_{2}(e^{-2\pi t}e^{-2\pi\I \vep})\right) \right] \right)
$$
leading to
$$
|I_{+}^{3}(r)| \leq 2\exp(2A+2B\pi/r)
\int_{0}^{1} e^{-\pi rt} \dte t =
\frac{2}{\pi r}\exp(2A+2B\pi/r)\left(1-e^{-\pi r}\right).
$$
By letting $C_{-}^{1} = [-\I,0]$, $C_{-}^{2} = [1-\I,-\I]$, and
$C_{-}^{3} = [1,1-\I]$ and
$$
I_{-}^{i}(r) =
\int_{C_{-}^{i}} (\tan(\pi r x) + \I)g_{r}(x)\dte x
$$
we find (now by the use of (\ref{eq:tan-})) an upper bound for
$|I_{-}^{i}(r)|$ identical with the upper bound for
$|I_{+}^{i}(r)|$, $i=1,2,3$, with the exception that
$\exp(2A+2B\pi/r)$ should be replaced by $\exp(2A+4B\pi/r)$ in
these bounds. To conclude we have shown that there exists a
constant $K_{1}$ independent of $r$ and $\vep$ such that
$$
|J_{+}^{''}(r,\vep) + J_{-}^{''}(r,\vep)| \leq
\sum_{i=1}^{3} \left( |I_{+}^{i}| + |I_{-}^{i}| \right) \leq
\frac{K_{1}}{r}
$$
for all $r \in \Z_{\geq 2}$.
\end{proofa}

\begin{proofa}{\it of\,} (\ref{eq:figure8estimate2})\,\, By the
remarks prior to the proof of (\ref{eq:figure8estimate1}) we
have
$$
J_{\pm}^{'}(r,\vep) = \mp \int_{\vep}^{1-\vep}
\exp\left(r\Phi(t) \right) h_{\gamma}(t) \dte t,
$$
where
$$
h_{\gamma}(t) = \exp\left( I_{\gamma}(\pi-2\pi t)
- I_{\gamma}(-\pi + 2\pi t) \right) -1.
$$
The integrand is continuous on $[0,1]$. Therefore
$$
\left| J_{\pm}^{'}(r,\vep) \right| \leq \int_{0}^{1}
\exp\left(r\rea\left(\Phi(t)\right) \right)
\left| h_{\gamma}(t) \right| \dte t
$$
for all $\vep \in ]0,\frac{1}{4r}[$. By
(\ref{eq:dilogimacircle}) and (\ref{eq:volume}) and the remarks
prior to (\ref{eq:volume}) we have
$$
\rea\left( \Phi(t) \right) = \frac{1}{2\pi} \ima
\left( \Li_{2}(e^{-2\pi\I t}) - \Li_{2}(e^{2\pi\I t}) \right)
= -\frac{1}{\pi}\Cl_{2}(2\pi t) \leq \frac{1}{2\pi}\Vol(K).
$$
Therefore
$$
|J_{\pm}^{'}(r,\vep)| \leq
\exp\left( \frac{r}{2\pi}\Vol(K)\right)
\int_{0}^{1} \left|h_{\gamma}(t)\right|\dte t.
$$
By definition we have
$$
h_{\gamma}(t) = \sum_{n=1}^{\infty} \frac{1}{n!}
\left( I_{\gamma}(\pi-2\pi t) -
I_{\gamma}(-\pi+2\pi t)\right)^{n}.
$$
From \reflem{lem:Sapproximation} we get
$$
\left| I_{\gamma}(\pi-2\pi t) -
I_{\gamma}(-\pi+2\pi t) \right| \leq
\left(C f(t) + D\right) \gamma
$$
for $t \in ]0,1[$, where $C$ and $D$ are positive constants
independent of $\gamma$ and $t$, and
$$
f(t)=\frac{1}{t}+\frac{1}{1-t}.
$$
Since $f : ]0,1[ \to \R$ is bigger than or equal to $4$ we can
choose $C$ so big that
$$
\left| I_{\gamma}(\pi-2\pi t) - I_{\gamma}(-\pi+2\pi t) \right|
\leq C f(t) \gamma
$$
for $t \in ]0,1[$. From \reflem{lem:Sapproximation} we also
have
$$
\left| h_{\gamma}(t) \right| \leq
\exp\left( |I_{\gamma}(\pi-2\pi t)
-I_{\gamma}(-\pi+2\pi t)| \right) \leq \exp(4A+4B\pi/r)
$$
for $t \in [0,1]$, where $A$ and $B$ are as in
\reflem{lem:Sapproximation}, so
$\frac{\pi}{r} \exp(4A+4B\pi/r)$ is an upper bound for both of
the integrals $\int_{0}^{\gamma} |h_{\gamma}(t)|\dte t$ and
$\int_{1-\gamma}^{1} |h_{\gamma}(t)|\dte t$ (for $r \geq 4$).
Left is to evaluate (for $r \geq 7$)
$$
\int_{\gamma}^{1-\gamma} |h_{\gamma}(t)|\dte t \leq
\sum_{n=1}^{\infty} \frac{C^{n}\gamma^{n}}{n!}
\int_{\gamma}^{1-\gamma} f(t)^{n} \dte t.
$$
By using that
$$
f(t)^{n} = \sum_{k=0}^{n} \left(
\begin{array}{c} n \\ k \end{array} \right)
\frac{1}{t^{k}}\frac{1}{(1-t)^{n-k}}
$$
we get
\begin{eqnarray*}
\int_{\gamma}^{1-\gamma} f(t)^{n}\dte t &=&
\sum_{k=0}^{n} \left(
\begin{array}{c} n \\ k \end{array} \right)
\int_{\gamma}^{1-\gamma} \frac{1}{t^{k}}
\frac{1}{(1-t)^{n-k}} \dte t \\
&\leq & \sum_{k=0}^{n} \left(
\begin{array}{c} n \\ k \end{array} \right) \left( 2^{n-k}
\int_{\gamma}^{1/2} \frac{1}{t^{k}}\dte t + 2^{k}
\int_{1/2}^{1-\gamma} \frac{1}{(1-t)^{n-k}} \dte t \right) \\
&\leq & 2^{n} \sum_{k=0}^{n} \left(
\begin{array}{c} n \\ k \end{array} \right) \left(
\int_{\gamma}^{1/2} \frac{1}{t^{k}}\dte t +
\int_{\gamma}^{1/2} \frac{1}{t^{n-k}} \dte t \right) \\
&\leq & 2^{n+1} \left(\sum_{k=0}^{n} \left(
\begin{array}{c} n \\ k \end{array} \right) \right)
\int_{\gamma}^{1/2} \frac{1}{t^{n}}\dte t
= 2^{2n+1} \int_{\gamma}^{1/2} \frac{1}{t^{n}}\dte t.
\end{eqnarray*}
Here
$$
\int_{\gamma}^{1/2} \frac{1}{t}\dte t =
-\Log(2) - \Log(\gamma) \leq \Log(r)
$$
and
$$
\int_{\gamma}^{1/2} \frac{1}{t^{n}}\dte t = \frac{1}{n-1}
\left( \frac{1}{\gamma^{n-1}} - 2^{n-1}\right)
\leq \frac{1}{n-1} \frac{1}{\gamma^{n-1}}
$$
for $n \geq 2$. Therefore
\begin{eqnarray*}
\int_{\gamma}^{1-\gamma} |h_{\gamma}(t)|\dte t &\leq&
2\gamma\left( 4C\Log(r)
+ \sum_{n=2}^{\infty} \frac{(4C)^{n}}{(n-1)n!} \right) \\*
&\leq & \frac{2\pi}{r}\left(4C\Log(r) + \exp(4C)-4C-1 \right).
\end{eqnarray*}
We conclude that there exists a constant $K_{2}$ independent of
$r$ and $\vep$ such that
$$
|J_{+}^{'}(r,\vep) + J_{-}^{'}(r,\vep)| \leq K_{2}
\frac{\Log(r)}{r} \exp\left( \frac{r}{2\pi}\Vol(K) \right),
$$
for all $r \in \Z_{\geq 2}$.
\end{proofa}

\appendix{\,\,\,The case $M_{0}$}

\noindent The manifold $M_{0}$ is the mapping torus of a torus
with monodromy matrix
$$
\left(
\begin{array}{cc}
 2 & 1 \\
 1 & 1 \\
\end{array} \right),
$$
see \cite[p.~366]{[18]}. The invariant $\tau_{r}(M_{0})$ has
been calculated by Jeffrey \cite[Theorem 4.1]{[13]}. This
theorem gives the large $r$ asymptotics of the invariant as
well. In fact, we have
\begin{equation}\label{eq:Jeffrey}
\tau_{r}(M_{0}) = \frac{1}{2} - \frac{1}{2\sqrt{5}}
- \frac{1}{\sqrt{5}}
\left( \exp\left(2\pi\I r \left(-\frac{1}{5}\right)\right)
+ \exp\left(2\pi\I r \frac{1}{5}\right) \right)
\end{equation}
which at the same time can be taken as the large $r$
asymptotics of the invariant. Let us relate this result to our
contour integral formula for the invariant
$\tau_{r}(M_{0}) = \bar{\tau}_{r}(M_{0})$. By
\reflem{lem:contour}
$$
\tau_{r}(M_{0}) = \frac{ri}{4}
\int_{C_{r}^{1}} \cot(\pi r x)
\left( \int_{C_{r}^{2}} \tan(\pi ry) f_{0,r}( x,y) \dte y
\right) \dte x,
$$
where
$$
f_{0,r}(x,y) = \sin(\pi x) e^{-2\pi i r xy}
\frac{S_{\pi/r}(-\pi+2\pi(x - y))}
{S_{\pi/r}(-\pi+2\pi(x + y))}.
$$
Following the discussion in Sec.~4.3 the relevant (shifted)
phase functions to consider are given by
$$
\Psi_{0}^{a,b}(x,y) = ax + by -xy + \frac{1}{4\pi^{2}}\left(
\Li_{2}(e^{2\pi i(x+y)}) - Li_{2}(e^{2\pi i(x-y)}) \right),
$$
where $a,b$ are certain integers. If we put
$\Psi=\Psi_{0}^{a,b}$, then by (\ref{eq:critical}) $(x,y)$ is a
critical point of $\Psi$ if and only if
\begin{eqnarray}\label{eq:criticalpzero}
0 &=& 2\pi i(a-y) + \Log(1-zw)- \Log(1-zw^{-1}), \nonumber \\*
0 &=& 2\pi i (b-x) + \Log(1-zw)+\Log(1-zw^{-1}),
\end{eqnarray}
where $z=e^{2\pi ix}$ and $w=e^{2\pi iy}$ as usual, and this
set of equations implies that
\begin{eqnarray}\label{eq:stationarypzero}
w-z &=& 1-zw, \nonumber \\*
(1-zw)(w-z) &=& zw,
\end{eqnarray}
compare with (\ref{eq:stationary1}) and (\ref{eq:stationary2}).
We note that the first of these equations is equivalent to
$$
w-1 = z(1-w)
$$
so $w=1$ or $z=-1$. For $z=-1$ we get $w^{2}+3w+1=0$ so
$w=\frac{-3 \pm \sqrt{5}}{2}$. For $w=1$ we find that
$z^{2}-3z+1=0$ so $z=\frac{3 \pm \sqrt{5}}{2}$. However, only
the point $(z,w)=((3-\sqrt{5})/2,1)$ satisfies that
$zw,zw^{-1} \notin [1,\infty[$. For this point we find that
$y \in \Z$ and from (\ref{eq:criticalpzero}) we get that $y=a$.
But then
$$
\Psi(x,y) = ax + by - xy + \frac{1}{4\pi^{2}}\left(
\Li_{2}(z)-\Li_{2}(z)\right) =  yx + ba - xy = ba \in \Z.
$$
Let $\bar{\rho}_{\theta,\vep}$ be the nonabelian
$\SU(2)$--representations of $\pi$ from \refprop{prop:SU(2)},
where $\theta \in [-1/3,-1/6] \cup [1/6,1/3]$ and
$\vep \in \{ \pm \}$. Here $\bar{\rho}_{\theta,\vep}$ and
$\bar{\rho}_{-\theta,\vep}$ are conjugate. By (\ref{eq:p/q})
$\bar{\rho}_{\theta,\vep}$ extends to a representation of
$\pi_{1}(M_{0})$ if and only if
$$
L_{\vep}(\theta) = 1.
$$
But this happens if and only if $\theta=\pm 1/4$ for both
$\vep=+$ and $\vep=-$, i.e.\ the set of conjugacy classes of
nonabelian representations of $\pi_{1}(M_{0})$ into $\SU(2)$ is
$\left\{ [\bar{\rho}_{1/4,-}], [\bar{\rho}_{1/4,+}] \right\}$.

By (\ref{eq:CSabelianp0}) the flat reducible
$\SU(2)$--connections on $M_{0}$ all have a Chern--Simons
invariant equal to zero, so we conclude that the image set of
the $\SU(2)$ Chern--Simons functional on $M_{0}$ has at most
three elements. By \cite[Theorem 5.6 and precedent text]{[18]}
we can therefore conclude that the set of Chern--Simons
invariants of flat $\SU(2)$--connections on $M_{0}$ is
\begin{equation}\label{eq:CSM0}
\left\{ 0 \bmod{\Z} , -\frac{1}{5} \bmod{\Z} ,
\frac{1}{5} \bmod{\Z} \right\}.
\end{equation}
In particular, the set of Chern--Simons invariants of
$[\bar{\rho}_{1/4,+}]$ and $[\bar{\rho}_{1/4,-}]$ is equal to
$\left\{-\frac{1}{5} \bmod{\Z},\frac{1}{5} \bmod{\Z} \right\}$.
We note that (\ref{eq:Jeffrey}) and (\ref{eq:CSM0}) prove the
AEC for the invariants $\tau_{r}(M_{0})$.

By (\ref{eq:Lspecialvalues}) and (\ref{eq:beta}) we have
$\beta_{\pm}(1/4)=f_{\pm}$, where $f_{-}=0$ and $f_{+}=1$. By
\refprop{prop:CS1} and (\ref{eq:beta}) we find
$$
\CS(\bar{\rho}_{1/4,\pm}) = \pm \frac{1}{6} - \frac{1}{\pi}
\int_{1/6}^{1/4} \Arg(L_{\pm}(t)) \dte t \pmod{\Z},
$$
where we use that $L_{\pm}(t) \in S^{1}$. By (\ref{eq:Lsym1})
we have
$$
\int_{1/6}^{1/4} \Arg(L_{-}(t)) \dte t =
- \int_{1/6}^{1/4} \Arg(L_{+}(t)) \dte t
$$
so we finally get
$$
\CS(\bar{\rho}_{1/4,\pm}) = \pm \left( \frac{1}{6}
- \frac{1}{\pi}
\int_{1/6}^{1/4} \Arg(L_{+}(t)) \dte t \right) \pmod{\Z}.
$$
Note that $\ima(L_{+}(t)) < 0$ on $]1/6,1/4[$ so
$$
-\frac{1}{12} = \frac{1}{\pi}(-\pi)(1/4 - 1/6) <
\frac{1}{\pi} \int_{1/6}^{1/4} \Arg(L_{+}(t)) \dte t < 0.
$$
Therefore
$$
1/6 - \frac{1}{\pi}\int_{1/6}^{1/4} \Arg(L_{+}(t)) \dte t
\in ]\frac{1}{6},\frac{1}{4}[.
$$
Since this value $\bmod{\Z}$ belongs to the set
$\{ \pm 1/5 \bmod{\Z}\}$ we conclude that it is equal to $1/5$
so
\begin{equation}\label{eq:Specialintegral1/4}
\frac{1}{\pi} \int_{1/6}^{1/4} \Arg(L_{-}(t)) \dte t =
- \frac{1}{\pi}\int_{1/6}^{1/4} \Arg(L_{+}(t)) \dte t =
\frac{1}{30}.
\end{equation}

Let us finally calculate the value of $\Psi=\Psi_{0}^{a,b}$ in
the critical points corresponding to the solutions
$(z,w)=((3+\sqrt{5})/2,1)$ and
$(z,w)=\left(-1,\frac{-3 \pm \sqrt{5}}{2}\right)$ to
(\ref{eq:stationarypzero}). Since the set of solutions to
(\ref{eq:stationarypzero}) is in one to one correspondence with
the set of conjugacy classes of nonabelian
$\SL(2,\C)$--representations of
$\pi_{1}(M_{0})$ by the proof of \refthm{thm:comparison} and
since the subset of solutions $(z,w)$ with $z \in S^{1}$ and
$w \in ]-\infty,0[$ corresponds to the set of conjugacy classes
of $\SU(2)$--representations we see that the points
$(z,w) = \left(-1,\frac{-3 \pm \sqrt{5}}{2}\right)$ correspond
to nonabelian $\SU(2)$--representations of $\pi_{1}(M_{0})$
while the points
$(z,w) = \left(\frac{3 \pm \sqrt{5}}{2},1\right)$ correspond to
nonabelian $\SL(2,\C)$--representations of $\pi_{1}(M_{0})$
which are not equivalent to $\SU(2)$--representations.

For $(z,w) = (e^{2\pi i x},e^{2\pi iy}) = ((3+\sqrt{5})/2,1)$
we find again that $y=a$ and then $\Psi(x,y) = ab$ (independent
of the choice of branch of the dilogarithm along $]1,\infty[$).

Finally, let
$(z,w) = (e^{2\pi i x},e^{2\pi iy}) =
\left(-1,\frac{-3 \pm \sqrt{5}}{2}\right)$.
The real values of the right-hand sides of
(\ref{eq:criticalpzero}) do not depend on $a$ and $b$. Taking
the imaginary values of these equations we get
\begin{eqnarray*}
0 &=& 2\pi\left(a-\rea(y)\right) +
\ima\left(\Log(1+w)- \Log(1+w^{-1})\right), \\*
0 &=& 2\pi (b-x) + \ima\left(\Log(1+w)+\Log(1+w^{-1})\right).
\end{eqnarray*}
The second of these equations is equivalent to
$$
x-b = \frac{1}{2\pi} \left( \Arg(1+w) +\Arg(1+w^{-1}) \right)
= \frac{1}{2}
$$
for both $w = \frac{-3 \pm \sqrt{5}}{2}$, and the first is
equivalent to
$$
a-\rea(y) =
\frac{1}{2\pi} \left( \Arg(1+w^{-1}) -\Arg(1+w) \right) =
\left\{
\begin{array}{cc} \frac{1}{2}, & w=\frac{-3+\sqrt{5}}{2}, \\
-\frac{1}{2}, & w=\frac{-3-\sqrt{5}}{2}.\end{array} \right.
$$
We have
$$
\Psi(x,y) = ax + by - xy + \frac{1}{4\pi^{2}}
\left( \Li_{2}(|w|) - \Li_{2}(1/|w|) \right),
$$
and since this is real by \refcor{cor:imaPsi}, we get
$$
\Psi(x,y) = ax + (b-x)\rea(y) + \frac{1}{4\pi^{2}}
\rea\left( \Li_{2}(|w|) - \Li_{2}(1/|w|) \right).
$$
Here
\begin{eqnarray*}
ax + (b-x)\rea(y) &=& (b-x)\left( \rea(y) - a \right) + ab =
ab + \frac{1}{2}\left( a - \rea(y) \right) \\*
&=& ab + \left\{
\begin{array}{cc} \frac{1}{4}, & w=\frac{-3+\sqrt{5}}{2}, \\
-\frac{1}{4}, & w=\frac{-3-\sqrt{5}}{2}.\end{array}\right.
\end{eqnarray*}
For $z \in \C \sm [0,\infty[$ we have
$$
\Li_{2}(z) + \Li_{2}(1/z) = - \frac{\pi^{2}}{6}
- \frac{1}{2}\Log^{2}(-z).
$$
This identity e.g.\ follows by differentiating the difference
of the two sides in the identity, showing that this difference
is constant on $\C \sm [0,\infty[$, and then evaluating in
$z=-1$ using that $\Li_{2}(-1)=-\pi^{2}/12$. Therefore
\begin{equation}\label{eq:dilogreciprocityreal}
\Li_{2}(t) = -\frac{\pi^{2}}{6} - \frac{1}{2}\Log^{2}(-t)
- \Li_{2}(t^{-1})
\end{equation}
for $t > 1$ for a branch of $\Li_{2}$ continuously extended
across $]1,\infty[$. Let $w_{\pm}=\frac{-3 \pm \sqrt{5}}{2}$.
We note that
\begin{equation}\label{eq:w*}
w_{+}w_{-} = 1.
\end{equation}
Moreover, by \cite[Formula (1.20) p.~7]{[21]},
\begin{equation}\label{eq:dilogvalue}
\Li_{2} \left( \frac{3-\sqrt{5}}{2} \right) =
\frac{\pi^{2}}{15}
- \Log^{2}\left( \frac{1+\sqrt{5}}{2} \right).
\end{equation}
Note also that
\begin{equation}\label{eq:w**}
\left( \frac{1+\sqrt{5}}{2} \right) ^{2} =
\frac{3+\sqrt{5}}{2}.
\end{equation}
Assume that $w=w_{-}=\frac{-3-\sqrt{5}}{2}$. Then, by
(\ref{eq:dilogreciprocityreal}),
$$
\Li_{2}(|w|) = -\frac{\pi^{2}}{6} - \frac{1}{2}\Log^{2}(-|w|)
- \Li_{2}(|w|^{-1}),
$$
where $\Log(-|w|)=\Log(|w|)+i\pi$, so
$$
\Li_{2}(|w|) = \frac{\pi^{2}}{2} - \frac{\pi^{2}}{6} -
\frac{1}{2}\Log^{2}(|w|) - \Li_{2}(|w|^{-1}) -i\pi\Log(|w|).
$$
But then
$$
\Li_{2}(|w|)-\Li_{2}(|w|^{-1}) = \frac{\pi^{2}}{2}
- \frac{\pi^{2}}{6} - \frac{1}{2}\Log^{2}(|w|)
- 2\Li_{2}(|w|^{-1})-i\pi\Log(|w|).
$$
Here $|w|^{-1} =\frac{3-\sqrt{5}}{2}$ by (\ref{eq:w*}) so by
(\ref{eq:dilogvalue}) and (\ref{eq:w**}) we get
$$
\Li_{2}(|w|)-\Li_{2}(|w|^{-1}) = \frac{\pi^{2}}{2}
- \frac{\pi^{2}}{6} - \frac{2\pi^{2}}{15} - i\pi\Log(|w|) =
\frac{\pi^{2}}{5} - i\pi\Log(|w|).
$$
By (\ref{eq:w*}) we conclude that
\begin{equation}\label{eq:dilog1/4}
\Li_{2}(|w_{\pm}|)-\Li_{2}(|w_{\pm}|^{-1}) =
\mp \frac{\pi^{2}}{5} - i\pi\Log(|w_{\pm}|).
\end{equation}
For $w=w_{\pm}$ we therefore get
$$
\Psi(x,y) = ab \pm \frac{1}{4} \mp \frac{1}{20} =
ab \pm\frac{1}{5}.
$$
Note here that we use the identities
(\ref{eq:Specialintegral1/4}) and (\ref{eq:dilog1/4}) in the
proof of \refthm{thm:main} to handle the cases
$\theta=\pm 1/4$, so we can not here refer to this theorem.

We have above shown that the set of values $\bmod{\Z}$ of
$\Psi_{0}^{a,b}$ in its critical points is identical with the
set (\ref{eq:CSM0}) of Chern--Simons invariants of flat
$\SU(2)$--connections on $M_{0}$ for arbitrary $a,b \in \Z$.
The determinant of the Hessian $H$ of $\Psi_0^{a,b}$ in a
critical point $(x,y)$ is given
by
$$
\det(H) = 1 - 2\left(z + \frac{1}{z} \right)
$$
according to (\ref{eq:detHessian3}). Thus all the critical
points are non-degenerate. It is interesting to note
that to obtain the leading order large $r$ asymptotics of
$\tau_r(M_0)$ using the saddle point method we have
to use critical points corresponding to
$\SL(2,C)$--representations of $\pi_1(M_0)$ which are not
equivalent to $\SU(2)$--representations in order to get the
part of that leading asymptotics being associated to the
reducible flat $\SU(2)$--connections on $M_0$. By
\refconj{conj:mainconj} this phenomenon should only occur for
$p/q=0$. Recall here that the case $p/q=0$ is in any case
special; only in that case the moduli space of reducible flat
$\SU(2)$--connections on $M_{p/q}$ is not discrete.

For $(z,w)=\left(-1,\frac{-3 \pm \sqrt{5}}{2}\right)$ we see
that $\det(H)=5$. Moreover, $z=e^{2\pi ix}=-1$ implies that
$x \in \frac{1}{2} +\Z$, so $\sin(\pi x) =\pm 1$. Thus we see
that the right-hand side of
the formula (\ref{eq:asympqi}) gives the part of the leading
order large $r$ asymptotics of $\tau_r(M_0)$ corresponding to
the irreducible flat $\SU(2)$--connections on $M_0$ if we let
$m_{\bar \rho}=4$ for both points in $\mM_{0}'$.

The
identity (\ref{eq:dilog1/4}) was proved by using the explicit
value of the dilogarithm in $(3-\sqrt{5})/2$, cf.\
(\ref{eq:dilogvalue}). We note that only very few explicit
values of the dilogaritm are known, see \cite[Chap~1]{[21]}.

\appendix{\,\,\,Proof of \reflem{lem:Qfunctions}}

\noindent Let us begin by showing the identitites
(\ref{eq:Qidentity1}) and (\ref{eq:Qidentity2}). Let
$\theta \in I$ and let $u = u_{\pm}(\theta)$,
$Q_{i}=Q_{i}^{\pm}(\theta)$ and $t=e^{4\pi i\theta}$ and get
\begin{eqnarray*}
Q_{1}Q_{2} &=& \left(1-(1+u)^{-1}t \right)\left(1-(1+u)t\right)
= 1 - t\left(1+u + (1+u)^{-1}\right) + t^{2} \\*
&=& -t(1+u)^{-1}
\left( (1+u)^{2} - (t+t^{-1})(u+1) + 1 \right) = t,
\end{eqnarray*}
where the last equality follows by the fact that $\phi(t,u)=0$,
where $\phi$ is given by (\ref{eq:phifunction}).

To show (\ref{eq:Qidentity2}) we observe that
\begin{eqnarray*}
\frac{Q_{1}Q_{3}}{Q_{2}} &=& \frac{Q_{3}Q_{1}^{2}}{Q_{1}Q_{2}}
= t^{-1} Q_{3}Q_{1}^{2} \\*
&=& t^{-1} (1+u) \left( 1 - t(1+u)^{-1}\right)^{2}
= t^{-1}(1+u) + t(1+u)^{-1} -2
\end{eqnarray*}
by (\ref{eq:Qidentity1}). Now $\phi(t,u)=0$ implies that
$$
t(1+u)^{-1} = t^{2} + 1 - t - t(u+1)
$$
leading to the identity
$$
\frac{Q_{1}Q_{3}}{Q_{2}} =
- 1 + t^{-1} - 2t + t^{2} + u(t^{-1}-t) = L_{\pm}(\theta).
$$

To prove the identities (\ref{eq:Qlog1}) and (\ref{eq:Qlog2})
it is necessary to examine the arguments of the functions
$Q_{i}^{\pm}(\theta)$, $i=1,2$. By \reflem{lem:lambdasymmetry}
we note that
\begin{equation}\label{eq:Q1Q2}
Q_{1}^{\pm}(\theta) = Q_{2}^{\mp}(\theta)
\end{equation}
for $\theta \in I$.

Assume first that $\theta \in [1/6,1/3]$. An elementary
calculation shows that
\begin{eqnarray}\label{eq:1+u}
&&1 + u_{+}(\theta) \in [-1, (\sqrt{5}-3)/2], \nonumber \\
&&1+u_{-}(\theta) \in [-(3+\sqrt{5})/2,-1]
\end{eqnarray}
for $\theta \in [1/6,1/3]$ and $1+u_{\pm}(\theta)=-1$ if and
only if $\theta \in \{1/6,1/3\}$. In particular,
$1+u_{\pm}(\theta)$ is negative. We therefore get
$$
\ima\left( Q_{1} \right) = |1+u|^{-1}\sin(4\pi\theta) \left\{
\begin{array}{cl} >0, & \theta \in [1/6,1/4[, \\
=0, & \theta=1/4,\\
<0, & \theta \in ]1/4,1/3]. \end{array} \right.
$$
By (\ref{eq:Q}) $\rea(Q_{1}) = 1 -(1+u)^{-1}\cos(4\pi \theta)$.
Since $1+u$ is negative we see that $\rea(Q_{1})$ have the same
sign as $\cos(4\pi\theta)-(1+u)$, where $\sign(0)=0$ as usual.
By (\ref{eq:1+u}) we conclude that
$$
\rea\left(Q_{1}^{-}(\theta)\right) > 0
$$
for all $\theta \in [1/6,1/3]$. Let us next consider
$Q_{1}^{+}$. First note that
$$
\cos(4\pi\theta) - u_{+}(\theta) - 1 = \frac{1}{2}
- \sqrt{\cos^{2}(4\pi\theta) -\cos(4\pi\theta)-\frac{3}{4}}.
$$
We therefore get that $\rea\left(Q_{1}^{+}(\theta)\right)$ has
the opposite sign as $\cos^{2}(4\pi\theta)-\cos(4\pi\theta)-1$.
By the assumption on $\theta$ we have
$\cos(4\pi\theta) \in [-1,-1/2]$, and we therefore get
$$
\rea\left(Q_{1}^{+}(\theta)\right) \left\{ \begin{array}{cl}
>0, & \theta \in [1/6,\theta_{0}[ \cup ]1/2-\theta_{0},1/3],\\
=0, & \theta \in \{\theta_{0},1/2-\theta_{0} \},\\
<0, & \theta \in ]\theta_{0},1/2-\theta_{0}[,\end{array}\right.
$$
where $\theta_{0} \in ]1/6,1/4[$ is the unique element such
that $\cos(4\pi \theta_{0}) = (1-\sqrt{5})/2 \in ]-1,-1/2[$ is
the negative solution to $t^{2}-t-1=0$.

Let $\psi_{i}^{\pm}(\theta) \in ]-\pi,\pi]$ be the principal
argument of $Q_{i}^{\pm}(\theta)$. Then the above analysis
shows, also using (\ref{eq:Q1Q2}), that
\begin{equation}\label{eq:arguments1}
\psi_{2}^{-}(\theta) = \psi_{1}^{+}(\theta) \left\{
\begin{array}{ll}
\in ]0,\frac{\pi}{2}[,&\theta \in [\frac{1}{6},\theta_{0}[, \\
= \frac{\pi}{2},& \theta = \theta_{0}, \\
\in ]\frac{\pi}{2},\pi[,&
\theta \in ]\theta_{0},\frac{1}{4}[, \\
= \pi,& \theta = \frac{1}{4}, \\
\in ]-\pi,-\frac{\pi}{2} [,&
\theta \in ]\frac{1}{4},\frac{1}{2} -\theta_{0}[, \\
= -\frac{\pi}{2},& \theta = \frac{1}{2} -\theta_{0}, \\
\in ]-\frac{\pi}{2},0[,&
\theta \in ]\frac{1}{2}-\theta_{0},\frac{1}{3}],
\end{array} \right.
\end{equation}
and
\begin{equation}\label{eq:arguments2}
\psi_{2}^{+}(\theta) = \psi_{1}^{-}(\theta) \left\{
\begin{array}{ll}
\in ]0,\frac{\pi}{2}[,&
\theta \in [\frac{1}{6},\frac{1}{4}[, \\
= 0,& \theta = \frac{1}{4}, \\
\in ]-\frac{\pi}{2},0[,&
\theta \in ]\frac{1}{4},\frac{1}{3}], \end{array} \right.
\end{equation}
so
$$
\psi_{1}^{\pm}(\theta) + \psi_{2}^{\pm}(\theta) \in \left\{
\begin{array}{ll}
]0,\pi[,& \theta \in [\frac{1}{6},\theta_{0}[, \\
]\frac{\pi}{2},\frac{3\pi}{2}[,&
\theta \in ]\theta_{0},\frac{1}{4}[, \\
]\frac{\pi}{2},\pi],& \theta = \frac{1}{4}, \\
]-\frac{3\pi}{2},-\frac{\pi}{2} [,&
\theta \in ]\frac{1}{4},\frac{1}{2} -\theta_{0}[, \\
]-\pi,0[,&
\theta \in [\frac{1}{2}-\theta_{0},\frac{1}{3}].
\end{array} \right.
$$
By (\ref{eq:Qidentity1}) we have
$$
\psi_{1}^{\pm}(\theta) + \psi_{2}^{\pm}(\theta) \in
4\pi \theta + 2\pi \Z,
$$
so we conclude that
$$
\psi_{1}^{\pm}(\theta) + \psi_{2}^{\pm}(\theta) = 4\pi \theta
$$
for all $\theta \in [1/6,1/3]$ proving (\ref{eq:Qlog1}) for
these $\theta$.

Next assume that $\theta \in [-1/3,-1/6]$. First observe that
\begin{equation}\label{eq:Qsym1}
Q_{i}^{\pm}(-\theta) = \overline{Q_{i}^{\pm}(\theta)},
\hspace{.2in} \theta \in I \sm \{ \pm 1/4 \},
\end{equation}
and
\begin{equation}\label{eq:Qsym2}
Q_{i}^{\pm}(-1/4) = Q_{i}^{\pm}(1/4)
\end{equation}
for $i=1,2$.
By (\ref{eq:Qsym2}) we immediately get that (\ref{eq:Qlog1})
holds for $\theta=-1/4$. For
$\theta \in [-1/3,-1/6] \sm \{-1/4\}$, (\ref{eq:Qlog1}) follows
by (\ref{eq:Qsym1}) and the fact that
$\Log(\bar{p}) = \overline{\Log(p)}$ for
$p \in \C \sm ]-\infty,0]$.

Note that (\ref{eq:Qlog2}) is true if we choose
$e_{\pm}(\theta) \in \Z$ such that
$$
\psi_{1}^{\pm}(\theta) + \pi - \psi_{2}^{\pm}(\theta)
- e_{\pm}(\theta)2\pi \in ]-\pi,\pi].
$$
By (\ref{eq:arguments1}) and (\ref{eq:arguments2}) we have that
$\psi_{1}^{\pm}(\theta) - \psi_{2}^{\pm}(\theta)
\in ]-\pi,\pi]$,
and we conclude that we have to put $e_{\pm}(\theta)=0$ if and
only if $\psi_{1}^{\pm}(\theta) \leq \psi_{2}^{\pm}(\theta)$
and $e_{\pm}(\theta)=1$ elsewhere. By (\ref{eq:arguments1}) and
(\ref{eq:arguments2}) we conclude that $e_{-}(1/4) = 0$ and
$e_{+}(1/4)=1$.

Assume that $\theta \in [1/6,1/3] \sm \{1/4\}$. Then
$\psi_{1}^{\pm}(\theta)$ and $\psi_{2}^{\pm}(\theta)$ both
belong to either $]-\pi,0[$ or to $]0,\pi[$. We use this fact
together with the fact that $\cot :]m\pi,(m+1)\pi[ \to \R$ is
strictly decreasing for any $m \in \Z$. In fact,
$\cot(\psi_{i}) = \frac{\rea(Q_{i})}{\ima(Q_{i})}$ so
\begin{eqnarray*}
\cot(\psi_{1}) &=& \cot(4\pi\theta)
+ \frac{|1+u|}{\sin(4\pi\theta)}, \\
\cot(\psi_{2}) &=& \cot(4\pi\theta)
+ \frac{|1+u|^{-1}}{\sin(4\pi\theta)}.
\end{eqnarray*}
By this we find that
$$
\sign\left(
\psi_{1}^{\pm}(\theta) - \psi_{2}^{\pm}(\theta) \right)
= \sign\left( \sin(4\pi\theta)
\left(1 - \left(1+u_{\pm}(\theta)\right)^{2} \right)\right).
$$
Since $(1+u_{+}(\theta))^{2} \leq 1$ and
$(1+u_{-}(\theta))^{2} \geq 1$ with equalities if and only if
$\theta \in \{1/6,1/3\}$, we get
$\psi_{1}^{\pm}(\theta) = \psi_{2}^{\pm}(\theta)$ for
$\theta \in \{1/6,1/3\}$ and
$$
\psi_{1}^{\vep}(\theta) < \psi_{2}^{\vep}(\theta)
$$
for $\theta \in ]1/6,1/4[$ and $\vep=-$ or for
$\theta \in ]1/4,1/3[$ and $\vep=+$. Moreover,
$$
\psi_{1}^{\vep}(\theta) > \psi_{2}^{\vep}(\theta)
$$
for $\theta \in ]1/6,1/4[$ and $\vep=+$ or for
$\theta \in ]1/4,1/3[$ and $\vep=-$. We therefore get that
(\ref{eq:Qlog2}) is true if we put $e_{\pm}(\theta)=0$ for
$\theta \in \{1/6,1/3\}$ and
$$
e_{+}(\theta) = \left\{ \begin{array}{cl}
1,& \theta \in ]1/6,1/4[,\\
0,& \theta \in ]1/4,1/3[, \end{array} \right.
$$
and let $e_{-}(\theta)=1-e_{+}(\theta)$ for
$\theta \in ]1/6,1/3[ \sm \{1/4\}$.

Let us finally consider the case $\theta \in [-1/3,-1/6]$.
First observe that
\begin{equation}\label{eq:Qsym3}
Q_{3}^{\pm}(-\theta) = Q_{3}^{\pm}(\theta)
\end{equation}
for all $\theta \in I$. By this and (\ref{eq:Qsym2}) we
conclude that $e_{\pm}(-1/4) = e_{\pm}(1/4)$. For
$\theta \in [-1/3,-1/6] \sm \{-1/4\}$ we get by
(\ref{eq:Qsym1}) and (\ref{eq:Qsym3}) that
\begin{eqnarray*}
&&\Log\left(Q_{1}^{\pm}(\theta)\right)
+ \Log\left(Q_{3}^{\pm}(\theta)\right)
- \Log\left(Q_{2}^{\pm}(\theta)\right)
= \Log\left(Q_{3}^{\pm}(\theta)\right)
- \overline{\Log\left(Q_{3}^{\pm}(\theta)\right)} \\*
&& \hspace{.5in}
+ \overline{\Log\left(Q_{1}^{\pm}(-\theta)\right)
+ \Log\left(Q_{3}^{\pm}(-\theta)\right)
- \Log\left(Q_{2}^{\pm}(-\theta)\right)} \\
&& \hspace{.3in} =
2i\ima\left(\Log\left(Q_{3}^{\pm}(\theta)\right)\right)
+ \overline{\left( \Log\left(
\frac{Q_{1}^{\pm}(-\theta)Q_{3}^{\pm}(-\theta)}
     {Q_{2}^{\pm}(-\theta)} \right)
+ e_{\pm}(-\theta)2\pi i \right)}.
\end{eqnarray*}
By using that $Q_{3}$ is negative and that
$L_{\pm}(\theta) \in \C \sm ]-\infty,0]$ for
$\theta \in ]1/6,1/3[$ we get for
$\theta \in ]-1/3,-1/6[ \sm \{-1/4\}$ that
\begin{eqnarray*}
\Log\left(Q_{1}^{\pm}(\theta)\right)
+ \Log\left(Q_{3}^{\pm}(\theta)\right)
- \Log\left(Q_{2}^{\pm}(\theta)\right)
&=& \Log\left(
\frac{Q_{1}^{\pm}(\theta)Q_{3}^{\pm}(\theta)}
     {Q_{2}^{\pm}(\theta)} \right) \\
&& + \left(1-e_{\pm}(-\theta)\right)2\pi i,
\end{eqnarray*}
so we conclude that $e_{\pm}(\theta)=1-e_{\pm}(-\theta)$ for
these $\theta$. Finally we get for $\theta \in \{-1/3,-1/6\}$
that $Q_{1}^{\pm}(\theta)=Q_{2}^{\pm}(\theta)$ so
(\ref{eq:Qlog2}) is satisfied for $e_{\pm}(\theta)=0$ for these
$\theta$.

\end{document}